\documentclass[10pt]{article}  
\usepackage{amsmath}
\usepackage{amssymb}
\usepackage{theorem}

\theoremheaderfont{\scshape}
\newtheorem{thm}{Theorem}[section]

\theorembodyfont{\normalfont}
\newtheorem{rem}{Remark}[section]

\newcommand{\nit}{\mathbb{N}}
\newcommand{\rit}{\mathbb{R}}
\newcommand{\sit}{\mathbb{S}}

\newcommand{\eps}{\varepsilon}
\newcommand{\avec}{\mathbf{a}}
\newcommand{\aavec}{\mathbf{A}}
\newcommand{\bvec}{\mathbf{b}}
\newcommand{\bbvec}{\mathbf{B}}
\newcommand{\eevec}{\mathbf{E}}
\newcommand{\xvec}{\mathbf{x}}
\newcommand{\xze}{\xvec^0}
\newcommand{\xun}{\xvec^1}

\newcommand{\xk}{\xvec^k}
\newcommand{\yvec}{\mathbf{y}}
\newcommand{\yze}{\yvec^0}
\newcommand{\yun}{\yvec^1}
\newcommand{\yde}{\yvec^2}

\newcommand{\yi}{\yvec^i}

\newcommand{\yj}{\yvec^j}
\newcommand{\yjmun}{\yvec^{j-1}}
\newcommand{\ykmun}{\yvec^{k-1}}

\newcommand{\yk}{\yvec^k}

\newcommand{\tmsoe}{\frac{t-s}{\eps}}

\newcommand{\xxvec}{\mathbf{X}}
\newcommand{\xxep}{\xxvec_\eps}
\newcommand{\xxze}{{\xxvec^0}}
\newcommand{\xxzeep}{{\xxvec_\eps}}
\newcommand{\xxun}{{\xxvec^1}}

\newcommand{\xxde}{{\xxvec^2}}

\newcommand{\xxi}{{\xxvec^i}}

\newcommand{\xxk}{{\xxvec^k}}
\newcommand{\xxkmun}{{\xxvec^{k-1}}}

\newcommand{\xxkep}{{\xxvec^k_\eps}}
\newcommand{\xxkmunep}{{\xxvec^{k-1}_\eps}}

\newcommand{\yyvec}{\mathbf{Y}}
\newcommand{\yyze}{{\yyvec^0}}
\newcommand{\yyzeep}{{\yyvec_\eps}}
\newcommand{\yyun}{{\yyvec^1}}

\newcommand{\yyde}{{\yyvec^2}}

\newcommand{\yyj}{{\yyvec^j}}
\newcommand{\yyjmun}{{\yyvec^{j-1}}}
\newcommand{\yyk}{{\yyvec^k}}
\newcommand{\yykmun}{{\yyvec^{k-1}}}

\newcommand{\yykep}{{\yyvec^k_\eps}}
\newcommand{\yykmunep}{{\yyvec^{k-1}_\eps}}

\newcommand{\zvec}{\mathbf{z}}
\newcommand{\zzvec}{\mathbf{Z}}
\newcommand{\uvec}{\mathbf{u}}
\newcommand{\uuvec}{\mathbf{U}}
\newcommand{\wvec}{\mathbf{w}}

\newcommand{\vvec}{\mathbf{v}}
\newcommand{\vvvec}{\mathbf{V}}
\newcommand{\evec}{\mathbf{E}}

\newcommand{\caln}{{\cal N}}
\newcommand{\calm}{{\cal M}}

\newcommand{\calr}{{\cal R}}

\newcommand{\basevecun}{\mathbf{e}_1}
\newcommand{\basevecde}{\mathbf{e}_2}
\newcommand{\basevectr}{\mathbf{e}_3}
\newcommand{\cxun}{{x_1}}
\newcommand{\cxde}{{x_2}}
\newcommand{\cxtr}{{x_3}}

\newcommand{\cyun}{{y_1}}
\newcommand{\cyde}{{y_2}}

\newcommand{\cuun}{{u_1}}
\newcommand{\cude}{{u_2}}

\newcommand{\czun}{{z_1}}
\newcommand{\czde}{{z_2}}

\newcommand{\Nabla}{\nabla\hspace{-2pt}}
\newcommand{\ov}{\overline}
\newcommand{\fracp}[2]{\frac{\partial #1}{\partial #2}}

\newcommand{\ds}{\displaystyle}
\numberwithin{equation}{section}
\newenvironment{proof}{\noindent{\it Proof. }}{\hfill\rule{2mm}{2mm}\vskip3mm \par} 
\author{
E. Fr\'enod\thanks{LMAM, Universit\'e de Bretagne Sud, Centre Yves Coppens, Campus de Tohannic,
 F-56000, Vannes}
} 
\title{Application of the averaging method to the gyrokinetic plasma
}
\begin{document}
\maketitle

{\small
{\bf Acknowledgements:} I would like to thank E. Sonnendr\"ucker for the stimulating
discussions we had together and from which the problem treated in this paper
emerged. 

~

{\bf Abstract:} we show that the solution to an oscillatory-singularly 
perturbed ordinary differential equation may be asymptotically expanded 
into a sum of oscillating terms.
Each of those terms writes as an oscillating operator acting on the 
solution to a non oscillating ordinary differential equation with 
an oscillating correction added to it.
\newline
The expression of the non oscillating ordinary differential equations 
are defined by a recurrence relation.
\newline
We then apply this result to problems where charged particles are submitted to large magnetic field.
}
\section{Introduction and results}
{\bf Purpose} \newline
The goal of this paper is to manage and justify the asymptotic two scale expansion as 
$\eps\rightarrow 0$ of 
the solution $\xxep(t)=\xxep(t;\xvec,s) $ to the following singularly perturbed 
dynamical system: 
\begin{equation}\label{sdsp}
\frac{d\xxep}{dt} = \avec (t,\tmsoe,\xxep) +\frac1\eps \bvec (t,\xxep),
\;\; \xxep(s;\xvec,s)= \xvec,
\end{equation}
and to apply this to models of charge particles submitted to a strong magnetic field
which is possibly non uniform.

~

More precisely, $C_b^{m}$ standing for the space of functions that have continuous and
bounded derivatives until the order $m$, we assume that 
\begin{gather} \ds 
\avec(\cdot,\cdot,\cdot) \in (C_b^{m+1}(\rit\times\rit\times\rit^d) )^d,
\text{ }
\theta\mapsto\avec(t,\theta,\xvec)
\text{ is $2\pi-$periodic for every $t\in\rit$ and $\xvec\in\rit^d$},
\label{regab}
\\
\bvec(\cdot,\cdot) \in (C_b^{m+2}(\rit\times\rit^d))^d,
\label{divb}
\end{gather}
for some $m\geq 0.$ In particular, it implies that $|\Nabla_x\cdot\bvec| \leq C,$
uniformly on $\rit\times\rit^d$.
We also suppose that the solution $\zzvec(t,\theta;\zvec)$ to
\begin{equation}\label{sdz}
\frac{\partial\zzvec }{\partial\theta} = \bvec(t,\zzvec), \;\;
\zzvec(t,0;\zvec) = \zvec,
\end{equation}
is known and $2\pi-$periodic in $\theta$ for every $t\in\rit$
and $\zvec\in\rit^d$.

Under those assumptions, we prove that $\xxep(\cdot;\xvec,s): \rit\rightarrow\rit^d $ 
admits the following expansion:
\begin{equation}\label{da}
\xxep(t;\xvec,s)=\xxze(t,\tmsoe;\xvec,s)+
\eps\xxun (t,\tmsoe;\xvec,s)+ \eps^2\xxde (t,\tmsoe;\xvec,s)+\dots,
\end{equation} 
as $\eps\rightarrow 0$,
where the $\xxi(t,\theta;\xvec,s)$ are $2\pi-$periodic in $\theta.$
The initial condition in (\ref{sdsp}) is assigned to $\xxze$ by:
\begin{equation}\label{aci}
\xxze(s,0;\xvec,s)= \xvec \text{ and } \xxi(s,0;\xvec,s)=0 
\text{ for } i\geq 1.
\end{equation}
{\bf Motivations} \newline
The target we have in mind while writing this paper is the conception of 
pic-methods for plasmas submitted to large magnetic field. Direct applications of this
are the simulation of magnetic confinement fusion or isotope resonant separation 
experiments.

~

A plasma submitted to a strong magnetic field could be modelled by several 
Vlasov equations, according to experimental conditions. 
For instance if we are interested 
in describing the global behaviour of a plasma in a magnetic confinement 
fusion experiment, the following equation
\begin{equation}\label{vlglob}
\begin{array}{l}
\ds\fracp{f^{\eps}}{t}+\vvec\cdot\Nabla_xf^{\eps}+
\Big((\evec(t,\xvec)+\frac{\caln}{\eps})+
\vvec\times(\bbvec(t,\xvec)+\frac{\calm}{\eps})\Big)
\cdot\Nabla_vf^{\eps} =0,
~~~
\ds f^\eps_{|t=s}=f_0,
\end{array}
\end{equation}
for a small parameter $\eps$ 
would be a relevant model equation. In equation (\ref{vlglob}), $f^\eps(t,\xvec,\vvec)$ 
is at time $t$, the
particle density of the plasma in position $\xvec$ and with velocity $\vvec$.
The vector $\calm \in \sit^2$ gives the direction of the strong magnetic field.
This direction may depend on $t$ and $\xvec$. the vector
$\caln\in\sit^2$ with $\caln \perp \calm$ is the direction of a 
strong electric field ($\caln$ may also be $0$).
The vector fields $\evec$ and $\bbvec$ are electric and magnetic fields that
may contain a self consistent part.
We shall call the scaling leading to this equation the ``Guiding Centre Regime''.

If now we are interested in understanding what happens close to the plasma boundary in a
tokamak, we would prefer to use the so called ``Finite Larmor Radius Regime''.
In other words, we would consider an equation of the type:
\begin{equation}\label{vlRLF}
\begin{array}{l}\ds
\frac{\partial f^\eps}{\partial t} +\vvec_\parallel\cdot\Nabla_x f^\eps 
+\frac{\vvec_\perp}{\eps}\cdot\Nabla_x f^\eps 
+\Big(\evec(t,\xvec)+\vvec\times \frac{\calm}{\eps}\Big)\cdot\Nabla_v f^\eps=0,
~~~
f^\eps_{|t=s}= f_0,
\end{array}
\end{equation} 
where the meaning of $f^\eps(t,\xvec,\vvec)$ is the same as above and where
$\vvec_\parallel=(\vvec\cdot\calm)\calm$,
$\vvec_\perp=\vvec-\vvec_\parallel$.

Lastly, for isotope resonant separation we would be interested in regarding
\begin{equation}\label{vlCR}
\begin{array}{l}
\ds \frac{\partial f^\eps}{\partial t}+ \vvec\cdot\Nabla_x f^\varepsilon
+ \Big(\evec(t,\tmsoe,\xvec)+
\vvec\times (\bbvec(t,\tmsoe,\xvec)+\frac{\cal M}{\varepsilon})\Big)
\cdot\Nabla_v f^\varepsilon =0,
~~~
\ds f^\varepsilon_{|t=s}=f_0,
\end{array}
\end{equation}
where the  electric and magnetic fields  $\evec(t,\tmsoe,\xvec)$ and  
$\bbvec(t,\tmsoe,\xvec)$ 
oscillate with the same frequency $1/(2\pi\eps)$ as the cyclotron frequency
of the plasma particles.

We refer to Fr\'enod and Sonnendr\"ucker 
\cite{frenod/sonnendrucker:PROC1998,frenod/sonnendrucker:1997,
frenod/sonnendrucker:CRAS99,frenod/sonnendrucker:1998, 
frenod/sonnendrucker:1999},
Fr\'enod, Raviart and Sonnendr\"ucker \cite{FRS:1999} and
Fr\'enod and Watbled \cite{FW}
where those models are explained and asymptotically analysed.
In particular, they give indication on how to get the self consistent
part of the electric field in some examples and under suitable assumptions.
Complementary works concerning the asymptotic behaviour of those kind of equations
are also led in Golse and Saint Raymond \cite{GSRcras,GSR1},
Saint Raymond \cite{saintraymond:2000},
Brenier \cite{brenier:2000}, Grenier \cite{grenier:1997,grenier2:1997},
Jabin \cite{jabin:2000}, 
Schochet \cite{schochet:1994}, 
Joly, M\'etivier and Rauch \cite{JMRDDD}.
We also refer to mathematical or physical works where similar
methods are used: \cite{JMR93, JMR96, serre91, frenod:1994, frenod/hamdache:1996, 
littlejohn:1981, lee:1983, dubin/etal:1983, cohen:1985,lochak/meunier:1988,grad:1971}.

~

It is relatively easy to see that equations (\ref{vlglob}), 
(\ref{vlRLF}) and (\ref{vlCR})
enter the framework of a singularly perturbed convection equation reading:
\begin{equation} \label{spce}
\begin{array}{l}
\ds \fracp{u_\eps}{t} +\avec\cdot\nabla u_\eps 
+\frac{1}{\eps}\bvec\cdot\nabla u_\eps =0, 
~~~~
\ds{u_\eps}_{|t=s} = u_0,
\end{array}
\end{equation}
for $\xvec\in\rit^d$ and $t>s$. At least formally, the solution $u_\eps(t,\xvec)$ 
is given by
\begin{equation}\label{ueps}
u_\eps(t,\xvec)= u_0(\xxep(s;\xvec,t)).
\end{equation}

Numerical methods to solve equation (\ref{vlglob}), 
(\ref{vlRLF}) or (\ref{vlCR}),
or similar ones, will be applications of generic methods solving (\ref{spce}).
Since $u_\eps$ contains high frequency oscillations, using in state any well 
known numerical method to solve (\ref{spce}) would compel to use a very small time step.

Yet, in order to relax this constraint, we may follow, at least, two strategies.
The first one could be based on the work presented in   
Fr\'enod, Raviart and Sonnendr\"ucker \cite{FRS:1999}
where we proved that $u_\eps$ writes
$u_\eps(t,\xvec) =\sum_{i \geq 0} \eps^i U^i(t,\tmsoe,\xvec),$
and where we determined the equations satisfied by every terms $U^i$. Since those
equations are independent of $\eps$, the computation of the first terms  $U^i$ could
be led using a standard numerical method. Then an approximation of $u_\eps$ would
be given by:
$u_\eps(t,\xvec) \simeq \sum_{i=0}^{k} \eps^i U^i(t,\tmsoe,\xvec),$
for some $k\in \nit$.

The second strategy consists in using an approximation 
$\xxep(t;\xvec,s) \simeq \sum_{i=0}^{k} \eps^i\xxi(t,\tmsoe;\xvec,s),$
of (\ref{da}), for some $k\in \nit$, after computing $\xxze, \xxun,\dots, \xxk$. 
Then the expression (\ref{ueps})
yields an approximation of $u_\eps$:
$u_\eps(t,\xvec)\simeq u_0 \big(\sum_{i=0}^{k} \eps^i \xxi(s,-\tmsoe;\xvec,t)\big).$
This motivates the present work.

~
 
We give now some references concerning perturbed ordinary differential equations.
First, we think to the Lindstedt-Poincar\'e method explained 
in Poincar\'e \cite{poincareMNMC} where 
the steady state periodic solutions to a perturbed  second order 
ordinary differential equation is studied.
Then the Krylov-Bogoliubov-Mitropolsky method, see \cite{BogoMitro} and \cite{KryBo}, 
allows to describe the transitory
behaviour of the solution to a perturbed ordinary differential equation to 
a periodic solution.

We also cite the works of Verhulst \cite{Ver} and 
Sanders and Verhulst \cite{SanVer} (see also Mickens \cite{Mic})
where the Method of Averaging  is developed to treat 
perturbed ordinary differential equations and adiabatic invariants 
in Hamiltonian systems.

Lastly, we mention the work initiated by Tikhonov 
\cite{tikhonov:1948,tikhonov:1950,tikhonov:1952}
and developed by Vasilieva and others (see the review paper of  Vasilieva 
\cite{vasilieva:1994} and the references in it).

~

This work is also related to 
homogenisation methods based on weak$-*$ convergence
and two scale convergence. The most important ideas about those tools
may be found in Tartar \cite{tartar:1977},  
Bensoussan, Lions and Papanicolaou \cite{bensoussan/lions/papanicolaou:1978},
Sanchez-Palencia \cite{sanchez-palencia:1978}, N'Guetseng \cite{nguetseng:1989},
Allaire \cite{allaire:1992} and 
Fr\'enod, Raviart and Sonnendr\"ucker \cite{FRS:1999}.

~

\noindent{\bf Theorems} \newline
We first justify the expansion (\ref{da}) until order 0.
\begin{thm}\label{thmX0}
Setting
\begin{equation}\label{a0tilde}
\tilde{\avec}^0(t,\yze) = \frac{1}{2\pi}\int_0^{2\pi}
\tilde \alpha^0(t,\theta,\yze)\,d\theta,
\end{equation}
with
\begin{equation}\label{al0}
\tilde \alpha^0(t,\theta,\yze) = \big\{\Nabla_z \zzvec(t,\theta;\yze)\big\}^{-1} 
\big\{\avec(t,\theta,\zzvec(t,\theta;\yze))-\fracp{\zzvec}{t}(t,\theta;\yze)\big\},
\end{equation}
under assumptions (\ref{regab}), (\ref{divb}) with $m=0$ and (\ref{sdz}), for any  
$\xvec\in\rit^d$, $s\in\rit$, $T\in \rit$ and any $\eps>0$,
the solution $\xxep(\cdot;\xvec,s)$ of (\ref{sdsp}) exists on $[s,s+T]$, is unique and the 
sequence $(\xxep(\cdot;\xvec,s))$ satisfies:
\begin{equation}\label{app0}
\lim_{\eps\rightarrow 0} \sup_{t\in[s,s+T]}
\big|\xxep(t;\xvec,s)-\xxze(t,\tmsoe;\xvec,s) \big|  = 0,
\end{equation}
$|\cdot| $ standing for the Euclidean norm on $\rit^d$,
where $\xxze$ satisfies
\begin{equation}\label{X0-Y0}
\xxze(t,\theta;\xvec,s) = \zzvec(t,\theta,\yyze(t;\xvec,s)),
\end{equation}
and where $\yyze$ is the solution to 
\begin{equation}\label{eqY0}
\frac{d\yyze}{dt} = \tilde{\avec}^0(t,\yyze) ,\;\; \yyze(s;\xvec,s)=\xvec.
\end{equation}
\end{thm}
In the Theorem above, $\Nabla_z \zzvec(t,\theta;\zvec)$ stands for the Jacobian matrix of
$\zvec\mapsto\zzvec(t,\theta;\zvec)$. 

~

In order to set this Theorem and to justify the expansion (\ref{da}) 
for higher orders, using a Van der Pol transformation, we define $\yyzeep$ 
being such that
\begin{equation}\label{defy0eps}
\xxzeep(t;\xvec,s)= \zzvec(t,\tmsoe;\yyzeep(t,\xvec,s))=[\zzvec(\yyzeep)]_\eps,
\end{equation}
where for any function $f$ we write
$[f(\yyzeep,\dots,\yykep)]_\eps$ for
$f(t,\tmsoe,\yyzeep(t;\xvec,s),\dots,\yykep(t;\xvec,s))$
and $[f(\yyze,\dots,\yyk)]_\eps$ for
$f(t,\tmsoe,\yyze(t;\xvec,s),\dots,\yyk(t;\xvec,s))$.
It is an easy game to show that
\begin{equation}\label{eqy0ep}
\frac{d\yyzeep}{dt} = \tilde \alpha^0(t,\tmsoe,\yyzeep)= [\tilde \alpha^0(\yyzeep)]_\eps; 
\;\; \yyzeep(s;\xvec,s)= \xvec.
\end{equation}
Indeed, derivating (\ref{defy0eps}) and using (\ref{eqy0ep}) we get
\begin{equation}\label{eqx0ep}
\begin{array}{l}
\ds \frac{d\xxzeep}{dt}= \frac{1}{\eps} \big[\frac{\partial\zzvec }{\partial\theta}(\yyzeep)\big]_\eps
+\big[\frac{\partial\zzvec }{\partial t}(\yyzeep)\big]_\eps +
\big\{\Nabla_z\zzvec(\yyzeep)\}_\eps\{\tilde \alpha^0(\yyzeep)\big\}_\eps
=[\avec(\xxzeep)]_\eps+\frac{1}{\eps}\bvec(\xxzeep),
\end{array}
\end{equation}
and we have the initial condition $\xxzeep(s;\xvec,s)=\xvec.$
\newline
We also use the following notation. For a vector field $\zzvec$, 
The $i-$th component of $\{\Nabla_x^{\,k}\zzvec \}\{\xze,\dots,\xk\}$,
for $i=1,\dots,d$, is given by:
\begin{equation}\label{defnab}
\begin{array}{c}
\ds \big(\{\Nabla_x^{\,k}\zzvec \}\{\xze,\xun,\dots,\xk\}\big)_i =
\sum_{l_1,\dots, l_k=1}^{d}\frac{\partial^k\zzvec_i}{\partial x_{l_1}
\dots\partial x_{l_k}}
\xze_{l_1}\dots\xk_{l_k}.
\end{array}
\end{equation}
In order to simplify we shall sometimes denote
$
\{\Nabla_x^{\,k}\zzvec \}\{\xze,\xze,\dots,\xze\} \text{ by }
\{\Nabla_x^{\,k}\zzvec \}\{\xze\}^k
$.

Now for $k\geq 0$ we recursively define
\begin{equation}\label{defAitild}
\tilde{\aavec}^k(t,\theta,\yze,\dots,\yk)= 
\frac{1}{\theta}\int_0^{\theta}\tilde \alpha^k(t,\sigma,\yze,\dots,\yk)\,d\sigma -
\tilde{\avec}^k(t,\yze,\dots,\yk),
\end{equation}
where $\tilde\alpha^0$ and $\tilde{\avec}^0$ are given by (\ref{al0}) and (\ref{a0tilde})
and where for $k\geq 1$ we have:
\begin{equation}\label{ali}
\begin{array}{l}
\ds \tilde \alpha^k(t,\theta,\yze,\dots,\yk)=
\big\{\Nabla_{y^0} \tilde \alpha^0(t,\theta;\yze)\big\}
\{\yk+\theta\tilde{\aavec}^{k-1}(t,\theta,\yze,\dots,\ykmun)\}+
\\ \ds ~
\frac12 \big\{\Nabla_{y^0}^{\,2} \tilde \alpha^0(t,\theta;\yze)\big\}
\Big(
\sum_{j=1}^{k-1}\{\yj+ \theta\tilde{\aavec}^{j-1}(t,\theta,\yze,\dots,\yjmun)\,,\,
                  \yvec^{k-j}+ \theta\tilde{\aavec}^{k-j-1}(t,\theta,\yze,\dots,\yvec^{k-j-1})\}
\Big)+
\\ \ds ~
\dots +\frac1{k!}  \big\{\Nabla_{y^0}^{\,k} \tilde \alpha^0(t,\theta;\yze)\big\}
\{\yun+\theta\tilde{\aavec}^{0}(t,\theta,\yze)\}^k-
\\ \ds
~ \hfill \theta \Big(
\sum_{j=0}^{k-1} \{\Nabla_{y^j}\tilde{\aavec}^{k-1}(t,\theta,\yze,\dots,\ykmun)\}
\{\tilde{\avec}^j(t,\yze,\dots,\yj)\}+
\frac{\partial\tilde{\aavec}^{k-1}}{\partial t}(t,\theta,\yze,\dots,\ykmun)
\Big),
\end{array}
\end{equation}
and \vspace{-5pt}
\begin{equation}\label{aitilde}
\tilde{\avec}^k(t,\yze,\dots,\yk) = \frac{1}{2\pi}\int_0^{2\pi}
\tilde\alpha^k(t,\theta,\yze,\dots,\yk)\,d\theta.
\end{equation}

~

With those notations, we can state the asymptotic expansion of $\yyzeep$,
\begin{equation}\label{expYeps}
\yyzeep = \yyze + \eps (\yyun + [\theta \tilde \aavec^{0}]_\eps)
+ \eps^2 (\yyde + [\theta \tilde \aavec^{1}]_\eps)+ \dots \, .
\end{equation}
In other words, defining for $k\geq1$,
\begin{equation}\label{defykeps}
\begin{array}{cl}
\ds\yykep &\ds= \ds\frac{1}{\eps^k}
\big(
\yyzeep-\yyze -\eps (\yyun + [\theta \tilde \aavec^{0}]_\eps)
-\dots
-\eps^{k-1}(\yykmun+ [\theta \tilde \aavec^{k-2}]_\eps)
\big)-[\theta \tilde{\aavec}^{k-1}]_\eps,
\\\ds
      &\ds=\frac{1}{\eps}\big(\yykmunep - \yykmun \big)
-[\theta \tilde{\aavec}^{k-1}]_\eps,
\end{array}
\end{equation}
we have the following Theorem.
\begin{thm}\label{thmappk}
Under assumptions (\ref{regab}), (\ref{divb}) and (\ref{sdz}) for any $\xvec\in\rit^d$, 
$s\in\rit$ and $T\in \rit$, the sequences
$(\yykep(\cdot;\xvec,s))$, for $k=0,\dots,m$, are bounded in $L^\infty([s,s+T])$ and we have
\begin{equation}\label{appk}
\begin{array}{l}\ds
\lim_{\eps\rightarrow 0} \sup_{t\in[s,s+T]}
\big|\yyzeep(t;\xvec,s)-\yyze(t;\xvec,s) \big| = 0,
\text{ and }
\lim_{\eps\rightarrow 0} \sup_{t\in[s,s+T]}
\big|\yykep(t;\xvec,s)-\yyk(t;\xvec,s) \big| = 0,
\end{array}
\end{equation}
for $k\geq 1$, where $\yyze$ is solution to (\ref{eqY0}) and  where 
$\yyk$ is the solution to \vspace{-5pt}
\begin{equation}\label{eqYi}
~~~~\frac{d\yyk}{dt} = \tilde{\avec}^k(t,\yyze,\dots,\yyk) ,\;\; \yyk(s;\xvec,s)=0.
\end{equation}
\end{thm}
As a consequence of this Theorem and of (\ref{defy0eps}), defining
$(\xxkep(t;\xvec,s))$ by:
\begin{equation}\label{defxkeps}
\ds \xxkep=\ds\frac{1}{\eps^k}\big(\xxep- 
[\xxze]_\eps-\dots-
\eps^{k-1}[\xxkmun]_\eps\big)\vspace{2pt}
=\ds\frac{1}{\eps}\big(\xxkmunep - [\xxkmun]_\eps\big),
\end{equation}
for $k\geq 1$, we have the rigorous justification of asymptotic
expansion (\ref{da}).
\begin{thm}\label{thmXi}
Under assumptions (\ref{regab}), (\ref{divb}) and (\ref{sdz}) for any $\xvec\in\rit^d$, 
$s\in\rit$ and $T\in \rit$, the sequences $(\xxkep(\cdot;\xvec,s))$, 
for $k=1,\dots,m$, are bounded in $L^\infty([s,s+T])$ and we have:
\begin{equation}\label{AjoutCvgXkeps}
\begin{array}{l}
\ds \lim_{\eps\rightarrow 0}  \sup_{t\in[s,s+T]}
\big|\xxkep(t;\xvec,s)-\xxk(t,\tmsoe;\xvec,s) \big|  = 0,
\end{array}
\end{equation}
where $\xxk$ writes
\begin{multline}\label{Xi-Yi}
\xxk(t,\theta;\xvec,s) =\big\{\Nabla_z \zzvec(t,\theta;\yyze(t;\xvec,s))\big\}
\big\{ \yyk(t;\xvec,s)+
\theta\tilde{\aavec}^{k-1}(t,\theta,\yyze(t;\xvec,s),\dots,\yykmun(t;\xvec,s)) \big\}+
\\
\frac12\big\{\Nabla_z^{\,2} \zzvec(t,\theta;\yyze(t;\xvec,s))\big\}
\Big( 
\sum_{j=1}^{k-1}
\big\{\yyj(t;\xvec,s)+ \theta\tilde{\aavec}^{j-1}
(t,\theta,\yyze(t;\xvec,s),\dots,\yyjmun(t;\xvec,s))
\,,\,~~~~~~~~~~~~
\\ ~~~~~~~~~~~~~~~~~~~~~~~~~~~~
\yyvec^{k-j}(t;\xvec,s)+ \theta\tilde{\aavec}^{k-j-1}
(t,\theta,\yyze(t;\xvec,s),\dots,\yyvec^{k-j-1}(t;\xvec,s))\big\}
\Big)+
\\
\dots +\frac1{k!}\big\{\Nabla_z^{\,k} \zzvec(t,\theta;\yyze(t;\xvec,s))\big\}
\{\yyun(t;\xvec,s)+\theta\tilde{\aavec}^{0}(t,\theta,\yyze(t;\xvec,s))\}^k,
\end{multline}
with $\yyk$ solution to (\ref{eqYi}).
\end{thm}
We end this result list by giving explicit expressions. 
Using (\ref{al0}) in (\ref{ali}) we have the following 
expression of $\tilde \alpha^1$:
\begin{equation}\label{alun}
\begin{array}{l}
\ds \tilde \alpha^1(t,\theta,\yze,\yun)=
\big\{\Nabla_z \zzvec(t,\theta;\yze)\big\}^{-1}
\Big\{
\big\{\Nabla_x\avec(t,\theta,\zzvec(t,\theta;\yze))\big\} 
\big\{\Nabla_z \zzvec(t,\theta;\yze)\big\} 
\big\{\yun+\theta \tilde{\aavec}^{0}(t,\theta;\yze)\big\}-
\\ \ds ~~
\big\{\frac{\partial\Nabla_z \zzvec}{\partial t}(t,\theta;\yze) \big\}
\{\yun+\theta\tilde{\aavec}^{0}(t,\theta,\yze)\}-
\big\{\Nabla_z^{\,2} \zzvec(t,\theta;\yze)\big\}
\{\tilde{\alpha}^0(t,\yze) ,\yun+\theta\tilde{\aavec}^{0}(t,\theta,\yze)\}
\Big\}-
\\ \ds
~ \hfill \theta \Big(
 \{\Nabla_{y^0}\tilde{\aavec}^{0}(t,\theta,\yze)\}
\{\tilde{\avec}^0(t,\yze)\}+
\frac{\partial\tilde{\aavec}^{0}}{\partial t}(t,\theta,\yze)
\Big).
\end{array}
\end{equation}
On another hand, if $\bvec$ is linear and independent of $t$,
$\tilde \alpha^2$ is given by
\begin{equation}\label{defal2+}
\begin{array}{l}
\ds \tilde \alpha^2(t,\theta,\yze,\yun,\yde)=
\big\{\Nabla_z \zzvec(\theta;\yze)\big\}^{-1}
\Big\{
\big\{\Nabla_x\avec(t,\theta,\zzvec(\theta;\yze))\big\} 
\big\{\Nabla_z \zzvec(\theta;\yze)\big\} 
\big\{\yde+\theta \tilde{\aavec}^{1}(t,\theta;\yze,\yun)\big\}+
\\ \ds ~~
\frac12 
\big\{\Nabla_x^{\,2}\avec(t,\theta,\zzvec(\theta;\yze))\big\} 
\Big\{\big\{\Nabla_z \zzvec(\theta;\yze)\big\} 
\big\{\yun+\theta \tilde{\aavec}^{0}(t,\theta;\yze)\big\}\Big\}^2-
\\ \ds
~ \hfill \theta \Big(
 \{\Nabla_{y^0}\tilde{\aavec}^{1}(t,\theta,\yze,\yun)\}
\{\tilde{\avec}^0(t,\yze)\}+
 \{\Nabla_{y^1}\tilde{\aavec}^{1}(t,\theta,\yze,\yun)\}
\{\tilde{\avec}^1(t,\yze,\yun)\}+
\frac{\partial\tilde{\aavec}^{1}}{\partial t}(t,\theta,\yze,\yun)
\Big).
\end{array}
\end{equation}

~

The paper is organised as follow: in the second section we briefly 
prove the Theorems. Section \ref {SecApp} is then devoted to applications 
to models describing magnetic 
confinement fusion and isotope resonant separation experiments.
Among those applications, one involves a non uniform strong magnetic field,
in this case the computations are led using Maple.

\begin{rem}
With very little changes we can apply the previous results to the case when 
$\bvec\equiv\bvec(t,\theta,\xvec)$ also depends on $\theta$, as soon as the
regularity of $\bvec$ is enough and the assumption (\ref{sdz}) is realized.
\end{rem}

\section{Proof of the Theorems}\label{secae}
\subsection{Sketch of the proof of Theorem \ref{thmappk}} 
The key point to prove  the results of this paper is Theorem \ref{thmappk},
in the same spirit as in Schochet \cite{schochet:1994}. 
This result is classical and  known as Single Phase Averaging. 
We do not give a detailed
proof of it but only a formal sketch of it using expansion (\ref{expYeps})
of $\yyzeep$. 

First, because of the definition (\ref{al0}) of $\tilde \alpha^0$ and 
the assumptions (\ref{regab}) - (\ref{sdz})
we deduce that the function $(t,\yze)\mapsto\tilde \alpha^0(t,\tmsoe,\yze)$
is regular enough to ensure that, for any $\eps>0$,
$\yyzeep(\cdot;\xvec,s)$ exists, is unique on $[s,s+T]$ and remains in a
bounded set of $\rit^d$ independent of $\eps$.
\newline
Then expanding $\tilde{\alpha}^0(\yyzeep)$,
using $\yyzeep=\yyze+ 
\sum_{j\geq 1} \eps^j (\yyj+[\theta\tilde{\aavec}^{j-1}]_\eps)$, 
we obtain
\begin{equation}\label{expa}
\begin{array}{l}
\ds \tilde{\alpha}^0(\yyzeep)=\tilde{\alpha}^0(\yyze) +
\eps\big\{\Nabla_{y^0}\tilde{\alpha}^0(\yyze)\big\}
\big\{\yyun+ [\theta\tilde{\aavec}^{0}]_\eps\big\}+
\eps^2\bigg( \big\{\Nabla_{y^0}\tilde{\alpha}^0(\yyze)\big\}
\{\yyde+ [\theta\tilde{\aavec}^{1}]_\eps\}+~~~~
\\ \ds ~~~
\frac12\big\{\Nabla_{y^0}^{\,2}\tilde{\alpha}^0(\yyze)\big\}
\big\{\yyun+ [\theta\tilde{\aavec}^{0}]_\eps \big\}^2\bigg)+ 
\dots +
\eps^k\bigg(\big\{\Nabla_{y^0}\tilde{\alpha}^0(\yyze)\big\}
\big\{\yyk + [\theta\tilde{\aavec}^{k-1}]_\eps\big\}+
\\ \ds ~~~
\frac12\big\{ \Nabla_{y^0}^{\,2}\tilde{\alpha}^0(\yyze)  \big\}
\big(\sum_{j=1}^{k-1} \{\yyj+ [\theta\tilde{\aavec}^{j-1}]_\eps ,
\yyvec^{k-j}+ [\theta\tilde{\aavec}^{k-j-1}]_\eps \}\big)+
\\ \ds ~~~ \hfill
\dots +\frac{1}{k!}\big\{\Nabla_{y^0}^{\,k}\tilde{\alpha}^0(\yyze)\big\}
\big\{\yyun+ [\theta\tilde{\aavec}^{0}]_\eps \big\}^k
\bigg)+\dots~.  \hspace{-7pt}
\end{array}
\end{equation}
Plugging the expansion (\ref{expYeps}) in the dynamical system (\ref{eqy0ep})
and identifying the terms of the same order, we have
\begin{equation}\label{eqx0}
\fracp{\yyze}{t}+\fracp{[\theta\tilde{\aavec}^{0}]}{\theta} = 
\tilde \alpha^0(t,\theta,\yyze) ,
\end{equation}
at the order $-1$, and at the order $k-1$ for $k\geq1$ we have
\begin{equation}\label{eqxk}
\fracp{\yyk}{t}+\fracp{[\theta\tilde{\aavec}^{k}]}{\theta} =
\tilde \alpha^k(t,\theta,\yyze,\dots\yyk).
\end{equation}
Then equation (\ref{eqY0}) for $\yyze$, definition (\ref{defAitild}) 
of $\tilde{\aavec}^k$ and equation (\ref{eqYi}) for $\yyk$ follow,
ending the sketch of the proof of Theorem \ref{thmappk}.
~~~~~{\hfill\rule{2mm}{2mm}\vskip3mm \par}

\subsection{Proof of Theorem \ref{thmX0}}
Existence and uniqueness of $\xxep$ are consequences of existence 
and uniqueness of $\yyzeep$. Then the regularity of $\zzvec$ and
Theorem \ref{thmappk} allows to pass to the limit in (\ref{defy0eps})
in order to deduce (\ref{app0}) and (\ref{eqY0}) and thus  prove
the Theorem.
~~~~~{\hfill\rule{2mm}{2mm}\vskip3mm \par}

\subsection{Proof of Theorem \ref{thmXi}}
From Theorem \ref{thmappk}, we deduce
\begin{equation}
\yyzeep = \yyze + \eps (\yyun + [\theta \tilde \aavec^{0}]_\eps)
+ \dots  +\eps^{k-1} (\yykmun + [\theta \tilde \aavec^{k-2}]_\eps)+
\eps^{k} (\yykep + [\theta \tilde \aavec^{k-1}]_\eps),
\end{equation}
and then, from (\ref{defy0eps}) we deduce
\begin{multline}
\xxep=[\zzvec(\yyzeep)]_\eps=[\zzvec(\yyze)]_\eps+ 
\eps\big\{ \Nabla_z\zzvec(\yyzeep) \big\}_\eps 
\big\{\yyun + [\theta\tilde{\aavec}^{0}]_\eps \big\}+
\eps^2
\bigg(
\big\{ \Nabla_z\zzvec(\yyzeep) \big\}_\eps
\big\{\yyde + [\theta\tilde{\aavec}^{1}]_\eps \big\}+
\\ 
\frac12\big\{ \Nabla_z^{\,2}\zzvec(\yyzeep) \big\}_\eps
\big\{\yyun + [\theta\tilde{\aavec}^{0}]_\eps \big\}^2
\bigg)+\dots + 
\eps^k\bigg(
\big\{ \Nabla_z\zzvec(\yyzeep) \big\}_\eps
\big\{\yykep + [\theta\tilde{\aavec}^{k-1}]_\eps \big\}+
\\ 
\frac12\big\{ \Nabla_z^{\,2}\zzvec(\yyzeep) \big\}_\eps
\Big(
\sum_{j=1}^{k-1}\{\yyj+ [\theta\tilde{\aavec}^{j-1}]_\eps,
                  \yyvec^{k-j}+ [\theta\tilde{\aavec}^{k-j-1}]_\eps\}
\Big)+
\\ 
\dots+
\frac1{k!}\big\{\Nabla_z^{\,k} \zzvec(\yyze) \big\}
\{\yyun+[\theta\tilde{\aavec}^{0}]_\eps\}^k
\bigg)+o(\eps^k).
\end{multline}
Comparing this expansion with
\begin{equation}
\xxzeep = [\xxze]_\eps + \eps[\xxun]_\eps + \dots  +\eps^{k-1}[\xxkmun]_\eps +
\eps^{k} \xxkep,
\end{equation}
obtained from (\ref{defxkeps}) and making the process
$\eps\rightarrow 0$ give finally Theorem \ref{thmXi}.
~~~~~{\hfill\rule{2mm}{2mm}\vskip3mm \par}

\section{Application to the gyrokinetic plasma}\label{SecApp}
The dynamical systems associated with equations (\ref{vlglob}), (\ref{vlRLF})
and (\ref{vlCR}) are in the form of (\ref{da}), i.e.:
\begin{equation}\label{Rsd}
\frac{d}{dt}\begin{pmatrix}\xxep \\ \vvvec_\eps \end{pmatrix}=
\avec(t,\tmsoe,\xxep,\vvvec_\eps) + \frac 1\eps \bvec(t,\xxep,\vvvec_\eps)\, ,\; \;
\begin{pmatrix}\xxep(s;\xvec,\vvec,s)\\\vvvec_\eps(s;\xvec,\vvec,s)\end{pmatrix}=
\begin{pmatrix} \xvec\\ \vvec \end{pmatrix}
\end{equation}
with variable $(\xvec,\vvec)\in \rit^3\times\rit^3$ in place of $\xvec$ and with ad-hoc
fields $\avec$ and  $\bvec$. Then we can apply our result saying that the solution 
$(\xxep(t;\xvec,\vvec,s), \vvvec_\eps(t;\xvec,\vvec,s))$
can be expanded as
\begin{equation}\label{Rexp}
\begin{array}{l}\ds 
\xxep(t;\xvec,\vvec,s)=\xxze(t,\tmsoe;\xvec,\vvec,s)+
\eps\xxun (t,\tmsoe;\xvec,\vvec,s)+ \eps^2\xxde (t,\tmsoe;\xvec,\vvec,s)+\dots,
\\ \ds
\vvvec_\eps(t;\xvec,\vvec,s)=\vvvec^0(t,\tmsoe;\xvec,\vvec,s)+
\eps\vvvec^1 (t,\tmsoe;\xvec,\vvec,s)+ \eps^2\vvvec^2 (t,\tmsoe;\xvec,\vvec,s)+\dots,
\end{array}
\end{equation}
and that this expansion may be justified until any order if the regularity of
the fields is enough.

~

In this section, we first lead the computations in the case of Isotope Resonant Separation
Regime with the restriction that $\calm$ is constant and $\bbvec(t,\theta,\xvec)=0$ until
the order 1. Then, with the same restrictions, we treat the case 
of the Guiding Centre Regime until the order 2 and we give the result for the 
Finite Larmor Radius Regime until the order 0.
The forth example concerns the 
Guiding Centre Regime with a variable strong magnetic field and a constant  $\bbvec$
until the order 1.

From the physical point of view, the last example is relevant to understand 
the behaviour of a plasma in a tokamak. 
This example shows that the generic computations made before,
coupled with the use of Maple, enable to deduce the result relatively 
comfortably.

\subsection{Isotope Resonant Separation Regime with constant strong magnetic field}
In the case of Isotope Resonance Separation Regime, i.e. of equation (\ref{vlCR}),
$\avec$ and $\bvec$ are:
\begin{equation}\label{R4.11}
\avec(t,\theta,\xvec,\vvec) =
\begin{pmatrix} \vvec \\ 
\eevec(t,\theta,\xvec)+\vvec\times \bbvec(t,\theta,\xvec) \end{pmatrix}
\text{, ~ } 
\bvec(t,\xvec,\vvec) =\bvec(\vvec)=
\begin{pmatrix} 0 \\ \vvec\times\calm \end{pmatrix}.
\end{equation}
For simplicity, we restrict to the case when $\bbvec(t,\theta,\xvec)=0$
and when $\calm=\basevecun$ is a constant vector, 
($\basevecun$,$\basevecde$,$\basevectr$) being the frame of  $\rit^3$.
Then $\zzvec(t,\theta;\zvec,\wvec)$ and
$\{\Nabla_{z,w} \zzvec(t,\theta;\zvec,\wvec)\}^{-1}$ 
are given by 
\begin{equation}\label{R4.12}
\zzvec(t,\theta;\zvec,\wvec)=\begin{pmatrix} \zvec \\ 
R(\theta)\wvec \end{pmatrix}\; ,\;\;
\big\{\Nabla_{z,w} \zzvec(t,\theta;\zvec,\wvec)\big\}^{-1}=
\begin{pmatrix} I&0\\0&R(-\theta)\end{pmatrix},
\end{equation}
where $R(\theta)$ is
the matrix of the rotation of angle $-\theta$ around $\calm$,
\begin{equation}\label{defR}
R(\theta)=
\begin{pmatrix} 
 1 & 0 & 0 \\
 0 & \cos\theta & \sin\theta\\
 0 & -\sin\theta & \cos\theta
\end{pmatrix}.
\end{equation}
We have the following result.
\begin{thm}
If $\eevec(t,\theta,\xvec)$ is $C^1_b(\rit\times\rit\times\rit^3)$
and $2\pi-$periodic in $\theta$,
the first term of the expansion (\ref{Rexp}) of the solution 
$(\xxep(t;\xvec,\vvec,s),\vvvec_\eps(t;\xvec,\vvec,s))$ to
\begin{equation}\label{irsrds}
\frac{d\xxep}{dt} = \vvvec_\eps \, , \; \;
\frac{d\vvvec_\eps}{dt} = \eevec(t,\tmsoe,\xxep)+\frac 1\eps \vvvec_\eps \times \calm
\, , \; \;
\xxep(s;\xvec,\vvec,s)=\xvec \, ,\; \; \vvvec_\eps(s;\xvec,\vvec,s)=\vvec,
\end{equation}
is given by  \vspace{-10pt}
\begin{equation}
\xxze(t,\theta;\xvec,\vvec,s) = \yyze(t;\xvec,\vvec,s) \, ,\;\; 
\vvvec^0 (t,\theta;\xvec,\vvec,s) = R(\theta) \uuvec^0(t;\xvec,\vvec,s),
\end{equation}
where $(\yyze(t;\xvec,\vvec,s),\uuvec^0(t;\xvec,\vvec,s))$ is solution to
\begin{equation}
\frac{d\yyze}{dt} = \uuvec^0_\parallel
\, , \;\;
\frac{d\uuvec^0}{dt} = 
\frac{1}{2\pi}\int_{0}^{2\pi} \hspace{-5pt} R(-\theta)\eevec(t,\theta,\yyze)\, d\theta
\, , \;\;
\yyze(s;\xvec,\vvec,s)= \xvec
\, , \;\;
\uuvec^0(s;\xvec,\vvec,s)= \vvec.
\end{equation}
\end{thm}
\begin{proof}
In view of (\ref{al0}) and (\ref{a0tilde}), here we have
\begin{equation}
\tilde \alpha^0(t,\theta,\yze,\uvec^0) =
\begin{pmatrix} R(\theta)\uvec^0 \\ R(-\theta)\eevec(t,\theta,\yze)  \end{pmatrix}
\, ,\; \;
\tilde{\avec}^0(t,\yze,\uvec^0)=
\begin{pmatrix} 
\uvec^0_\parallel 
\\ \ds
\frac{1}{2\pi}\int_{0}^{2\pi} \hspace{-5pt} 
R(-\theta)\eevec(t,\theta,\yze)\, d\theta
\end{pmatrix} \, .
\end{equation}
Then the proof of the Theorem is straightforward.
\end{proof}
In order to obtain the system satisfied by the second term 
$(\yyun,\uuvec^1)$ of the expansion we notice
that $\calr(\theta)= -R(\frac{\pi}{2}+\theta)+R(\frac{\pi}{2})$ is such that 
$\int_0^\theta R(\sigma) \, d\sigma = \theta P +\calr(\theta),$
with $P$ the matrix of the orthogonal projection onto $\calm$. We have
\begin{equation}\label{defcalR}
\calr(\theta)=
\begin{pmatrix} 
 0& 0 & 0 \\
 0 & \sin\theta & 1-\cos\theta\\
 0 & \cos\theta-1 & \sin\theta
\end{pmatrix}\text{ and }
P=
\begin{pmatrix} 
 1 & 0 & 0 \\
 0 & 0 & 0 \\
 0 & 0 & 0
\end{pmatrix}.
\end{equation}
Then, from (\ref{defAitild}), we get:
\begin{equation}
\theta \tilde \aavec^0(t,\theta,\yze,\uvec^0) =
\begin{pmatrix}
\calr(\theta) \uvec^0
\\ \ds 
\big(\int_0^\theta \hspace{-3mm}d\sigma
-\frac\theta{2\pi} \hspace{-1mm}\int_0^{2\pi} \hspace{-4.5mm}d\sigma\big)
\big(R(-\sigma)\eevec(t,\sigma,\yze)  \big)
\end{pmatrix},
\end{equation}
where we denote 
$\int_0^\theta f(\sigma) \,d\sigma
-\frac{\theta}{2\pi} \int_0^{2\pi} f(\sigma) \,d\sigma$
by
$\big(\int_0^\theta d\sigma
-\frac{\theta}{2\pi} \int_0^{2\pi} d\sigma\big)(f(\sigma)).$
We also get
\begin{equation}
\begin{array}{l}
\tilde\alpha^1(t,\theta,\yze,\uvec^0,\yun,\uvec^1) = 
\begin{pmatrix}\ds
R(\theta)
\big(
\uvec^1+\Big(\int_0^\theta \hspace{-3mm}d\sigma
-\frac\theta{2\pi} \hspace{-1mm}\int_0^{2\pi} \hspace{-4.5mm}d\sigma\Big)
\Big(R(-\sigma)\eevec(t,\sigma,\yze)  \Big)
\big)
\\ \ds
R(-\theta) \Nabla_x\eevec(t,\theta,\yze)
(\yun+\calr(\theta) \uvec^0)
\end{pmatrix} -
\\ \ds ~~~~~~~~~
\Big(\int_0^\theta \hspace{-3mm}d\sigma
-\frac\theta{2\pi} \hspace{-1mm}\int_0^{2\pi} \hspace{-4.5mm}d\sigma\Big)
\Big[
\begin{pmatrix}\ds
R(\sigma)
\big( \frac{1}{2\pi}\int_{0}^{2\pi} \hspace{-5pt} 
R(-\varsigma)\eevec(t,\varsigma,\yze)\, d\varsigma \big)
\\ \ds
R(-\sigma) \Nabla_x\eevec(t,\sigma,\yze)
\uvec^0_\parallel
\end{pmatrix}
+
\begin{pmatrix}\ds
0 \\ \ds R(-\sigma) \fracp{\eevec}{t}(t,\sigma,\yze)
\end{pmatrix}
\Big]
\end{array}
\end{equation}
Then using
$\big(\int_0^\theta d\sigma-\frac\theta{2\pi} \int_0^{2\pi} d\sigma\big)
\big( R(\sigma)\big) = \calr(\theta) $
and
$\frac 1{2\pi} \int_0^{2\pi}\calr(\theta) \, d\theta = R(\frac \pi 2)-P,$
we obtain
\begin{multline}
\tilde\avec^1(t,\yze,\uvec^0,\yun,\uvec^1) =
\\ \ds ~~~~~~
\begin{pmatrix}\ds
\uvec^1_\parallel+
\frac 1{2\pi} \int_0^{2\pi} \Big(\int_0^\theta \hspace{-3mm}d\sigma
-\frac\theta{2\pi} \hspace{-1mm}\int_0^{2\pi} \hspace{-4.5mm}d\sigma\Big)
\Big(R(\theta-\sigma)\eevec(t,\sigma,\yze)  \Big)d\theta
\displaybreak[0]\\ \ds 
\big(\frac{1}{2\pi}\int_0^{2\pi}
R(-\theta) \Nabla_x\eevec(t,\theta,\yze)
d\theta\big)\yun
+
\big(\frac{1}{2\pi}\int_0^{2\pi}
R(-\theta) \Nabla_x\eevec(t,\theta,\yze)\calr(\theta)
d\theta\big)\uvec^0
\end{pmatrix} -
\\ \ds ~~~~~~~~~
\begin{pmatrix}\ds
 (R(\frac{\pi}{2})-P)
\big( \frac{1}{2\pi}\int_{0}^{2\pi} \hspace{-5pt} 
R(-\varsigma)\eevec(t,\varsigma,\yze)\, d\varsigma \big)
\\ \ds
\big(\frac{1}{2\pi}\int_0^{2\pi}
\Big(
\int_0^\theta \hspace{-3mm}d\sigma
-\frac\theta{2\pi} \hspace{-1mm}\int_0^{2\pi} \hspace{-4.5mm}d\sigma
\Big)
\Big( R(-\sigma) \Nabla_x\eevec(t,\sigma,\yze) \Big)
d\theta\big)\uvec^0_\parallel
\end{pmatrix}
-
\\ \ds \hfill
\begin{pmatrix}\ds
0 \\ \ds 
\frac{1}{2\pi}\int_0^{2\pi}
\Big(
\int_0^\theta \hspace{-3mm}d\sigma
-\frac\theta{2\pi} \hspace{-1mm}\int_0^{2\pi} \hspace{-4.5mm}d\sigma
\Big)
\Big(R(-\sigma) \fracp{\eevec}{t}(t,\sigma,\yze )\Big)d\theta
\end{pmatrix}.~~~
\end{multline}

Hence we can state the following Theorem.
\begin{thm}
If $\eevec(t,\theta,\xvec)$ is $C^2_b(\rit\times\rit\times\rit^3)$
and $2\pi-$periodic in $\theta$,
the second term of the expansion (\ref{Rexp}) of the solution 
$(\xxep(t;\xvec,\vvec,s),\vvvec_\eps(t;\xvec,\vvec,s))$ to
(\ref{irsrds}) is given by
\begin{equation}
\begin{array}{l}\ds
\xxun(t,\theta;\xvec,\vvec,s) =\yyun(t;\xvec,\vvec,s)
+\calr(\theta) \uuvec^0(t;\xvec,\vvec,s),
\\ \ds
\vvvec^1(t,\theta;\xvec,\vvec,s) = R(\theta) \uuvec^1(t;\xvec,\vvec,s)+R(\theta)   
\big(\int_0^\theta \hspace{-3mm}d\sigma
-\frac\theta{2\pi} \hspace{-1mm}\int_0^{2\pi} \hspace{-4.5mm}d\sigma\big)
\big(R(-\sigma)\eevec(t,\sigma,\yyze(t;\xvec,\vvec,s))  \big),
\end{array}
\end{equation}
where $(\yyun,\uuvec^1)$ is solution to
\begin{equation}
\begin{array}{l}\ds
\frac{d\yyun}{dt} =
\uuvec^1_\parallel+
\frac{1}{2\pi}\int_0^{2\pi}
\Big(\int_0^\theta \hspace{-3mm}d\sigma
-\frac\theta{2\pi} \hspace{-1mm}\int_0^{2\pi} \hspace{-4.5mm}d\sigma\Big)
\Big(R(\theta-\sigma)\eevec(t,\sigma,\yyze)  \Big) d\theta
-
~~~\\ \hfill \ds 
(R(\frac{\pi}{2})-P)
\big( \frac{1}{2\pi}\int_{0}^{2\pi} \hspace{-5pt} 
R(-\varsigma)\eevec(t,\varsigma,\yyze)\, d\varsigma \big),
\\ \ds 
\frac{d\uuvec^1}{dt} =
\big(\frac{1}{2\pi}\int_0^{2\pi}
R(-\theta) \Nabla_x\eevec(t,\theta,\yyze)
d\theta\big)\yyun +
\big(\frac{1}{2\pi}\int_0^{2\pi}
R(-\theta) \Nabla_x\eevec(t,\theta,\yyze)\calr(\theta)
d\theta\big)\uuvec^0-
~~~\\ ~~~~~~~~~~~~\ds 
\big(\frac{1}{2\pi}\int_0^{2\pi}
\Big(
\int_0^\theta \hspace{-3mm}d\sigma
-\frac\theta{2\pi} \hspace{-1mm}\int_0^{2\pi} \hspace{-4.5mm}d\sigma
\Big)
\Big( R(-\sigma) \Nabla_x\eevec(t,\sigma,\yyze) \Big)
d\theta\big)\uuvec^0_\parallel-
\\ \ds \hfill 
\frac{1}{2\pi}\int_0^{2\pi}
\Big(
\int_0^\theta \hspace{-3mm}d\sigma
-\frac\theta{2\pi} \hspace{-1mm}\int_0^{2\pi} \hspace{-4.5mm}d\sigma
\Big)
\Big( R(-\sigma) \fracp{\eevec}{t}(t,\sigma,\yyze )\Big)d\theta,
\\ \ds
\yyun(s;\xvec,\vvec,s) =0, \;\;\; \uuvec^1(s;\xvec,\vvec,s) =0.
\end{array}
\end{equation}
\end{thm}
\subsection{Guiding Centre Regime with constant strong magnetic field}
In the case of the Guiding Centre Regime, i.e. of
equation (\ref{vlglob}), we have 
\begin{equation}\label{R4.7}
\avec(t,\theta,\xvec,\vvec) =\avec(t,\xvec,\vvec) =
\begin{pmatrix} \vvec \\ \eevec(t,\xvec)+\vvec\times \bbvec(t,\xvec) \end{pmatrix},
\text{ and }
\bvec(t,\xvec,\vvec) =\bvec(\vvec)=
\begin{pmatrix} 0 \\ \caln + \vvec\times\calm \end{pmatrix},
\end{equation}

As previously, we make the restriction
$\bbvec(t,\theta,\xvec)=\caln=0$ and $\calm=\basevecun$. 
Since this situation is similar to the previous one with the only difference
that $\eevec(t,\xvec)$ does not depend on $\theta$ we can directly deduce 
the following Theorem.
\begin{thm} 
If $\eevec(t,\xvec)$ is $C^2_b(\rit\times\rit^3)$,
the first and second  terms of the expansion (\ref{Rexp}) of the solution 
$(\xxep(t;\xvec,\vvec,s),$ $\vvvec_\eps(t;\xvec,\vvec,s))$ to
\begin{equation}\label{gcrds}
\frac{d\xxep}{dt} = \vvvec _\eps\, , \; \;
\frac{d\vvvec_\eps}{dt} = \eevec(t,\xxep) + \frac 1\eps \vvvec _\eps \times \calm
\, , \; \;
\xxep(s;\xvec,\vvec,s)=\xvec \, ,\; \; \vvvec_\eps(s;\xvec,\vvec,s)=\vvec,
\end{equation}
are given by
\begin{equation}
\xxze(t,\theta;\xvec,\vvec,s) = \yyze(t;\xvec,\vvec,s) \, ,\;\; 
\vvvec^0 (t,\theta;\xvec,\vvec,s) = R(\theta) \uuvec^0(t;\xvec,\vvec,s),
\end{equation}
and
\begin{equation}
\begin{array}{l}\ds
\xxun(t,\theta;\xvec,\vvec,s) =\yyun(t;\xvec,\vvec,s)
+\calr(\theta)\uuvec^0(t;\xvec,\vvec,s),
\\ \ds
\vvvec^1(t,\theta;\xvec,\vvec,s) = R(\theta) \uuvec^1(t;\xvec,\vvec,s)
+ \calr(\theta)
\eevec(t,\yyze(t;\xvec,\vvec,s)),
\end{array}
\end{equation}
where  $(\yyze(t;\xvec,\vvec,s),\uuvec^0(t;\xvec,\vvec,s))$ is solution to
\begin{equation}
\frac{d\yyze}{dt} = \uuvec^0_\parallel
\, , \;\;
\frac{d\uuvec^0}{dt} = 
\eevec_\parallel(t,\yyze)
\, , \;\;
\yyze(s;\xvec,\vvec,s)= \xvec
\, , \;\;
\uuvec^0(s;\xvec,\vvec,s)= \vvec,
\end{equation}
and where $(\yyun,\uuvec^1)$ is solution to
\begin{equation}
\begin{array}{l}\ds
\frac{d\yyun}{dt} =
\uuvec^1_\parallel + (R(\frac{\pi}{2})-P) \eevec(t,\yyze)
= \uuvec^1_\parallel + \eevec(t,\yyze) \times \calm,
\\ \ds 
\frac{d\uuvec^1}{dt} =
P\, \Nabla_x\eevec(t,\yyze)
\yyun
+
\frac 12 \text{\rm tr\hspace{-1pt}} \big((I-P)\Nabla_x\eevec(t,\yyze)\big)
 \big(R(-\frac{\pi}{2})-P\big)\uuvec^0+
~~~~~~~\\ \ds \hfill ~~~~~~~~~
\frac 12 \text{\rm tr\hspace{-1pt}} 
\big((R(-\frac{\pi}{2})-P)\Nabla_x\eevec(t,\yyze) \big)
(I-P)\uuvec^0-
(R(-\frac{\pi}{2})-P)\Nabla_x\eevec(t,\yyze)
\uuvec^0_\parallel-
~~~~~\\ \ds \hfill 
(R(-\frac{\pi}{2})-P)\fracp{\eevec}{t}(t,\yyze ),
\\ \ds
\yyun(s;\xvec,\vvec,s) =0, \;\;\; \uuvec^1(s;\xvec,\vvec,s) =0.
\end{array}
\end{equation}
\end{thm}
In the computations leading to this Theorem, we use, among other formula,
\begin{equation}
\Big(
\int_0^\theta \hspace{-3mm}d\sigma
-\frac\theta{2\pi} \hspace{-1mm}\int_0^{2\pi} \hspace{-4.5mm}d\sigma
\Big)
\Big( R(-\sigma)\Big) = -\calr(-\theta),\;\; 
\frac 1{2\pi} \int_0^{2\pi}\hspace{-4pt}-\calr(-\theta) \, d\theta = R(-\frac \pi 2)-P
,\;\; 
-R(\theta) \calr(-\theta) =\calr(\theta),
\end{equation}
and
\begin{equation}
 \begin{array}{l}\ds
\frac{1}{2\pi}\int_0^{2\pi}
R(-\theta) \Nabla_x\eevec(t,\yyze)\calr(\theta)
d\theta =
\frac 12 \text{tr\hspace{-1pt}} \big((I-P)\Nabla_x\eevec(t,\yyze)\big)
 \big(R(-\frac{\pi}{2})-P\big)
+
~~~~~~~~~~~ \\ \ds \hfill
\frac 12 \text{\rm tr\hspace{-1pt}} 
\big((R(-\frac{\pi}{2})-P)\Nabla_x\eevec(t,\yyze) \big)
(I-P).
 \end{array}
\end{equation}

~

Now we turn to the characterisation of $(\xxde,\vvvec^2)$.
For this purpose we notice that here
\begin{equation}\label{gcral1}
\begin{array}{l}\ds
\tilde\alpha^1(t,\theta,\yze,\uvec^0,\yun,\uvec^1)=
\\ ~~~
\begin{pmatrix}\ds
R(\theta) \uvec^1 + \calr(\theta) \eevec(t,\yze )
\\ \ds
R(-\theta)\Nabla_x\eevec(t,\yze) \yun + 
R(-\theta) \Nabla_x\eevec(t,\yze)\calr(\theta)\uvec^0 +
\calr(-\theta) \Nabla_x\eevec(t,\yze) \uvec^0_\parallel +
\calr(-\theta) \fracp{\eevec}{t}(t,\yze) 
\end{pmatrix}.
\end{array}
\end{equation}
In order to get now $\theta\tilde\aavec^1$ we need first to compute
\begin{equation} 
\big(\int_0^\theta \hspace{-3mm}d\sigma
-\frac\theta{2\pi} \hspace{-1mm}\int_0^{2\pi} \hspace{-4.5mm}d\sigma\big)
\big(\calr(\sigma)\big) = I- R(\theta)\, ,\;\;
\big(\int_0^\theta \hspace{-3mm}d\sigma
-\frac\theta{2\pi} \hspace{-1mm}\int_0^{2\pi} \hspace{-4.5mm}d\sigma\big)
\big(\calr(-\sigma\big) )= I- R(-\theta).
\end{equation}
Secondly, 
\begin{equation}
\begin{array}{l}\ds
\big(\int_0^\theta \hspace{-3mm}d\sigma
-\frac\theta{2\pi} \hspace{-1mm}\int_0^{2\pi} \hspace{-4.5mm}d\sigma\big)
\big(R(-\sigma) \Nabla_x\eevec(t,\yze)\calr(\sigma)\big)= 
\vspace{-10pt}
\\ ~~~~~ \ds
P \Nabla_x \eevec(t,\yze) (I-R(\theta)) +
(I-P)\Nabla_x \eevec(t,\yze)
\begin{pmatrix}
0 & 0 & 0 \\
0 & \ds\frac{\sin^2\theta}{2} & \ds\sin\theta- \frac{\sin2\theta}{4} \\
0 & \ds\frac{\sin2\theta}{4}- \sin\theta & \ds\frac{\sin^2\theta}{2}
\end{pmatrix}+
~~~~~~~~~~~~~\\ ~~~~~ \hfill \ds
\big( R(-\frac\pi 2)-P \big)\Nabla_x \eevec(t,\yze)
\begin{pmatrix}
0 & 0 & 0 \\
0 & \ds-\frac{\sin2\theta}{4} & \ds 1-\cos\theta-\frac{\sin^2\theta}{2} \\
0 & \ds\frac{\sin^2\theta}{2}+\cos\theta-1& \ds\frac{-\sin2\theta}{4}
\end{pmatrix}\, ,
\end{array}\vspace{-10pt}
\end{equation}
which also reads
\begin{equation}
\begin{array}{l}\ds
\big(\int_0^\theta \hspace{-3mm}d\sigma
-\frac\theta{2\pi} \hspace{-1mm}\int_0^{2\pi} \hspace{-4.5mm}d\sigma\big)
\big(R(-\sigma) \Nabla_x\eevec(t,\yze)\calr(\sigma)\big)=
P \Nabla_x \eevec(t,\yze) (I-R(\theta)) +
~~~~~~~~~~~~\\ \hfill \ds
\frac{1}{2}\big(-\calr(-\theta)+ R(\frac\pi 2)-P)
\Nabla_x \eevec(t,\yze) 
\big(\calr(\theta) + R(\frac\pi 2)-P \big) +
\frac{1}{2} \big(R(-\frac\pi 2)-P \big) \Nabla_x \eevec(t,\yze) \big(R(\frac\pi 2)-P \big).
\end{array}
\end{equation}
Hence integrating (\ref{gcral1}) we have
\begin{multline}
\ds
\theta\tilde\aavec^1(t,\theta,\yze,\uvec^0,\yun,\uvec^1)=
\begin{pmatrix}\ds
\calr(\theta) \uvec^1 +(I-R(\theta))\eevec(t,\yze)
\\ \ds
-\calr(-\theta)\Nabla_x \eevec(t,\yze)\yun -
 (I-R(-\theta)) 
\Big(\Nabla_x \eevec(t,\yze)\uvec^0_\parallel + \fracp{\eevec}{t}(t,\yze)\Big)
\end{pmatrix}+
~\displaybreak[0]\\ ~~~~~\hfill \ds
\begin{pmatrix}\ds
0
\\ \ds
\Big(P \Nabla_x \eevec(t,\yze) (I-R(\theta)) +
\frac{1}{2}\big(-\calr(-\theta)+ R(\frac\pi 2)-P)
\Nabla_x \eevec(t,\yze) 
\big(\calr(\theta) + R(\frac\pi 2)-P \big)\Big)\uvec^0
\end{pmatrix}+
~\\ ~~~~~\hfill \ds
\begin{pmatrix}\ds
0
\\ \ds
\Big(\frac{1}{2} \big(R(-\frac\pi 2)-P \big) 
\Nabla_x \eevec(t,\yze) \big(R(\frac\pi 2)-P \big)\uvec^0
\end{pmatrix}.
\end{multline}

In order to obtain the expression of $\tilde\alpha^2$,
we need to compute 
\begin{multline}\ds
\big\{\Nabla_{z,w} \zzvec(t,\theta;\yze,\uvec^0)\big\}^{-1}
\big\{\Nabla_{x,v}\avec(t,\theta,\zzvec(t,\theta;\yze,\uvec^0))\big\} 
\big\{\Nabla_{z,w} \zzvec(t,\theta;\yze,\uvec^0)\big\} 
\\ \ds \hfill
\Bigg\{\hspace{-3pt}\begin{pmatrix} \yde\\ \uvec^2 \end{pmatrix}\hspace{-3pt}
+\theta \tilde{\aavec}^{1}(t,\theta,\yze,\uvec^0,\yun,\uvec^1)\Bigg\} =
~~~~~\displaybreak[0] \\~~~~~~\ds
\begin{pmatrix}
0 & R(\theta) \\
R(-\theta)\Nabla_x \eevec(t,\yze) & 0 
\end{pmatrix}
\bigg(\hspace{-3pt}\begin{pmatrix} \yde\\ \uvec^2 \end{pmatrix}\hspace{-3pt}
+\theta \tilde{\aavec}^{1}(t,\theta,\yze,\uvec^0,\yun,\uvec^1)\bigg) =
\displaybreak[0] \\
\\~~~~~~~~~\ds
\begin{pmatrix} \ds 
  \begin{array}{c}\ds
    R(\theta)\uvec^2+\calr(\theta)\Nabla_x \eevec(t,\yze)\yun -
    (R(\theta)-I)
    \Big(\Nabla_x \eevec(t,\yze)\uvec^0_\parallel + \fracp{\eevec}{t}(t,\yze) \Big)+
    \\ 
    \Big(
    P \Nabla_x \eevec(t,\yze) (I-R(\theta))+\frac 12 \big( R(\frac\pi 2)-P)
    \Nabla_x \eevec(t,\yze) \big(\calr(\theta) + R(\frac\pi 2)-P \big)+
    \\
    \frac{1}{2} 
    \big(R(\theta-\frac\pi 2)-P \big) \Nabla_x \eevec(t,\yze) \big(R(\frac\pi 2)-P \big)
    \Big)\uvec^0 
  \end{array}
  \\
  \\ \ds
  R(-\theta)\Nabla_x \eevec(t,\yze) 
  \Big(\yde+\calr(\theta)\uvec^1+(I-R(\theta))\eevec(t,\yze) \Big) 
\end{pmatrix},~
\end{multline}

\begin{equation}
\begin{array}{l}\ds
\big\{\Nabla_{x,v}^{\,2}\avec(t,\theta,\zzvec(t,\theta;\yze,\uvec^0))\big\} 
\Bigg\{\Big\{\Nabla_{z,w} \zzvec(t,\theta;\yze,\uvec^0) \Big\}
\Big\{
\hspace{-3pt}\begin{pmatrix} \yun\\ \uvec^1 \end{pmatrix}\hspace{-3pt}
+\theta \tilde{\aavec}^{0}(t,\theta,\yze,\uvec^0) 
\Big\}
\Bigg\}^2=
~~~~\\~~~~~~~~~\ds
\begin{pmatrix}
  0
  \\
  \Big\{\Nabla_x^{\,2} \eevec(t,\yze) \Big\}^2
  \Big( \{\yun,\yun \} +2 \{\yun,\calr(\theta)\uvec^0 \} 
      + \{\calr(\theta)\uvec^0,\calr(\theta)\uvec^0 \} \Big)
\end{pmatrix}.
\end{array}
\end{equation}
A direct but heavy computation also leads the derivatives 
$\ds  \Nabla_{\yze,\uvec^0} \tilde\aavec^1(t,\theta,\yze,\uvec^0,\yun,\uvec^1)
\,\tilde\avec^0(t,\yze,\uvec^0)$,
$\ds  \Nabla_{\yun,\uvec^1} \tilde\aavec^1(t,\theta,\yze,\uvec^0,\yun,\uvec^1)
\,\tilde\avec^1(t,\yze,\uvec^0,\yun,\uvec^1)$,
and $\ds \fracp{\tilde\aavec^1}{t}(t,\theta,\yze,\uvec^0,\yun,\uvec^1)$.

Then in view of (\ref{defal2+}), we get
$\tilde\alpha^2$, which, integrating with respect to $\theta$ leads to
$\tilde\avec^2$ and to the equation satisfied by $(\xxde,\vvvec^2)$.
Then we have the following Theorem.
\begin{thm}
If $\eevec(t,\xvec)$ is $C^3_b(\rit\times\rit^3)$,
the third term of the expansion (\ref{Rexp}) of the solution 
$(\xxep(t;\xvec,\vvec,s),$ $\vvvec_\eps(t;\xvec,\vvec,s))$ to (\ref{gcrds})
is given by
\begin{equation}
\begin{array}{l}\ds
\xxde(t,\theta;\xvec,\vvec,s)= \yyde(t;\xvec,\vvec,s)
+\calr(\theta) \uuvec^1 +(I-R(\theta))\eevec(t,\yyze(t;\xvec,\vvec,s)),
\\ \ds
\vvvec^2(t,\theta;\xvec,\vvec,s) = R(\theta) \uuvec^2(t;\xvec,\vvec,s)+
\calr(\theta) \Nabla_x \eevec(t,\yyze(t;\xvec,\vvec,s)) \yyun(t;\xvec,\vvec,s)-
~~~~\\~~~~~~ \ds
(R(\theta)-I)
 \Nabla_x \eevec(t,\yyze(t;\xvec,\vvec,s)
\Big(\uuvec^0_\parallel(t;\xvec,\vvec,s)+\fracp{\eevec}{t}(t,\yyze(t;\xvec,\vvec,s))
\Big)+
\\~~~~~~ \ds
\Big(
P \Nabla_x \eevec(t,\yyze(t;\xvec,\vvec,s)) (I-R(\theta))+\frac 12 \big( R(\frac\pi 2)-P)
\Nabla_x \eevec(t,\yyze(t;\xvec,\vvec,s)) \big(\calr(\theta) + R(\frac\pi 2)-P \big)+
\\~~~~~~ \hfill \ds 
\frac{1}{2} 
\big(R(\theta-\frac\pi 2)-P \big) \Nabla_x \eevec(t,\yyze(t;\xvec,\vvec,s)) 
\big(R(\frac\pi 2)-P \big)
\Big)\uuvec^0(t;\xvec,\vvec,s).
\end{array}
\end{equation}
Moreover $(\yyde,\uuvec^2)$ is solution to
\begin{equation}
\begin{array}{l}\ds
\frac{d\yyde}{dt}=
\uuvec^2_\parallel+(R(\frac{\pi}{2})-P)\Nabla_x \eevec(t,\yyze)\yyun +
(I-P)
\Big(\Nabla_x \eevec(t,\yyze)\uuvec^0_\parallel + \fracp{\eevec}{t}(t,\yyze) \Big)+
\\ \ds ~~
\Big(
P \Nabla_x \eevec(t,\yyze) (I-P)+ \big( R(\frac\pi 2)-P)
\Nabla_x \eevec(t,\yyze) \big(R(\frac\pi 2)-P \big)\Big)\uuvec^0-
\\ \ds ~~
\frac 12 \text{\rm tr\hspace{-1pt}} \big((I-P)\Nabla_x\eevec(t,\yyze)\big)
\big(I-P)\big)\uuvec^0-
\frac 12 \text{\rm tr\hspace{-1pt}} 
\big((R(-\frac{\pi}{2})-P)\Nabla_x\eevec(t,\yyze) \big)
(R(\frac{\pi}{2})-P)\uuvec^0,
\end{array}
\end{equation}
\begin{multline}\ds
\frac{d\uuvec^2}{dt}=
P \Nabla_x \eevec(t,\yyze) \yyde+
\\ \ds ~~~~
\frac 12 \text{\rm tr\hspace{-1pt}} \big((I-P)\Nabla_x\eevec(t,\yyze)\big)
\big(R(-\frac{\pi}{2})-P\big)\uuvec^1+
\frac 12 \text{\rm tr\hspace{-1pt}} 
\big((R(-\frac{\pi}{2})-P)\Nabla_x\eevec(t,\yyze) \big)(I-P)\uuvec^1+
\hfill
\displaybreak[0] \\ \ds ~~~~
P \Nabla_x \eevec(t,\yyze)\eevec(t,\yyze)
-\Big(P \Nabla_x \eevec(t,\yyze)P +
\frac 12 \text{\rm tr\hspace{-1pt}}\big((I-P) \Nabla_x\eevec(t,\yyze)\big)(I-P)+
\hfill
\displaybreak[0] \\ \ds ~~~~
\frac 12 \text{\rm tr\hspace{-1pt}}
\big( (R(-\frac{\pi}{2})-P)\Nabla_x\eevec(t,\yyze)\big)(R(-\frac{\pi}{2})-P)
\Big)\eevec(t,\yyze) +
\hfill
\displaybreak[0] \\ \ds ~~~~
\frac 12\Big\{\Nabla_x^{\,2} \eevec(t,\yyze) \Big\}^2
\Big( \{\yyun,\yyun \} +2 \{\yyun,\uuvec^0_\parallel \} + 
\{(I-P)\uuvec^0,(I-P)\uuvec^0 \}+
\hfill
\\ \ds ~
\hfill
\{(R(-\frac{\pi}{2})-P)\uuvec^0,(R(-\frac{\pi}{2})-P)\uuvec^0 \} \Big)+
\displaybreak[0] \\ \ds ~~~~
(P-R(-\frac \pi 2))\{\Nabla_x^{\,2} \eevec(t,\yyze)\}\{\yyun,\uuvec^0_\parallel\} +
(I-P)\Big(\{\Nabla_x^{\,2} \eevec(t,\yyze)\}\{\uuvec^0_\parallel,\uuvec^0_\parallel\}
+ \fracp{\Nabla_x\eevec}{t}(t,\yyze)\uuvec^0_\parallel\Big)-
\hfill
\displaybreak[0] \\ \ds ~~~~
P \{\Nabla_x^{\,2} \eevec(t,\yyze)\} \{(I-P)\uuvec^0,\uuvec^0_\parallel\} -
\hfill
\displaybreak[0] \\ \ds ~~~~
\frac{1}{4}(I-P)\{\Nabla_x^{\,2} \eevec(t,\yyze)\} 
\{(I-P)\uuvec^0,\uuvec^0_\parallel\}-
\frac{3}{4} \big(R(-\frac\pi 2)-P \big) 
\{\Nabla_x^{\,2} \eevec(t,\yyze)\}
\{ \big(R(\frac\pi 2)-P \big)\uuvec^0,\uuvec^0_\parallel\}-
\hfill
\displaybreak[0] \\ \ds ~~~~
\Big(
-(I-P) 
\Big(\Nabla_x \eevec(t,\yyze)\Big)+
P \Nabla_x \eevec(t,\yyze) (I-P)+
\hfill
\\ \ds ~
\hfill
\frac{1}{4}(I-P)\Nabla_x \eevec(t,\yyze)(I-P) +
\frac{3}{4} \big(R(-\frac\pi 2)-P \big) \Nabla_x \eevec(t,\yyze) 
\big(R(\frac\pi 2)-P \big)\Big)\eevec_\parallel(t,\yyze)+
\displaybreak[0] \\ \ds ~~~~
(P-R(-\frac \pi 2))\Nabla_x\eevec(t,\yyze) \uuvec^1_\parallel+
(P-R(-\frac \pi 2))\Nabla_x\eevec(t,\yze) (R(\frac \pi 2)-P) \eevec(t,\yyze)+
\hfill
\displaybreak[0] \\ \ds ~~~~
(P-R(-\frac \pi 2))\fracp{\Nabla_x \eevec}{t}(t,\yyze)\yyun -
(I-P) 
\Big(\fracp{\Nabla_x \eevec}{t}(t,\yyze)\uuvec^0_\parallel + 
\fracp{^2\eevec}{t^2}(t,\yyze)\Big)-
P \fracp{\Nabla_x \eevec}{t}(t,\yyze) (I-P) -
\displaybreak[0] \\ \ds ~~~~
\Big( \frac{1}{4}(I-P)\fracp{\Nabla_x \eevec}{t}(t,\yyze) (I-P) +
\frac{3}{4} \big(R(-\frac\pi 2)-P \big) 
\fracp{\Nabla_x \eevec}{t}(t,\yyze) \big(R(\frac\pi 2)-P \big)\Big)\uuvec^0.
\end{multline}
\end{thm}

\subsection{Finite Larmor Radius Regime with constant strong magnetic field}
In the case of Finite Larmor Radius Regime (\ref{vlRLF}) with $\calm=\basevecun$, 
we have
\begin{equation}\label{R4.9}
\avec(t,\theta,\xvec,\vvec) =\begin{pmatrix} \vvec_\parallel\\    
\evec(t,\xvec)\end{pmatrix}
\text{ and }
\bvec(t,\xvec,\vvec)=\bvec(\vvec) =\begin{pmatrix}  \vvec_\perp\\ 
\vvec\times\calm\end{pmatrix},
\end{equation}
\begin{equation}\label{R4.10}
\zzvec(t,\theta;\zvec,\wvec)=\begin{pmatrix} \zvec +\calr(\theta)\wvec \\ 
R(\theta)\wvec\end{pmatrix}, \;\;
\big\{\Nabla_{z,w}\zzvec (t,\theta;\xvec,\vvec)\big\}^{-1}=
\begin{pmatrix} I&\calr(-\theta)\\0&R(-\theta)\end{pmatrix}.
\end{equation}
with $R$ and $\calr$ defined by (\ref{defR}) and (\ref{defcalR}).

Here  we give the result for this case only for the order 0. We have
\begin{equation}
\tilde{\alpha}^0(t,\theta,\yze,\uvec^0)= \begin{pmatrix}\ds
\uvec^0_\parallel +\calr(-\theta) \evec(\yze+\calr(\theta)\uvec^0,t) \\ \ds
R(-\theta)\evec(\yze+\calr(\theta)\uvec^0,t) \end{pmatrix},
\end{equation}
{and}
\begin{equation}
\tilde{\avec}^0(t,\yze,\uvec^0)= \begin{pmatrix}\ds
\uvec^0_\parallel +
\frac1{2\pi}\int_0^{2\pi}
\calr(-\theta) \evec(\yze+\calr(\theta)\uvec^0,t)\,d\theta \\ \ds
\frac1{2\pi}\int_0^{2\pi}
R(-\theta)\evec(\yze+\calr(\theta)\uvec^0,t)\,d\theta \end{pmatrix}.
\end{equation}
Hence we have the following Theorem.
\begin{thm}
If we assume that $\eevec(t,\xvec)$ is $C^1_b(\rit\times\rit^3)$,
the first term of the expansion (\ref{Rexp}) of the solution 
$(\xxep(t;\xvec,\vvec,s),\vvvec_\eps(t;\xvec,\vvec,s))$ to
\begin{equation}\label{lrrds}
\frac{d\xxep}{dt} = {\vvvec_{\hspace{-3pt}\eps}}_\parallel 
+\frac 1\eps {\vvvec_{\hspace{-3pt}\eps}}_\perp\, , \; \;
\frac{d\vvvec_{\hspace{-3pt}\eps}}{dt} = \eevec(t,\xxep)
+\frac 1\eps \vvvec_{\hspace{-3pt}\eps} \times \calm
\, , \; 
\xxep(s;\xvec,\vvec,s)=\xvec \, , \; \vvvec_{\hspace{-3pt}\eps}(s;\xvec,\vvec,s)=\vvec,
\end{equation}
is given  by
\begin{equation}
\xxze(t,\theta;\xvec,\vvec,s) = \yyze(t;\xvec,\vvec,s)
+\calr(\theta) \uuvec^0(t;\xvec,\vvec,s) \, ,\;\; 
\vvvec^0 (t,\theta;\xvec,\vvec,s) = R(\theta) \uuvec^0(t;\xvec,\vvec,s),
\end{equation}
where $(\yyze(t;\xvec,\vvec,s),\uuvec^0(t;\xvec,\vvec,s))$ is solution to
\begin{equation}
\frac{d\yyze}{dt} = \uuvec^0_\parallel +
\frac1{2\pi}\int_0^{2\pi}
\calr(-\theta) \evec(\yyze+\calr(\theta)\uuvec^0,t)\,d\theta 
\, , \;
\frac{d\uuvec^0}{dt} = \frac1{2\pi}\int_0^{2\pi}
R(-\theta)\evec(\yyze+\calr(\theta)\uuvec^0,t)\,d\theta,
\end{equation}
with the initial conditions
\begin{equation}
\yyze(s;\xvec,\vvec,s)= \xvec
\, , \;\;
\uuvec^0(s;\xvec,\vvec,s)= \vvec.
\end{equation}
\end{thm}

\subsection{Guiding Centre Regime with variable strong magnetic field}
Here we study 
a situation representative of what happens inside a tokamak, i.e.,
the Guiding Centre Regime with a variable $\calm$,
with $\caln=0$ and  $\bbvec=\basevectr$. In other words we consider
the following system:
\begin{equation}\label{ct1}
\frac{d\xxep}{dt} = \vvvec _\eps\, , \; \;
\frac{d\vvvec_\eps}{dt} = \eevec(t,\xxep)+\vvvec _\eps \times \basevectr 
+ \frac 1\eps \vvvec_\eps \times \calm(\xxep)
\, , \; \;
\xxep(s;\xvec,\vvec,s)=\xvec \, ,\; \; \vvvec_\eps(s;\xvec,\vvec,s)=\vvec,
\end{equation}
where \vspace{-10pt}
\begin{equation}\label{ct2}
\calm(\xvec)=\frac{1}{\sqrt{\cxun^2+\cxde^2}}
\begin{pmatrix} -\cxde \\ \cxun \\ 0 \end{pmatrix}
=\rho^T(\xvec) \basevecun,
\text{ with } 
\rho(\xvec)= \frac{1}{\sqrt{\cxun^2+\cxde^2}}
             \begin{pmatrix} -\cxde & \cxun & 0 \\
                             -\cxun & \cxde & 0 \\
                               0    &    0  & 1\end{pmatrix}\hspace{-4pt},
\end{equation}
where $\xvec =(\cxun,\cxde,\cxtr)$ in the frame 
($\basevecun$,$\basevecde$,$\basevectr$) of  $\rit^3$.
In this case
\begin{equation}\label{ct3}
\avec(t,\theta,\xvec,\vvec) =\avec(t,\xvec,\vvec) =
\begin{pmatrix} \vvec \\ \eevec(t,\xvec)+\vvec\times \basevectr \end{pmatrix}\hspace{-4pt},
\text{ and }
\bvec(t,\xvec,\vvec) =\bvec(\xvec,\vvec)=
\begin{pmatrix} 0 \\ \vvec\times\calm(\xvec) \end{pmatrix}\hspace{-4pt},
\end{equation}
and thus, $R(\theta)$ being defined by (\ref{defR}),
\begin{equation}\label{ct4}
\zzvec(t,\theta;\zvec,\wvec)\hspace{-2pt}=\hspace{-2pt}\begin{pmatrix} \zvec \\ 
\rho^T\hspace{-2pt}(\zvec)R(\theta)\rho(\zvec)\wvec\end{pmatrix}\hspace{-4pt}, 
\text{ where }
\rho^T R \rho =\hspace{-2pt}
\begin{pmatrix} 
\frac{\czde^2+\czun^2\cos(\theta)}{\czun^2+\czde^2} & 
\frac{\czun\czde(\cos(\theta)-1)}{\czun^2+\czde^2} & 
-\frac{\czun\sin(\theta)}{\sqrt{\czun^2+\czde^2}}\\
\frac{\czun\czde(\cos(\theta)-1)}{\czun^2+\czde^2} & 
\frac{\czun^2+\czde^2\cos(\theta)}{\czun^2+\czde^2} & 
-\frac{\czde\sin(\theta)}{\sqrt{\czun^2+\czde^2}}\\
\frac{\czun\sin(\theta)}{\sqrt{\czun^2+\czde^2}}& 
\frac{\czde\sin(\theta)}{\sqrt{\czun^2+\czde^2}} & 
\cos(\theta) \\
\end{pmatrix}\hspace{-4pt},
\end{equation}
and we have the following Theorem.
\begin{thm} 
If $\eevec(t,\xvec)$ is $C^2_b(\rit\times\rit^3)$,
the first term of the expansion (\ref{Rexp}) of the solution 
$(\xxep(t;\xvec,\vvec,s),$ $\vvvec_\eps(t;\xvec,\vvec,s))$ to (\ref{ct1})
is given by
\begin{equation}
\xxze(t,\theta;\xvec,\vvec,s) = \yyze(t;\xvec,\vvec,s) \, ,\;\; 
\vvvec^0 (t,\theta;\xvec,\vvec,s) = 
\rho^T(\xvec)R(\theta)\rho(\xvec) \uuvec^0(t;\xvec,\vvec,s),
\end{equation}
where $(\yyze(t;\xvec,\vvec,s),\uuvec^0(t;\xvec,\vvec,s))$ is solution to
\begin{equation}
\begin{array}{l}
\ds 
\frac{d\yyze}{dt} = \overline{A}(\yyze)\uuvec^0,~~~~
\frac{d\uuvec^0}{dt} = \overline{\beta}(\yyze,\uuvec^0)+
\overline{A}(\yyze)\eevec_\parallel(t,\yyze)+
\uuvec^0\times\overline{A}(\yyze)\basevectr,
\\
\text{~~~~~ with } \overline{A}(\yvec)=\ds\frac{1}{\cyun^2+\cyde^2}
\begin{pmatrix}
  \cyde^2 & -\cyun\cyde & 0 \\
  -\cyun\cyde & \cyun^2 & 0 \\
  0&0&0
\end{pmatrix},
\text{ and } 
\overline{\beta}(\yvec,\uvec) = 
\begin{pmatrix} \ds
  \frac{(\cyde\cuun-\cyun\cude)\cude}{\cyun^2+\cyde^2}\\ \ds
  \frac{(\cyun\cude-\cyde\cuun)\cuun }{\cyun^2+\cyde^2}\\
 0
\end{pmatrix},
\\
\yyze(s;\xvec,\vvec,s)= \xvec,
\, , \;\;
\uuvec^0(s;\xvec,\vvec,s)= \vvec.
\end{array}
\end{equation}
The term $\xxun(t,\theta;\xvec,\vvec,s)$ is given by
\begin{equation}
\begin{array}{l} \ds 
{\xxun}_1(t,\theta;\xvec,\vvec,s) \hspace{-2pt}= \hspace{-2pt}
\frac{1}{{\Omega}^2}\Big({{\yyze}_1}{\Omega} ({\cos}({\theta})-1) {{\uuvec^0}_3}  + 
{{\yyze}_1} ({{\yyze}_1} {{\uuvec^0}_1} + 
{{\yyze}_2} {{\uuvec^0}_2}) \sin({\theta})
+ {{\yyun}_1} {\Omega}^2  \Big)  ,
\\ \ds 
{\xxun}_2(t,\theta;\xvec,\vvec,s) =\hspace{+1pt} 
\frac{1}{{\Omega}^2}\Big( {{\yyze}_2}{\Omega} ({\cos}({\theta})-1) {{\uuvec^0}_3}  + 
{{\yyze}_2} ({{\yyze}_1} {{\uuvec^0}_1} + {{\yyze}_2} {{\uuvec^0}_2}) \sin({\theta})
 + {{\yyun}_2} {\Omega}^2  \Big) ,
\\ \ds 
{\xxun}_3(t,\theta;\xvec,\vvec,s) \hspace{-2pt}= \hspace{-2pt}
\frac{1}{{\Omega}}\Big( (-{{\yyze}_1} {{\uuvec^0}_1} - {{\yyze}_2} {{\uuvec^0}_2}) 
({\cos}({\theta})-1) + {{\yyun}_3} {\Omega}  +
\sin({\theta}) {{\uuvec^0}_3} {\Omega} \Big),
\end{array}
\end{equation}
where ${\Omega}=\sqrt{{{\yyze}_1}^2+{{\yyze}_2}^2}$ and 
where $\yyun(t;\xvec,\vvec,s)$ is solution to
\begin{equation*}
\begin{array}{l} \ds 
\frac{d{\yyun}_1}{dt}= 
\bigg(\Big(-{\Omega}^2{{\yyze}_2}{{\uuvec^0}_2}+
2{{\yyze}_1}{{\yyze}_2}(-{{\uuvec^0}_1}{{\yyze}_2}+
{{\uuvec^0}_2}{{\yyze}_1})\Big){{\yyun}_1}+
\Big((2{{\uuvec^0}_1}{{\yyze}_2}-
{{\uuvec^0}_2}{{\yyze}_1}){\Omega}^2+
\\ \ds ~~~~
2{{\yyze}_2}^2(-{{\uuvec^0}_1}{{\yyze}_2}+
{{\uuvec^0}_2}{{\yyze}_1})\Big){{\yyun}_2}+
{{\yyze}_2}^2{{\uuvec^1}_1}{\Omega}^2-
{{\yyze}_1}{{\yyze}_2}{\Omega}^2{{\uuvec^1}_2}-
{\Omega}^3{{\yyze}_1}{\evec_3}-
{{\yyze}_2}{{\uuvec^0}_3}{\Omega}^3-
\\ \ds ~\hfill
2{{\yyze}_2}{{\uuvec^0}_3}(-{{\uuvec^0}_1}{{\yyze}_2}+
{{\uuvec^0}_2}{{\yyze}_1}){\Omega}\bigg)\bigg/{\Omega}^4,
\end{array}
\end{equation*}
\vspace{-10pt}
\begin{equation*}
\begin{array}{l} \ds 
\\ \ds 
\frac{d{\yyun}_2}{dt}= 
\bigg(\Big((2{{\uuvec^0}_2}{{\yyze}_1}-
{{\uuvec^0}_1}{{\yyze}_2}){\Omega}^2-
2{{\yyze}_1}^2(-{{\uuvec^0}_1}{{\yyze}_2}+
{{\uuvec^0}_2}{{\yyze}_1})\Big){{\yyun}_1}+
\Big(-{\Omega}^2{{\yyze}_1}{{\uuvec^0}_1}-
\\ \ds ~~~~
2{{\yyze}_1}{{\yyze}_2}(-{{\uuvec^0}_1}{{\yyze}_2}+
{{\uuvec^0}_2}{{\yyze}_1})\Big){{\yyun}_2}-
{{\yyze}_1}{{\yyze}_2}{\Omega}^2{{\uuvec^1}_1}+
{{\yyze}_1}^2{{\uuvec^1}_2}{\Omega}^2- 
\\ \ds ~ \hfill
{\Omega}^3{{\yyze}_2}{\evec_3}+
{\Omega}^3{{\yyze}_1}{{\uuvec^0}_3}+
2{{\yyze}_1}{{\uuvec^0}_3}(-{{\uuvec^0}_1}{{\yyze}_2}+
{{\uuvec^0}_2}{{\yyze}_1}){\Omega}\bigg)\bigg/{\Omega}^4,
\end{array}
\end{equation*}
\vspace{-10pt}
\begin{equation}
\begin{array}{l} \ds 
\\ \ds 
\frac{d{\yyun}_3}{dt}=
\bigg({\Omega}^2{{\yyze}_1}{\evec_1}+
{\Omega}^2{{\yyze}_2}{\evec_2}+
\Big({{\uuvec^0}_2}{{\yyze}_1}+{{\uuvec^0}_1}^2-
{{\uuvec^0}_1}{{\yyze}_2}+
{{\uuvec^0}_2}^2\Big){\Omega}^2-~~~~~~~~~~~~~~~~
\\ \ds \hfill
\Big({{\yyze}_1}{{\uuvec^0}_1}+
{{\yyze}_2}{{\uuvec^0}_2}\Big)^2\bigg)\bigg/{\Omega}^3,
\\ \ds 
\yyun(s;\xvec,\vvec,s)= 0,
\end{array}
\end{equation}
with  $\evec$ being evaluated in ${\yyze}$.\;
The result concerning $\vvvec^1(t,\theta;\xvec,\vvec,s)$ is given in 
the appendix \ref{appB}.
\end{thm}
This Theorem is the consequence of the following computations. 
Setting $A=\rho^T R \rho$ we have 
\begin{equation}
\big\{\Nabla_{z,w} \zzvec(t,\theta;\zvec,\wvec)\big\}=
\begin{pmatrix} I&0\\\Nabla_{z}(A\wvec)& A\end{pmatrix} \text{ and }
\big\{\Nabla_{z,w} \zzvec(t,\theta;\zvec,\wvec)\big\}^{-1}=
\begin{pmatrix} I&0\\A^T \Nabla_{z}(A\wvec)&A^T\end{pmatrix}.
\end{equation}
Then, since in this example $\fracp{\zzvec}{t}=0,$ and 
\begin{equation}
\avec(t,\zzvec(t,\theta;\zvec,\wvec))= \begin{pmatrix}
                                  A\wvec\\ \evec + A \wvec\times \basevectr
				       \end{pmatrix},
\end{equation}
applying (\ref{a0tilde}) and (\ref{al0}) we get \begin{equation}
\tilde \alpha^0(t,\theta,\zvec,\wvec)= \begin{pmatrix}
         A\wvec\\ -A^T {\Nabla_{z}(A\wvec)}(A\wvec)+A^T\evec +  \wvec\times A^T\basevectr
				       \end{pmatrix},
\end{equation}
and since $\int_0^{2\pi}A d\theta =\int_0^{2\pi}A^T d\theta = 2\pi \overline{A}$,
and $\int_0^{2\pi}-A^T {\Nabla_{z}(A\wvec)}(A\wvec) = \overline{\beta}$,
we finally get the first part of the Theorem.
 
The second part of the Theorem is obtained following the computation
program described  in the previous sections using Maple. 
\begin{appendix}
\section{Appendix : $\vvvec^1$ for the Guiding Centre Regime with variable strong magnetic field}
\label{appB}
The computations of this appendix have been realized using Maple5.5, on computers
of Medicis Centre, Ecole Polytechnique, Palaiseau, France.

~

The electric field $\evec$ and its derivatives being evaluated in ${\yyze}$,
and ${\Omega}$ standing for $\sqrt{({{\yyze}_1})^2+({{\yyze}_2})^2}$,
the term $\vvvec^1(t,\theta;\xvec,\vvec,s)$ of the expansion (\ref{Rexp}) 
of the solution to (\ref{ct1}) is given by~:
{\small
\begin{multline}
{\vvvec^1}_1(t;\xvec,\vvec,s)=
\Bigg(\bigg(\Big((2 {{\yyze}_1} {{\uuvec^0}_1} + {{\yyze}_2} {{\uuvec^0}_2}) {\Omega}^2  - 
2 {{\yyze}_1}^2  ({{\yyze}_1} {{\uuvec^0}_1} + {{\yyze}_2} {{\uuvec^0}_2})\Big)  {{\yyun}_1} +
\displaybreak[0] \\[-5pt]
\Big({{\yyze}_1} {{\uuvec^0}_2} {\Omega}^2  
- 2 {{\yyze}_1} {{\yyze}_2} ({{\yyze}_1} {{\uuvec^0}_1} + {{\yyze}_2} {{\uuvec^0}_2})\Big) {{\yyun}_2}+ 
{{\yyze}_1}^2  {{\uuvec^1}_1} {\Omega}^2  + {{\yyze}_1} {{\yyze}_2} {\Omega}^2  {{\uuvec^1}_2} + 
{\Omega}^3  {{\yyze}_1} {{\evec}_3} +
\displaybreak[0] \\[-5pt]
{\Omega}^3  {{\yyze}_2} {{\uuvec^0}_3} + 
{{\yyze}_2} {{\uuvec^0}_3} (-{{\uuvec^0}_1} {{\yyze}_2} + 
{{\yyze}_1} {{\uuvec^0}_2}) {\Omega}\bigg) {\cos}({\theta}) + 
\Big( (-{\Omega}^3  {{\uuvec^0}_3} + {\Omega} {{\uuvec^0}_3} {{\yyze}_1}^2 ) {{\yyun}_1} + 
\displaybreak[0] \\[-5pt]
{{\yyze}_1} {{\uuvec^0}_3} {{\yyze}_2} {\Omega} {{\yyun}_2} - 
{{\yyze}_1} {\Omega}^3  {{\uuvec^1}_3} + 
{\Omega}^2   {{\yyze}_1}^2   {{\evec}_1} + 
{\Omega}^2  {{\yyze}_1} {{\evec}_2} {{\yyze}_2} + 
{{\uuvec^0}_2} {\Omega}^4 + 
\displaybreak[0] \\[-2pt]
{{\uuvec^0}_2} (-{{\uuvec^0}_1} {{\yyze}_2} + 
{{\yyze}_1} {{\uuvec^0}_2}) {\Omega}^2 \Big) {\sin}({\theta}) + 
\Big(-{\Omega}^2  {{\yyze}_2} {{\uuvec^0}_2} + 
2 {{\yyze}_1} {{\yyze}_2} (-{{\uuvec^0}_1} {{\yyze}_2} + {{\yyze}_1} {{\uuvec^0}_2})\Big) {{\yyun}_1} +
\displaybreak[0] \\[-2pt]
 \Big((2 {{\uuvec^0}_1} {{\yyze}_2} - {{\yyze}_1} {{\uuvec^0}_2}) {\Omega}^2  + 
2 {{\yyze}_2}^2  (-{{\uuvec^0}_1} {{\yyze}_2} + {{\yyze}_1} {{\uuvec^0}_2})\Big) {{\yyun}_2}+
{{\yyze}_2}^2  {{\uuvec^1}_1} {\Omega}^2  - {{\yyze}_1} {{\yyze}_2} {\Omega}^2  {{\uuvec^1}_2} -  
\displaybreak[0] \\[-5pt]
{\Omega}^3  {{\yyze}_1} {{\evec}_3} - 
{\Omega}^3  {{\yyze}_2} {{\uuvec^0}_3} - 
{{\yyze}_2} {{\uuvec^0}_3} (-{{\uuvec^0}_1} {{\yyze}_2} + 
{{\yyze}_1} {{\uuvec^0}_2}) {\Omega}\Bigg) \bigg /  {\Omega}^4  ,
\end{multline}
}\vspace{-20pt}
{\small
\begin{multline}
\tiny
{\vvvec^1}_2(t;\xvec,\vvec,s)= 
\Bigg(\bigg( ({{\uuvec^0}_1} {{\yyze}_2} {\Omega}^2   - 
2 {{\yyze}_1} {{\yyze}_2} ({{\yyze}_1} {{\uuvec^0}_1} + {{\yyze}_2} {{\uuvec^0}_2})) {{\yyun}_1} +
\bigg(({{\yyze}_1} {{\uuvec^0}_1} + 2 {{\yyze}_2} {{\uuvec^0}_2}) {\Omega}^2  - 
\displaybreak[0] \\[-5pt]
2 {{\yyze}_2}^2  ({{\yyze}_1} {{\uuvec^0}_1} + {{\yyze}_2} {{\uuvec^0}_2})\bigg) {{\yyun}_2} +
{{\yyze}_1} {{\yyze}_2} {\Omega}^2  {{\uuvec^1}_1} + {{\yyze}_2}^2  {{\uuvec^1}_2} {\Omega}^2 + 
{\Omega}^3  {{\yyze}_2} {{\evec}_3} - 
{\Omega}^3  {{\yyze}_1} {{\uuvec^0}_3} - 
\displaybreak[0] \\[-5pt]
{{\yyze}_1} {{\uuvec^0}_3} (-{{\uuvec^0}_1} {{\yyze}_2} + 
{{\yyze}_1} {{\uuvec^0}_2}) {\Omega}\bigg) {\cos}({\theta}) + 
\bigg( {{\yyze}_1} {{\uuvec^0}_3} {{\yyze}_2} {\Omega} {{\yyun}_1} + (-{\Omega}^3  {{\uuvec^0}_3} + 
{\Omega} {{\uuvec^0}_3} {{\yyze}_2}^2 ) {{\yyun}_2} - 
\displaybreak[0] \\[-5pt]
{{\yyze}_2} {\Omega}^3  {{\uuvec^1}_3} + 
{\Omega}^2  {{\yyze}_1} {{\yyze}_2} {{\evec}_1} + 
{\Omega}^2  {{\yyze}_2}^2  {{\evec}_2} - 
{{\uuvec^0}_1} {\Omega}^4 - 
{{\uuvec^0}_1} (-{{\uuvec^0}_1} {{\yyze}_2} + 
{{\yyze}_1} {{\uuvec^0}_2}) {\Omega}^2 \bigg) {\sin}({\theta}) +
\displaybreak[0] \\[-5pt]
 \bigg((2 {{\yyze}_1} {{\uuvec^0}_2} - {{\uuvec^0}_1} {{\yyze}_2}) {\Omega}^2  - 
2 {{\yyze}_1}^2  (-{{\uuvec^0}_1} {{\yyze}_2} + 
{{\yyze}_1} {{\uuvec^0}_2})\bigg) {{\yyun}_1} +
\displaybreak[0] \\[-5pt] 
(-{\Omega}^2   {{\yyze}_1} {{\uuvec^0}_1} - 
2 {{\yyze}_1} {{\yyze}_2} (-{{\uuvec^0}_1} {{\yyze}_2} + {{\yyze}_1} {{\uuvec^0}_2})) {{\yyun}_2} -
{{\yyze}_1} {{\yyze}_2} {\Omega}^2   {{\uuvec^1}_1} + {{\yyze}_1}^2   {{\uuvec^1}_2} {\Omega}^2  - 
\displaybreak[0] \\[-5pt]
{\Omega}^3  {{\yyze}_2} {{\evec}_3} + 
{\Omega}^3  {{\yyze}_1} {{\uuvec^0}_3} + 
{{\yyze}_1} {{\uuvec^0}_3} (-{{\uuvec^0}_1} {{\yyze}_2} + 
{{\yyze}_1} {{\uuvec^0}_2}) {\Omega}\Bigg)  \bigg/  {\Omega}^4 ,
\end{multline}
}\vspace{-20pt}
{\small
\begin{multline}
{\vvvec^1}_3(t;\xvec,\vvec,s)= 
\Bigg(\bigg({\Omega}^3  {{\uuvec^1}_3} - {\Omega}^2  {{\yyze}_1} {{\evec}_1} - 
{\Omega}^2   {{\yyze}_2} {{\evec}_2} + 
({{\uuvec^0}_1} {{\yyze}_2} - {{\yyze}_1} {{\uuvec^0}_2}) {\Omega}^2  - 
\displaybreak[0] \\[-5pt]
(-{{\uuvec^0}_1} {{\yyze}_2} 
+ {{\yyze}_1} {{\uuvec^0}_2})^2 \bigg) {\cos}({\theta}) + 
\bigg(({{\uuvec^0}_1} {\Omega}^2  - {{\yyze}_1} ({{\yyze}_1} {{\uuvec^0}_1} + 
{{\yyze}_2} {{\uuvec^0}_2})) {{\yyun}_1} + 
\displaybreak[0] \\[-5pt]
({{\uuvec^0}_2} {\Omega}^2  - {{\yyze}_2} ({{\yyze}_1} {{\uuvec^0}_1} + 
{{\yyze}_2} {{\uuvec^0}_2})) {{\yyun}_2} + 
{{\yyze}_2} {\Omega}^2 {{\uuvec^1}_2} + {\Omega}^3  {{\evec}_3} + 
{{\yyze}_1} {\Omega}^2 {{\uuvec^1}_1}\bigg) {\sin}({\theta}) + 
\displaybreak[0] \\[-5pt]
{\Omega}^2  {{\yyze}_1} {{\evec}_1} + 
{\Omega}^2  {{\yyze}_2} {{\evec}_2} + 
({{\uuvec^0}_1}^2   + {{\yyze}_1} {{\uuvec^0}_2} - {{\uuvec^0}_1} {{\yyze}_2} + 
{{\uuvec^0}_2}^2  ) {\Omega}^2  - 
({{\yyze}_1} {{\uuvec^0}_1} + {{\yyze}_2} {{\uuvec^0}_2})^2 \Bigg)  \bigg/  {\Omega}^3 ,
\end{multline}
}
where $\uuvec^1(t;\xvec,\vvec,s))$ is solution to:
{\scriptsize

$\ds
\frac{d{\uuvec^1}_1}{dt}$ $=
$
$
1/4$ $\bigg(\Big(-4$ ${{\yyze}_2}$ ${\Omega}^8$ ${{\yyze}_1}$ $\fracp{\evec_2}{x_3}$ $+
4{{\yyze}_2}^2$ ${\Omega}^8$ $\fracp{\evec_1}{x_3}$ $\Big)$ ${{\yyun}_3}+\Big(-4{\Omega}^9$ ${{\yyze}_2}+
(-{{\yyze}_1}^2$ ${{\uuvec^0}_3}+4{{\yyze}_2}$ ${{\uuvec^0}_2}$ ${{\yyze}_1}-4{{\uuvec^0}_1}$ ${{\yyze}_2}^2)$ ${\Omega}^7+
(-{{\yyze}_1}^4$ ${{\uuvec^0}_1}-2{{\yyze}_1}$ ${{\yyze}_2}^3$ ${{\uuvec^0}_1}+2$ ${{\yyze}_2}^4$ ${{\uuvec^0}_1}+
2{{\yyze}_1}^2$ ${{\yyze}_2}^2$ ${{\uuvec^0}_2}+{{\yyze}_2}^2$ ${{\yyze}_1}^2$ ${{\uuvec^0}_3}-
{{\uuvec^0}_2}$ ${{\yyze}_2}$ ${{\yyze}_1}^3-2$ ${{\yyze}_1}^2$ ${{\uuvec^0}_1}$ ${{\yyze}_2}^2+
{{\yyze}_1}^3$ ${{\uuvec^0}_3}$ ${{\yyze}_2}-4$ ${{\uuvec^0}_2}$ ${{\yyze}_2}^3$ ${{\yyze}_1})$ ${\Omega}^5+
{{\yyze}_1}$ $(3$ ${{\yyze}_1}^2$ ${{\yyze}_2}^2+{{\yyze}_2}^4+{{\yyze}_1}^4+
{{\yyze}_2}$ ${{\yyze}_1}^3)$ $({{\yyze}_1}$ ${{\uuvec^0}_1}+{{\yyze}_2}$ ${{\uuvec^0}_2})$ ${\Omega}^3\Big)$ ${\evec_3}$ $+
\Big(4$ ${\Omega}^8$ ${{\yyze}_2}$ ${{\uuvec^0}_2}-{{\yyze}_1}^3$ ${\Omega}^6$ ${{\uuvec^0}_3}+
{{\yyze}_2}$ ${{\yyze}_1}^2$ $(-{{\yyze}_2}^2$ ${{\uuvec^0}_2}-{{\uuvec^0}_1}$ ${{\yyze}_2}$ ${{\yyze}_1}+
{{\yyze}_1}^2$ ${{\uuvec^0}_3}+{{\uuvec^0}_3}$ ${{\yyze}_1}$ ${{\yyze}_2})$ ${\Omega}^4+
{{\yyze}_1}^4$ ${{\yyze}_2}$ $({{\yyze}_2}+{{\yyze}_1})$ $({{\yyze}_1}$ ${{\uuvec^0}_1}+
{{\yyze}_2}$ ${{\uuvec^0}_2})$ ${\Omega}^2\Big)$ ${{\uuvec^1}_1}+\Big(-{\Omega}^7$ ${{\yyze}_1}$ ${{\uuvec^0}_3}+
{{\yyze}_1}$ $(8$ ${{\yyze}_1}^2$ ${{\uuvec^0}_3}-{{\yyze}_2}$ ${{\uuvec^0}_2}$ ${{\yyze}_1}-
{{\uuvec^0}_1}$ ${{\yyze}_2}$ ${{\yyze}_1}+5$ ${{\uuvec^0}_3}$ ${{\yyze}_1}$ ${{\yyze}_2}+
2$ ${{\yyze}_2}^2$ ${{\uuvec^0}_3})$ ${\Omega}^5-{{\yyze}_1}^2$ $(8$ ${{\yyze}_1}$ ${{\uuvec^0}_3}$ ${{\yyze}_2}^2-
{{\yyze}_1}$ ${{\yyze}_2}^2$ ${{\uuvec^0}_1}+7$ ${{\yyze}_1}^2$ ${{\yyze}_2}$ ${{\uuvec^0}_3}-
{{\yyze}_1}$ ${{\yyze}_2}^2$ ${{\uuvec^0}_2}+6$ ${{\yyze}_2}^3$ ${{\uuvec^0}_3}+7$ ${{\yyze}_1}^3$ ${{\uuvec^0}_3}-
{{\yyze}_1}^2$ ${{\yyze}_2}$ ${{\uuvec^0}_1}-{{\yyze}_2}^3$ ${{\uuvec^0}_2})$ ${\Omega}^3\Big)$ ${{\yyun}_1}^2+
\Big({{\yyze}_1}^2$ ${{\yyze}_2}^5$ ${{\uuvec^0}_2}^3-
8$ ${{\uuvec^0}_1}^2$ ${{\uuvec^0}_3}$ ${{\yyze}_2}^7+7$ ${{\uuvec^0}_1}^2$ ${{\uuvec^0}_3}$ ${{\yyze}_1}^7+
{{\yyze}_2}^2$ ${{\yyze}_1}^5$ ${{\uuvec^0}_1}^3+
6$ ${{\uuvec^0}_3}$ ${{\yyze}_1}^7$ ${{\uuvec^0}_2}^2+
{{\yyze}_1}^6$ ${{\uuvec^0}_1}^3$ ${{\yyze}_2}+{{\yyze}_2}^4$ ${{\yyze}_1}^3$ ${{\uuvec^0}_2}^3+
5$ ${{\uuvec^0}_2}^2$ ${{\yyze}_1}^5$ ${{\yyze}_2}^2$ ${{\uuvec^0}_3}+
{{\uuvec^0}_3}$ ${{\uuvec^0}_2}$ ${{\uuvec^0}_1}$ ${{\yyze}_1}^7-
5$ ${{\yyze}_1}^3$ ${{\yyze}_2}^4$ ${{\uuvec^0}_3}$ ${{\uuvec^0}_2}^2+
3$ ${{\yyze}_1}^5$ ${{\yyze}_2}^2$ ${{\uuvec^0}_2}$ ${{\uuvec^0}_1}^2+
19$ ${{\yyze}_1}^4$ ${{\uuvec^0}_3}$ ${{\uuvec^0}_2}^2$ ${{\yyze}_2}^3+
9$ ${{\uuvec^0}_1}^2$ ${{\uuvec^0}_3}$ ${{\yyze}_1}^6$ ${{\yyze}_2}-
9$ ${{\uuvec^0}_1}^2$ ${{\uuvec^0}_3}$ ${{\yyze}_2}^6$ ${{\yyze}_1}+
{{\uuvec^0}_1}^2$ ${{\uuvec^0}_3}$ ${{\yyze}_1}^4$ ${{\yyze}_2}^3-
16$ ${{\uuvec^0}_1}^2$ ${{\uuvec^0}_3}$ ${{\yyze}_2}^4$ ${{\yyze}_1}^3+
3$ ${{\yyze}_2}^4$ ${{\yyze}_1}^3$ ${{\uuvec^0}_2}^2$ ${{\uuvec^0}_1}+
{{\yyze}_1}^5$ ${{\yyze}_2}^2$ ${{\uuvec^0}_3}$ ${{\uuvec^0}_1}^2+
3$ ${{\yyze}_2}^3$ ${{\uuvec^0}_1}$ ${{\yyze}_1}^4$ ${{\uuvec^0}_2}^2-
15$ ${{\yyze}_2}^5$ ${{\yyze}_1}^2$ ${{\uuvec^0}_1}^2$ ${{\uuvec^0}_3}+
3$ ${{\yyze}_2}^3$ ${{\yyze}_1}^4$ ${{\uuvec^0}_2}$ ${{\uuvec^0}_1}^2-
5$ ${{\yyze}_2}^6$ ${{\uuvec^0}_3}$ ${{\uuvec^0}_2}^2$ ${{\yyze}_1}+
10$ ${{\yyze}_2}^5$ ${{\uuvec^0}_3}$ ${{\uuvec^0}_2}^2$ ${{\yyze}_1}^2+
8$ ${{\yyze}_2}$ ${{\uuvec^0}_3}$ ${{\uuvec^0}_2}^2$ ${{\yyze}_1}^6-
6$ ${{\yyze}_1}^2$ ${{\yyze}_2}^5$ ${{\uuvec^0}_3}$ ${{\uuvec^0}_2}$ ${{\uuvec^0}_1}+
19$ ${{\yyze}_1}^5$ ${{\yyze}_2}^2$ ${{\uuvec^0}_3}$ ${{\uuvec^0}_2}$ ${{\uuvec^0}_1}+
7$ ${{\yyze}_1}^6$ ${{\yyze}_2}$ ${{\uuvec^0}_3}$ ${{\uuvec^0}_2}$ ${{\uuvec^0}_1}+
3$ ${{\yyze}_1}^4$ ${{\yyze}_2}^3$ ${{\uuvec^0}_3}$ ${{\uuvec^0}_1}$ ${{\uuvec^0}_2}+
16$ ${{\yyze}_1}^3$ ${{\yyze}_2}^4$ ${{\uuvec^0}_3}$ ${{\uuvec^0}_1}$ ${{\uuvec^0}_2}\Big)$ ${\Omega}+
(-4$ ${\Omega}^8$ ${{\uuvec^0}_2}$ ${{\yyze}_1}-{{\yyze}_1}^2$ ${\Omega}^6$ ${{\yyze}_2}$ ${{\uuvec^0}_3}+
{{\yyze}_1}$ $\Big({{\yyze}_1}$ ${{\yyze}_2}^3$ ${{\uuvec^0}_3}-{{\uuvec^0}_2}$ ${{\yyze}_2}^4+
{{\uuvec^0}_2}$ ${{\yyze}_2}^3$ ${{\yyze}_1}-{{\yyze}_1}$ ${{\yyze}_2}^3$ ${{\uuvec^0}_1}+
{{\yyze}_1}^2$ ${{\uuvec^0}_1}$ ${{\yyze}_2}^2+{{\yyze}_2}^2$ ${{\yyze}_1}^2$ ${{\uuvec^0}_3}+
{{\uuvec^0}_2}$ ${{\yyze}_2}$ ${{\yyze}_1}^3+{{\yyze}_1}^4$ ${{\uuvec^0}_1}\Big)$ ${\Omega}^4-
{{\yyze}_1}^5$ $({{\yyze}_2}+
{{\yyze}_1})$ $({{\yyze}_1}$ ${{\uuvec^0}_1}+{{\yyze}_2}$ ${{\uuvec^0}_2})$ ${\Omega}^2)$ ${{\uuvec^1}_2}+
\Big(6$ ${\Omega}^6$ ${{\yyze}_1}^2$ ${{\yyze}_2}$ ${\evec_1}$ $-2$ ${{\yyze}_1}$ $(-{{\yyze}_2}^2+
2$ ${{\yyze}_1}^2)$ ${\Omega}^6$ ${\evec_2}$ $+4$ ${\Omega}^8$ ${{\yyze}_2}^2$ $\fracp{\evec_1}{x_2}$ $-
4$ ${\Omega}^8$ ${{\yyze}_1}$ $\fracp{\evec_2}{x_2}$ ${{\yyze}_2}-{\Omega}^8$ ${{\yyze}_1}$ ${{\uuvec^0}_1}+
{{\yyze}_1}$ $(-{{\yyze}_1}$ ${{\uuvec^0}_3}$ ${{\uuvec^0}_2}+3$ ${{\yyze}_2}$ ${{\uuvec^0}_3}^2+
{{\uuvec^0}_1}$ ${{\yyze}_1}$ ${{\uuvec^0}_2}+12$ ${{\yyze}_2}$ ${{\uuvec^0}_2}^2+2$ ${{\yyze}_1}^2$ ${{\uuvec^0}_1}-
{{\yyze}_2}$ ${{\uuvec^0}_2}$ ${{\yyze}_1}+4$ ${{\yyze}_2}$ ${{\uuvec^0}_1}^2+3$ ${{\uuvec^0}_1}$ ${{\yyze}_2}^2+
{{\yyze}_1}$ ${{\uuvec^0}_3}^2)$ ${\Omega}^6+(-{{\yyze}_1}^2$ ${{\yyze}_2}^2$ ${{\uuvec^0}_2}^2+
2$ ${{\yyze}_1}^2$ ${{\yyze}_2}^2$ ${{\uuvec^0}_3}$ ${{\uuvec^0}_2}+
10$ ${{\yyze}_1}$ ${{\yyze}_2}^3$ ${{\uuvec^0}_1}$ ${{\uuvec^0}_2}+
{{\uuvec^0}_3}$ ${{\uuvec^0}_1}$ ${{\yyze}_2}$ ${{\yyze}_1}^3-3$ ${{\yyze}_1}$ ${{\yyze}_2}^3$ ${{\uuvec^0}_2}^2-
8$ ${{\yyze}_2}$ ${{\uuvec^0}_3}^2$ ${{\yyze}_1}^3-{{\yyze}_1}^4$ ${{\uuvec^0}_2}$ ${{\uuvec^0}_1}-
3$ ${{\uuvec^0}_3}^2$ ${{\yyze}_1}^4+{{\uuvec^0}_2}$ ${{\yyze}_1}^3$ ${{\yyze}_2}$ ${{\uuvec^0}_3}+
2$ ${{\yyze}_1}^2$ ${{\yyze}_2}^2$ ${{\uuvec^0}_1}^2-9$ ${{\yyze}_2}^2$ ${{\uuvec^0}_3}^2$ ${{\yyze}_1}^2-
2$ ${{\yyze}_1}^3$ ${{\uuvec^0}_1}$ ${{\yyze}_2}$ ${{\uuvec^0}_2}+2$ ${{\yyze}_1}^3$ ${{\yyze}_2}$ ${{\uuvec^0}_1}^2-
4$ ${{\uuvec^0}_1}$ ${{\yyze}_2}^4$ ${{\uuvec^0}_2}-6$ ${{\uuvec^0}_1}^2$ ${{\yyze}_1}$ ${{\yyze}_2}^3-
8$ ${{\uuvec^0}_1}^2$ ${{\yyze}_2}^4-{{\yyze}_2}^3$ ${{\uuvec^0}_3}^2$ ${{\yyze}_1})$ ${\Omega}^4-
{{\yyze}_2}$ $(-8$ ${{\uuvec^0}_1}^2$ ${{\yyze}_2}^5-2$ ${{\yyze}_1}$ ${{\yyze}_2}^4$ ${{\uuvec^0}_2}^2-
7$ ${{\yyze}_2}$ ${{\uuvec^0}_3}^2$ ${{\yyze}_1}^4-{{\yyze}_2}^4$ ${{\uuvec^0}_3}^2$ ${{\yyze}_1}+
9$ ${{\yyze}_1}^2$ ${{\yyze}_2}^3$ ${{\uuvec^0}_2}^2-2$ ${{\yyze}_1}^3$ ${{\uuvec^0}_1}^2$ ${{\yyze}_2}^2-
4$ ${{\yyze}_1}^4$ ${{\yyze}_2}$ ${{\uuvec^0}_1}^2-9$ ${{\yyze}_1}^4$ ${{\yyze}_2}$ ${{\uuvec^0}_2}$ ${{\uuvec^0}_1}-
14$ ${{\yyze}_2}^3$ ${{\yyze}_1}^2$ ${{\uuvec^0}_1}^2-9$ ${{\uuvec^0}_3}^2$ ${{\yyze}_1}^3$ ${{\yyze}_2}^2-
7$ ${{\yyze}_1}^5$ ${{\uuvec^0}_3}^2+11$ ${{\yyze}_1}^3$ ${{\yyze}_2}^2$ ${{\uuvec^0}_2}^2+
{{\yyze}_2}$ ${{\uuvec^0}_3}$ ${{\yyze}_1}^4$ ${{\uuvec^0}_1}+{{\yyze}_2}^2$ ${{\uuvec^0}_3}$ ${{\yyze}_1}^3$ ${{\uuvec^0}_1}+
7$ ${{\uuvec^0}_1}$ ${{\uuvec^0}_2}$ ${{\yyze}_1}^5-14$ ${{\yyze}_1}^2$ ${{\uuvec^0}_1}$ ${{\yyze}_2}^3$ ${{\uuvec^0}_2}+
{{\uuvec^0}_2}$ ${{\uuvec^0}_3}$ ${{\yyze}_2}^2$ ${{\yyze}_1}^3+12$ ${{\yyze}_1}^3$ ${{\uuvec^0}_1}$ ${{\yyze}_2}^2$ ${{\uuvec^0}_2}+
2$ ${{\yyze}_1}$ ${{\yyze}_2}^4$ ${{\uuvec^0}_2}$ ${{\uuvec^0}_1}+{{\uuvec^0}_2}$ ${{\yyze}_2}^3$ ${{\uuvec^0}_3}$ ${{\yyze}_1}^2+
10$ ${{\uuvec^0}_2}^2$ ${{\yyze}_1}^5+8$ ${{\yyze}_1}^4$ ${{\yyze}_2}$ ${{\uuvec^0}_2}^2-
6$ ${{\uuvec^0}_3}^2$ ${{\yyze}_2}^3$ ${{\yyze}_1}^2)$ ${\Omega}^2+2$ ${{\yyze}_1}^4$ ${{\yyze}_2}$ $({{\yyze}_2}+
{{\yyze}_1})$ $(-{{\uuvec^0}_1}$ ${{\yyze}_2}+{{\uuvec^0}_2}$ ${{\yyze}_1})$ $({{\yyze}_1}$ ${{\uuvec^0}_1}+
{{\yyze}_2}$ ${{\uuvec^0}_2})\Big)$ ${{\yyun}_2}+\Big(-6$ ${{\yyze}_2}^5$ ${{\uuvec^0}_1}^3+
{{\yyze}_1}^4$ ${{\uuvec^0}_1}^2$ ${{\yyze}_2}^2-
2$ ${{\uuvec^0}_3}$ ${{\uuvec^0}_2}^2$ ${{\yyze}_1}^5+5$ ${{\yyze}_2}^2$ ${{\yyze}_1}^3$ ${{\uuvec^0}_2}^3-
5$ ${{\uuvec^0}_3}$ ${{\uuvec^0}_1}^2$ ${{\yyze}_1}^5-{{\uuvec^0}_1}$ ${{\yyze}_1}^6$ ${{\uuvec^0}_3}+
8$ ${{\yyze}_2}^5$ ${{\uuvec^0}_1}^2$ ${{\uuvec^0}_3}-{{\yyze}_1}^4$ ${{\yyze}_2}^2$ ${{\uuvec^0}_2}^2+
{{\yyze}_1}^3$ ${{\yyze}_2}^3$ ${{\uuvec^0}_1}^2-{{\yyze}_1}^3$ ${{\yyze}_2}^3$ ${{\uuvec^0}_2}^2-
4$ ${{\yyze}_2}$ ${{\uuvec^0}_2}^3$ ${{\yyze}_1}^4+6$ ${{\yyze}_2}^4$ ${{\uuvec^0}_1}^3$ ${{\yyze}_1}-
{{\yyze}_2}$ ${{\yyze}_1}^4$ ${{\uuvec^0}_1}^3-6$ ${{\yyze}_1}$ ${{\uuvec^0}_3}$ ${{\yyze}_2}^4$ ${{\uuvec^0}_2}^2+
12$ ${{\yyze}_1}^3$ ${{\uuvec^0}_3}$ ${{\yyze}_2}^2$ ${{\uuvec^0}_1}^2-
5$ ${{\yyze}_1}^4$ ${{\uuvec^0}_3}$ ${{\yyze}_2}$ ${{\uuvec^0}_1}^2-{{\yyze}_1}^2$ ${{\uuvec^0}_3}$ ${{\yyze}_2}^3$ ${{\uuvec^0}_2}^2-
{{\yyze}_1}^4$ ${{\uuvec^0}_1}$ ${{\yyze}_2}^2$ ${{\uuvec^0}_2}-19$ ${{\yyze}_1}^2$ ${{\uuvec^0}_1}^2$ ${{\yyze}_2}^3$ ${{\uuvec^0}_2}+
{{\yyze}_1}^3$ ${{\uuvec^0}_1}$ ${{\yyze}_2}^3$ ${{\uuvec^0}_2}+
4$ ${{\uuvec^0}_3}$ ${{\uuvec^0}_2}$ ${{\uuvec^0}_1}$ ${{\yyze}_2}^5+
14$ ${{\uuvec^0}_1}$ ${{\yyze}_1}^3$ ${{\yyze}_2}^2$ ${{\uuvec^0}_2}^2-
{{\yyze}_1}^3$ ${{\yyze}_2}^2$ ${{\uuvec^0}_2}$ ${{\uuvec^0}_1}^2+
{{\yyze}_1}^2$ ${{\uuvec^0}_1}$ ${{\yyze}_2}^4$ ${{\uuvec^0}_2}-
{{\yyze}_1}^5$ ${{\yyze}_2}$ ${{\uuvec^0}_2}$ ${{\uuvec^0}_1}-
6$ ${{\yyze}_1}^3$ ${{\yyze}_2}^2$ ${{\uuvec^0}_3}$ ${{\uuvec^0}_2}^2+
9$ ${{\uuvec^0}_1}^2$ ${{\yyze}_1}^2$ ${{\yyze}_2}^3$ ${{\uuvec^0}_3}+
12$ ${{\uuvec^0}_1}^2$ ${{\yyze}_2}^4$ ${{\uuvec^0}_3}$ ${{\yyze}_1}-
18$ ${{\yyze}_1}^2$ ${{\uuvec^0}_2}^2$ ${{\uuvec^0}_1}$ ${{\yyze}_2}^3+
18$ ${{\yyze}_1}$ ${{\yyze}_2}^4$ ${{\uuvec^0}_2}$ ${{\uuvec^0}_1}^2+
{{\yyze}_1}^5$ ${{\yyze}_2}$ ${{\uuvec^0}_2}$ ${{\uuvec^0}_3}-
5$ ${{\yyze}_1}^5$ ${{\yyze}_2}$ ${{\uuvec^0}_1}$ ${{\uuvec^0}_3}-
6$ ${{\yyze}_1}^3$ ${{\yyze}_2}^3$ ${{\uuvec^0}_1}$ ${{\uuvec^0}_3}-
2$ ${{\yyze}_2}^5$ ${{\yyze}_1}$ ${{\uuvec^0}_1}$ ${{\uuvec^0}_3}-
2$ ${{\yyze}_1}^2$ ${{\yyze}_2}^4$ ${{\uuvec^0}_1}$ ${{\uuvec^0}_3}+
5$ ${{\yyze}_1}^4$ ${{\yyze}_2}^2$ ${{\uuvec^0}_2}$ ${{\uuvec^0}_3}+
2$ ${{\yyze}_1}^2$ ${{\yyze}_2}^4$ ${{\uuvec^0}_2}$ ${{\uuvec^0}_3}-
2$ ${{\yyze}_1}^4$ ${{\yyze}_2}^2$ ${{\uuvec^0}_1}$ ${{\uuvec^0}_3}+
2$ ${{\yyze}_1}^3$ ${{\uuvec^0}_2}$ ${{\yyze}_2}^3$ ${{\uuvec^0}_3}-
3$ ${{\yyze}_2}$ ${{\yyze}_1}^4$ ${{\uuvec^0}_2}$ ${{\uuvec^0}_1}^2-
{{\yyze}_2}$ ${{\uuvec^0}_1}$ ${{\yyze}_1}^4$ ${{\uuvec^0}_2}^2-
17$ ${{\uuvec^0}_1}$ ${{\uuvec^0}_3}$ ${{\yyze}_1}^4$ ${{\yyze}_2}$ ${{\uuvec^0}_2}-
4$ ${{\uuvec^0}_1}$ ${{\uuvec^0}_2}$ ${{\yyze}_2}^4$ ${{\uuvec^0}_3}$ ${{\yyze}_1}-
17$ ${{\uuvec^0}_1}$ ${{\uuvec^0}_2}$ ${{\uuvec^0}_3}$ ${{\yyze}_2}^2$ ${{\yyze}_1}^3+
4$ ${{\uuvec^0}_2}$ ${{\yyze}_1}^6$ ${{\uuvec^0}_3}\Big)$ ${\Omega}^3+\Big(6$ ${{\uuvec^0}_1}^2$ ${{\yyze}_2}^4+
6$ ${{\yyze}_1}^2$ ${{\yyze}_2}^2$ ${{\uuvec^0}_2}^2+{{\yyze}_1}^4$ ${{\uuvec^0}_2}$ ${{\uuvec^0}_1}-
{{\yyze}_1}^3$ ${{\yyze}_2}$ ${{\uuvec^0}_1}^2-6$ ${{\uuvec^0}_1}^2$ ${{\yyze}_1}$ ${{\yyze}_2}^3+
{{\yyze}_1}^2$ ${{\yyze}_2}^2$ ${{\uuvec^0}_3}$ ${{\uuvec^0}_2}-
12$ ${{\yyze}_1}$ ${{\yyze}_2}^3$ ${{\uuvec^0}_1}$ ${{\uuvec^0}_2}+
6$ ${{\uuvec^0}_3}$ ${{\uuvec^0}_1}$ ${{\yyze}_2}$ ${{\yyze}_1}^3-
4$ ${{\uuvec^0}_2}$ ${{\yyze}_1}^3$ ${{\yyze}_2}$ ${{\uuvec^0}_3}+
{{\yyze}_1}^2$ ${{\uuvec^0}_1}^2$ ${{\yyze}_2}$ ${{\uuvec^0}_2}+
5$ ${{\yyze}_2}^2$ ${{\uuvec^0}_3}$ ${{\uuvec^0}_1}$ ${{\yyze}_1}^2+
11$ ${{\yyze}_1}^2$ ${{\uuvec^0}_1}$ ${{\yyze}_2}^2$ ${{\uuvec^0}_2}+
2$ ${{\uuvec^0}_1}$ ${{\yyze}_1}$ ${{\yyze}_2}^3$ ${{\uuvec^0}_3}-
3$ ${{\uuvec^0}_3}$ ${{\yyze}_1}^2$ ${{\yyze}_2}$ ${{\uuvec^0}_2}^2-
{{\yyze}_2}$ ${{\uuvec^0}_2}^3$ ${{\yyze}_1}^2+{{\yyze}_1}^3$ ${{\uuvec^0}_1}$ ${{\uuvec^0}_2}^2-
{{\uuvec^0}_3}$ ${{\yyze}_1}^3$ ${{\uuvec^0}_2}^2-2$ ${{\yyze}_1}$ ${{\uuvec^0}_3}^3$ ${{\yyze}_2}^2-
4$ ${{\uuvec^0}_3}$ ${{\uuvec^0}_2}$ ${{\yyze}_1}^4-5$ ${{\uuvec^0}_2}^2$ ${{\yyze}_1}^3$ ${{\yyze}_2}+
2$ ${{\yyze}_1}^2$ ${{\uuvec^0}_3}^3$ ${{\yyze}_2}-2$ ${{\yyze}_1}^3$ ${{\uuvec^0}_1}^2$ ${{\uuvec^0}_3}+
2$ ${{\yyze}_2}$ ${{\uuvec^0}_3}$ ${{\uuvec^0}_2}$ ${{\yyze}_1}^2$ ${{\uuvec^0}_1}+
6$ ${{\yyze}_1}$ ${{\yyze}_2}^2$ ${{\uuvec^0}_3}$ ${{\uuvec^0}_2}^2-
4$ ${{\yyze}_1}^2$ ${{\uuvec^0}_1}^2$ ${{\uuvec^0}_3}$ ${{\yyze}_2}+
{{\yyze}_1}^3$ ${{\uuvec^0}_1}$ ${{\uuvec^0}_3}$ ${{\uuvec^0}_2}\Big)$ ${\Omega}^5+
\Big(4$ ${\Omega}^7$ ${{\yyze}_2}$ ${{\uuvec^0}_3}$ ${{\yyze}_1}+
{{\yyze}_2}$ $(6$ ${{\yyze}_1}$ ${{\yyze}_2}^2$ ${{\uuvec^0}_2}-
5$ ${{\yyze}_1}^2$ ${{\yyze}_2}$ ${{\uuvec^0}_2}+{{\yyze}_1}^3$ ${{\uuvec^0}_1}-
6$ ${{\yyze}_1}^2$ ${{\yyze}_2}$ ${{\uuvec^0}_3}+6$ ${{\yyze}_1}$ ${{\yyze}_2}^2$ ${{\uuvec^0}_1}-
6$ ${{\yyze}_2}^3$ ${{\uuvec^0}_1}-5$ ${{\yyze}_1}^3$ ${{\uuvec^0}_3})$ ${\Omega}^5+
{{\yyze}_1}$ ${{\yyze}_2}$ $(3$ ${{\yyze}_1}^3$ ${{\uuvec^0}_3}$ ${{\yyze}_2}-
{{\yyze}_1}^2$ ${{\yyze}_2}^2$ ${{\uuvec^0}_2}+9$ ${{\yyze}_2}^2$ ${{\yyze}_1}^2$ ${{\uuvec^0}_3}-
{{\uuvec^0}_2}$ ${{\yyze}_2}^3$ ${{\yyze}_1}+5$ ${{\yyze}_1}^4$ ${{\uuvec^0}_3}-
{{\yyze}_1}^2$ ${{\uuvec^0}_1}$ ${{\yyze}_2}^2-{{\uuvec^0}_1}$ ${{\yyze}_2}$ ${{\yyze}_1}^3+
3$ ${{\uuvec^0}_3}$ ${{\yyze}_2}^4+2$ ${{\yyze}_1}$ ${{\yyze}_2}^3$ ${{\uuvec^0}_3})$ ${\Omega}^3\Big)$ ${\evec_2}$ $+
\Big({\Omega}^8$ ${{\yyze}_1}$ ${{\uuvec^0}_3}+(2$ ${{\yyze}_1}^2$ ${{\yyze}_2}$ ${{\uuvec^0}_3}+4$ ${{\yyze}_2}^3$ ${{\uuvec^0}_1}-
4$ ${{\yyze}_1}$ ${{\yyze}_2}^2$ ${{\uuvec^0}_2}+4$ ${{\yyze}_1}^2$ ${{\yyze}_2}$ ${{\uuvec^0}_2}+{{\yyze}_1}^3$ ${{\uuvec^0}_3}-
4$ ${{\yyze}_1}$ ${{\yyze}_2}^2$ ${{\uuvec^0}_1})$ ${\Omega}^6-{{\yyze}_1}$ $({{\yyze}_2}$ ${{\yyze}_1}^3+2$ ${{\yyze}_2}^4+
5$ ${{\yyze}_1}^2$ ${{\yyze}_2}^2+2$ ${{\yyze}_1}^4)$ ${{\uuvec^0}_3}$ ${\Omega}^4\Big)$ ${{\uuvec^1}_3}+
\Big(-4$ ${\Omega}^7$ ${{\yyze}_2}^2$ ${{\uuvec^0}_3}+{{\yyze}_1}$ $(6$ ${{\yyze}_1}$ ${{\yyze}_2}^2$ ${{\uuvec^0}_2}-
5$ ${{\yyze}_1}^2$ ${{\yyze}_2}$ ${{\uuvec^0}_2}+{{\yyze}_1}^3$ ${{\uuvec^0}_1}-6$ ${{\yyze}_1}^2$ ${{\yyze}_2}$ ${{\uuvec^0}_3}+
6$ ${{\yyze}_1}$ ${{\yyze}_2}^2$ ${{\uuvec^0}_1}-6$ ${{\yyze}_2}^3$ ${{\uuvec^0}_1}-5$ ${{\yyze}_1}^3$ ${{\uuvec^0}_3})$ ${\Omega}^5+
{{\yyze}_1}^2$ $(3$ ${{\yyze}_1}^3$ ${{\uuvec^0}_3}$ ${{\yyze}_2}-{{\yyze}_1}^2$ ${{\yyze}_2}^2$ ${{\uuvec^0}_2}+
9$ ${{\yyze}_2}^2$ ${{\yyze}_1}^2$ ${{\uuvec^0}_3}-{{\uuvec^0}_2}$ ${{\yyze}_2}^3$ ${{\yyze}_1}+5$ ${{\yyze}_1}^4$ ${{\uuvec^0}_3}-
{{\yyze}_1}^2$ ${{\uuvec^0}_1}$ ${{\yyze}_2}^2-{{\uuvec^0}_1}$ ${{\yyze}_2}$ ${{\yyze}_1}^3+
3$ ${{\uuvec^0}_3}$ ${{\yyze}_2}^4+2$ ${{\yyze}_1}$ ${{\yyze}_2}^3$ ${{\uuvec^0}_3})$ ${\Omega}^3\Big)$ ${\evec_1}$ $+
\Big(\big({{\yyze}_1}$ $(-2$ ${{\uuvec^0}_3}+{{\uuvec^0}_2})$ ${\Omega}^7+{{\yyze}_1}$ $(-{{\yyze}_2}$ ${{\uuvec^0}_2}$ ${{\yyze}_1}+
8$ ${{\yyze}_2}^2$ ${{\uuvec^0}_3}+4$ ${{\yyze}_1}^2$ ${{\uuvec^0}_3}-{{\uuvec^0}_1}$ ${{\yyze}_2}$ ${{\yyze}_1}-
2$ ${{\yyze}_2}^2$ ${{\uuvec^0}_2}+9$ ${{\uuvec^0}_3}$ ${{\yyze}_1}$ ${{\yyze}_2})$ ${\Omega}^5-
{{\yyze}_1}$ ${{\yyze}_2}$ $(7$ ${{\yyze}_1}^2$ ${{\yyze}_2}$ ${{\uuvec^0}_3}-{{\yyze}_1}$ ${{\yyze}_2}^2$ ${{\uuvec^0}_2}-
{{\yyze}_2}^3$ ${{\uuvec^0}_2}-{{\yyze}_1}^2$ ${{\yyze}_2}$ ${{\uuvec^0}_1}+6$ ${{\yyze}_2}^3$ ${{\uuvec^0}_3}+
8$ ${{\yyze}_1}^3$ ${{\uuvec^0}_3}+9$ ${{\yyze}_1}$ ${{\uuvec^0}_3}$ ${{\yyze}_2}^2-
{{\yyze}_1}$ ${{\yyze}_2}^2$ ${{\uuvec^0}_1})$ ${\Omega}^3\big)$ ${{\yyun}_2}+({{\yyze}_1}^2$ ${\Omega}^7-
{{\yyze}_1}^2$ ${{\yyze}_2}$ $({{\yyze}_2}+{{\yyze}_1})$ ${\Omega}^5)$ ${{\uuvec^1}_1}+({{\yyze}_1}$ ${{\yyze}_2}$ ${\Omega}^7-
{{\yyze}_1}$ ${{\yyze}_2}^2$ $({{\yyze}_2}+{{\yyze}_1})$ ${\Omega}^5)$ ${{\uuvec^1}_2}+(-{{\yyze}_1}$ $({{\yyze}_2}+
2$ ${{\yyze}_1})$ ${\Omega}^7+{{\yyze}_1}^2$ $({{\yyze}_1}$ ${{\yyze}_2}+
2$ ${{\yyze}_1}^2+3$ ${{\yyze}_2}^2)$ ${\Omega}^5)$ ${{\uuvec^1}_3}+({{\yyze}_1}$ $(-8$ ${{\yyze}_2}^2+
2$ ${{\yyze}_1}^2+3$ ${{\yyze}_1}$ ${{\yyze}_2})$ ${\Omega}^6-{{\yyze}_1}^3$ $({{\yyze}_1}$ ${{\yyze}_2}+
2$ ${{\yyze}_1}^2+3$ ${{\yyze}_2}^2)$ ${\Omega}^4)$ ${\evec_1}$ $+\big({{\yyze}_2}$ $(3$ ${{\yyze}_1}$ ${{\yyze}_2}-
4$ ${{\yyze}_2}^2+6$ ${{\yyze}_1}^2)$ ${\Omega}^6-{{\yyze}_2}$ ${{\yyze}_1}^2$ $({{\yyze}_1}$ ${{\yyze}_2}+
2$ ${{\yyze}_1}^2+3$ ${{\yyze}_2}^2)$ ${\Omega}^4\big)$ ${\evec_2}$ $+({\Omega}^8$ ${{\yyze}_1}-{{\yyze}_1}$ ${{\yyze}_2}$ $({{\yyze}_2}+
{{\yyze}_1})$ ${\Omega}^6)$ ${\evec_3}$ $+4$ ${\Omega}^8$ ${{\yyze}_2}^2$ $\fracp{\evec_1}{x_1}$ $-
4$ ${\Omega}^8$ ${{\yyze}_1}$ ${{\yyze}_2}$ $\fracp{\evec_2}{x_1}$ $-{{\yyze}_1}$ $(-{{\uuvec^0}_2}+
2$ ${{\uuvec^0}_1})$ ${\Omega}^8+\big(5$ ${{\yyze}_1}$ ${{\yyze}_2}$ ${{\uuvec^0}_2}^2-
2$ ${{\yyze}_1}$ ${{\yyze}_2}^2$ ${{\uuvec^0}_2}+2$ ${{\yyze}_1}^2$ ${{\yyze}_2}$ ${{\uuvec^0}_2}-
{{\yyze}_1}^2$ ${{\uuvec^0}_2}$ ${{\uuvec^0}_1}-{{\uuvec^0}_2}$ ${{\yyze}_1}^3-{{\yyze}_1}$ ${{\yyze}_2}^2$ ${{\uuvec^0}_1}-
2$ ${{\yyze}_2}^2$ ${{\uuvec^0}_3}^2-6$ ${{\yyze}_2}^2$ ${{\uuvec^0}_2}^2+2$ ${{\yyze}_1}^3$ ${{\uuvec^0}_1}-
8$ ${{\uuvec^0}_1}^2$ ${{\yyze}_2}^2-9$ ${{\uuvec^0}_1}$ ${{\yyze}_2}$ ${{\uuvec^0}_2}$ ${{\yyze}_1}+
4$ ${{\uuvec^0}_2}^2$ ${{\yyze}_1}^2\big)$ ${\Omega}^6+\big(-6$ ${{\yyze}_2}^2$ ${{\uuvec^0}_3}^2$ ${{\yyze}_1}^2+
{{\yyze}_2}^4$ ${{\uuvec^0}_3}^2-6$ ${{\yyze}_2}^3$ ${{\uuvec^0}_3}^2$ ${{\yyze}_1}+
14$ ${{\yyze}_1}$ ${{\yyze}_2}^3$ ${{\uuvec^0}_1}$ ${{\uuvec^0}_2}+4$ ${{\yyze}_1}^3$ ${{\uuvec^0}_1}$ ${{\yyze}_2}$ ${{\uuvec^0}_2}+
2$ ${{\uuvec^0}_2}^2$ ${{\yyze}_1}^4-13$ ${{\uuvec^0}_3}^2$ ${{\yyze}_1}^4-
2$ ${{\yyze}_1}^2$ ${{\uuvec^0}_1}$ ${{\yyze}_2}^2$ ${{\uuvec^0}_2}+3$ ${{\yyze}_1}^4$ ${{\uuvec^0}_2}$ ${{\uuvec^0}_1}-
13$ ${{\yyze}_2}$ ${{\uuvec^0}_3}^2$ ${{\yyze}_1}^3-8$ ${{\yyze}_1}^2$ ${{\yyze}_2}^2$ ${{\uuvec^0}_2}^2+
{{\uuvec^0}_2}$ ${{\yyze}_1}^3$ ${{\yyze}_2}$ ${{\uuvec^0}_3}-6$ ${{\yyze}_1}$ ${{\yyze}_2}^3$ ${{\uuvec^0}_2}^2+
4$ ${{\uuvec^0}_1}$ ${{\yyze}_2}^4$ ${{\uuvec^0}_2}+{{\uuvec^0}_3}$ ${{\uuvec^0}_1}$ ${{\yyze}_2}$ ${{\yyze}_1}^3\big)$ ${\Omega}^4+
\big(3$ ${{\yyze}_1}^4$ ${{\uuvec^0}_1}^2$ ${{\yyze}_2}^2+4$ ${{\yyze}_1}^4$ ${{\yyze}_2}^2$ ${{\uuvec^0}_2}^2-
3$ ${{\yyze}_1}^3$ ${{\yyze}_2}^3$ ${{\uuvec^0}_1}^2-5$ ${{\yyze}_1}^3$ ${{\yyze}_2}^3$ ${{\uuvec^0}_2}^2+
7$ ${{\yyze}_1}^4$ ${{\uuvec^0}_1}$ ${{\yyze}_2}^2$ ${{\uuvec^0}_2}-3$ ${{\yyze}_1}^3$ ${{\uuvec^0}_1}$ ${{\yyze}_2}^3$ ${{\uuvec^0}_2}+
4$ ${{\yyze}_1}^2$ ${{\uuvec^0}_1}$ ${{\yyze}_2}^4$ ${{\uuvec^0}_2}+6$ ${{\yyze}_1}^5$ ${{\yyze}_2}$ ${{\uuvec^0}_2}$ ${{\uuvec^0}_1}-
{{\yyze}_1}^5$ ${{\yyze}_2}$ ${{\uuvec^0}_1}$ ${{\uuvec^0}_3}-{{\yyze}_1}^4$ ${{\yyze}_2}^2$ ${{\uuvec^0}_2}$ ${{\uuvec^0}_3}-
{{\yyze}_1}^4$ ${{\yyze}_2}^2$ ${{\uuvec^0}_1}$ ${{\uuvec^0}_3}-{{\yyze}_1}^3$ ${{\uuvec^0}_2}$ ${{\yyze}_2}^3$ ${{\uuvec^0}_3}+
2$ ${{\yyze}_1}^5$ ${{\yyze}_2}$ ${{\uuvec^0}_1}^2+18$ ${{\yyze}_1}^3$ ${{\uuvec^0}_3}^2$ ${{\yyze}_2}^3+
6$ ${{\yyze}_1}$ ${{\uuvec^0}_3}^2$ ${{\yyze}_2}^5-6$ ${{\yyze}_1}^6$ ${{\uuvec^0}_1}$ ${{\uuvec^0}_2}-
5$ ${{\yyze}_1}$ ${{\uuvec^0}_1}^2$ ${{\yyze}_2}^5+5$ ${{\yyze}_1}$ ${{\yyze}_2}^5$ ${{\uuvec^0}_2}^2+
7$ ${{\yyze}_1}^2$ ${{\uuvec^0}_1}^2$ ${{\yyze}_2}^4-11$ ${{\yyze}_1}^5$ ${{\yyze}_2}$ ${{\uuvec^0}_2}^2-
8$ ${{\uuvec^0}_1}$ ${{\uuvec^0}_2}$ ${{\yyze}_2}^6+7$ ${{\yyze}_1}^2$ ${{\yyze}_2}^4$ ${{\uuvec^0}_3}^2+
21$ ${{\yyze}_1}^4$ ${{\uuvec^0}_3}^2$ ${{\yyze}_2}^2+13$ ${{\yyze}_1}^5$ ${{\yyze}_2}$ ${{\uuvec^0}_3}^2+
16$ ${{\yyze}_1}^2$ ${{\yyze}_2}^4$ ${{\uuvec^0}_2}^2-18$ ${{\yyze}_1}$ ${{\uuvec^0}_1}$ ${{\yyze}_2}^5$ ${{\uuvec^0}_2}+
{{\yyze}_2}^6$ ${{\uuvec^0}_2}^2-6$ ${{\uuvec^0}_2}^2$ ${{\yyze}_1}^6+8$ ${{\uuvec^0}_1}^2$ ${{\yyze}_2}^6+
13$ ${{\yyze}_1}^6$ ${{\uuvec^0}_3}^2\big)$ ${\Omega}^2+2$ ${{\yyze}_1}^5$ $({{\yyze}_2}+{{\yyze}_1})$ $(-{{\uuvec^0}_1}$ ${{\yyze}_2}+
{{\uuvec^0}_2}$ ${{\yyze}_1})$ $({{\yyze}_1}$ ${{\uuvec^0}_1}+{{\yyze}_2}$ ${{\uuvec^0}_2})\Big)$ ${{\yyun}_1}-
4$ ${\Omega}^9$ ${{\yyze}_1}$ $\fracp{\evec_3}{t}$ $-2$ ${{\uuvec^0}_3}$ $({{\yyze}_1}+2$ ${{\uuvec^0}_2})$ ${\Omega}^9+
{{\uuvec^0}_3}$ $\Big(-5$ ${{\yyze}_1}$ ${{\uuvec^0}_2}^2-3$ ${{\uuvec^0}_1}$ ${{\yyze}_2}$ ${{\yyze}_1}+
4$ ${{\yyze}_2}$ ${{\uuvec^0}_2}$ ${{\uuvec^0}_1}+{{\yyze}_1}^2$ ${{\uuvec^0}_1}\Big)$ ${\Omega}^7-
2$ ${{\yyze}_1}$ $\Big(-3$ ${{\uuvec^0}_1}$ ${{\yyze}_2}$ ${{\yyze}_1}-{{\yyze}_2}^2$ ${{\uuvec^0}_2}+
2$ ${{\yyze}_1}^2$ ${{\uuvec^0}_2}\Big)$ ${\Omega}^7$ $\fracp{\evec_3}{x_2}$ $+
6$ ${\Omega}^7$ ${{\yyze}_1}^2$ ${{\uuvec^0}_3}$ ${{\yyze}_2}$ $\fracp{\evec_2}{x_1}$ $+
6$ ${\Omega}^7$ ${{\yyze}_1}$ ${{\uuvec^0}_3}$ ${{\yyze}_2}^2$ $\fracp{\evec_2}{x_2}$ $+
2$ ${{\yyze}_1}$ $(-2$ ${{\yyze}_2}^2+{{\yyze}_1}^2)$ ${{\uuvec^0}_3}$ ${\Omega}^7$ $\fracp{\evec_1}{x_1}$ $+
2$ ${{\yyze}_1}$ $\Big(-2$ ${{\uuvec^0}_1}$ ${{\yyze}_2}^2+{{\yyze}_1}^2$ ${{\uuvec^0}_1}+
3$ ${{\yyze}_2}$ ${{\uuvec^0}_2}$ ${{\yyze}_1}\Big)$ ${\Omega}^7$ $\fracp{\evec_3}{x_1}$ $+
2$ ${\Omega}^9$ ${{\yyze}_1}$ $\fracp{\evec_3}{x_3}$ ${{\uuvec^0}_3}-
2$ $(-2$ ${{\yyze}_2}^2+{{\yyze}_1}^2)$ $({{\yyze}_1}$ ${{\uuvec^0}_1}+
{{\yyze}_2}$ ${{\uuvec^0}_2})$ ${\Omega}^7$ $\fracp{\evec_1}{x_3}$ $-
6$ ${{\yyze}_1}$ ${{\yyze}_2}$ $({{\yyze}_1}$ ${{\uuvec^0}_1}+
{{\yyze}_2}$ ${{\uuvec^0}_2})$ ${\Omega}^7$ $\fracp{\evec_2}{x_3}$ $+2$ ${{\yyze}_2}$ $(-2$ ${{\yyze}_2}^2+
{{\yyze}_1}^2)$ ${{\uuvec^0}_3}$ ${\Omega}^7$ $\fracp{\evec_1}{x_2}$ $+({{\yyze}_2}^2$ ${{\uuvec^0}_3}$ ${\Omega}^5$ ${{\yyze}_1}-
{\Omega}^7$ ${{\yyze}_1}$ ${{\uuvec^0}_3})$ ${{\yyun}_2}^2\bigg)$ $\bigg/$ ${\Omega}^{10},
$

~

$
\frac{d{\uuvec^1}_2}{dt}$ $=$ $
$
$
$ $1/4$ $\bigg(-4$ ${\Omega}^9$ ${{\yyze}_2}$ $\fracp{\evec_3}{t}$ $-2$ ${{\uuvec^0}_3}$ $({{\yyze}_2}-
4$ ${{\uuvec^0}_1})$ ${\Omega}^9-{{\uuvec^0}_3}$ $\Big(3$ ${{\uuvec^0}_1}$ ${{\yyze}_2}^2+
4$ ${{\yyze}_1}^2$ ${{\uuvec^0}_1}-4$ ${{\uuvec^0}_1}$ ${{\yyze}_1}$ ${{\uuvec^0}_2}+
{{\yyze}_2}$ ${{\uuvec^0}_2}^2-2$ ${{\yyze}_1}^2$ ${{\uuvec^0}_2}-{{\uuvec^0}_1}$ ${{\yyze}_2}$ ${{\yyze}_1}+
4$ ${{\yyze}_2}$ ${{\uuvec^0}_1}^2\Big)$ ${\Omega}^7+2$ $(-{{\yyze}_2}^2+2$ ${{\yyze}_1}^2)$ $({{\yyze}_1}$ ${{\uuvec^0}_1}+
{{\yyze}_2}$ ${{\uuvec^0}_2})$ ${\Omega}^7$ $\fracp{\evec_2}{x_3}$ $+
2$ ${\Omega}^9$ ${{\yyze}_2}$ $\fracp{\evec_3}{x_3}$ ${{\uuvec^0}_3}+2$ ${{\yyze}_2}$ $\big(-2$ ${{\uuvec^0}_1}$ ${{\yyze}_2}^2+
{{\yyze}_1}^2$ ${{\uuvec^0}_1}+3$ ${{\yyze}_2}$ ${{\uuvec^0}_2}$ ${{\yyze}_1}\big)$ ${\Omega}^7$ $\fracp{\evec_3}{x_1}$ $-
2$ ${{\yyze}_2}$ $\big(-3$ ${{\uuvec^0}_1}$ ${{\yyze}_2}$ ${{\yyze}_1}-{{\yyze}_2}^2$ ${{\uuvec^0}_2}+
2$ ${{\yyze}_1}^2$ ${{\uuvec^0}_2}\big)$ ${\Omega}^7$ $\fracp{\evec_3}{x_2}$ $+
\Big(-2$ ${{\yyze}_1}$ $(-2$ ${{\yyze}_2}^2+{{\yyze}_1}^2)$ ${\Omega}^6$ ${\evec_1}$ $-
6$ ${\Omega}^6$ ${{\yyze}_1}^2$ ${{\yyze}_2}$ ${\evec_2}$ $-4$ ${\Omega}^8$ ${{\yyze}_1}$ ${{\yyze}_2}$ $\fracp{\evec_1}{x_2}$ $+
4$ ${\Omega}^8$ ${{\yyze}_1}^2$ $\fracp{\evec_2}{x_2}$ $+(-2$ ${{\uuvec^0}_2}$ ${{\yyze}_1}+
{{\uuvec^0}_1}$ ${{\yyze}_2})$ ${\Omega}^8+(-2$ ${{\uuvec^0}_2}^2$ ${{\yyze}_1}^2-
6$ ${{\yyze}_1}^2$ ${{\uuvec^0}_1}^2+4$ ${{\uuvec^0}_2}$ ${{\yyze}_1}^3-2$ ${{\yyze}_1}^2$ ${{\uuvec^0}_3}^2-
{{\yyze}_2}^3$ ${{\uuvec^0}_1}-2$ ${{\yyze}_1}^2$ ${{\yyze}_2}$ ${{\uuvec^0}_1}-
{{\yyze}_2}$ ${{\uuvec^0}_3}$ ${{\uuvec^0}_2}$ ${{\yyze}_1}+3$ ${{\yyze}_1}$ ${{\yyze}_2}^2$ ${{\uuvec^0}_2}+
{{\yyze}_2}^2$ ${{\uuvec^0}_3}^2+{{\uuvec^0}_3}^2$ ${{\yyze}_1}$ ${{\yyze}_2}-
11$ ${{\uuvec^0}_1}$ ${{\yyze}_2}$ ${{\uuvec^0}_2}$ ${{\yyze}_1}+
3$ ${{\uuvec^0}_1}^2$ ${{\yyze}_2}^2)$ ${\Omega}^6+(-4$ ${{\uuvec^0}_2}^2$ ${{\yyze}_1}^4-
{{\yyze}_1}^4$ ${{\uuvec^0}_2}$ ${{\uuvec^0}_1}-8$ ${{\yyze}_2}^2$ ${{\uuvec^0}_3}^2$ ${{\yyze}_1}^2+
3$ ${{\yyze}_1}^3$ ${{\yyze}_2}$ ${{\uuvec^0}_1}^2+{{\yyze}_1}^2$ ${{\yyze}_2}^2$ ${{\uuvec^0}_3}$ ${{\uuvec^0}_2}+
10$ ${{\uuvec^0}_1}^2$ ${{\yyze}_1}$ ${{\yyze}_2}^3-2$ ${{\yyze}_1}$ ${{\yyze}_2}^3$ ${{\uuvec^0}_2}^2-
2$ ${{\yyze}_2}^4$ ${{\uuvec^0}_3}^2+{{\uuvec^0}_1}^2$ ${{\yyze}_2}^4-{{\yyze}_1}^2$ ${{\yyze}_2}^2$ ${{\uuvec^0}_2}^2+
5$ ${{\yyze}_1}$ ${{\yyze}_2}^3$ ${{\uuvec^0}_1}$ ${{\uuvec^0}_2}-7$ ${{\yyze}_2}^3$ ${{\uuvec^0}_3}^2$ ${{\yyze}_1}+
{{\yyze}_1}^4$ ${{\uuvec^0}_1}^2-{{\yyze}_2}$ ${{\uuvec^0}_3}^2$ ${{\yyze}_1}^3+
2$ ${{\uuvec^0}_2}$ ${{\yyze}_1}$ ${{\yyze}_2}^3$ ${{\uuvec^0}_3}+7$ ${{\yyze}_1}^2$ ${{\yyze}_2}^2$ ${{\uuvec^0}_1}^2-
11$ ${{\yyze}_1}^2$ ${{\uuvec^0}_1}$ ${{\yyze}_2}^2$ ${{\uuvec^0}_2}-{{\uuvec^0}_2}^2$ ${{\yyze}_1}^3$ ${{\yyze}_2}+
{{\uuvec^0}_3}^2$ ${{\yyze}_1}^4+{{\yyze}_2}^2$ ${{\uuvec^0}_3}$ ${{\uuvec^0}_1}$ ${{\yyze}_1}^2+
2$ ${{\uuvec^0}_1}$ ${{\yyze}_2}^4$ ${{\uuvec^0}_2})$ ${\Omega}^4+
{{\yyze}_2}$ $(7$ ${{\yyze}_2}$ ${{\uuvec^0}_3}^2$ ${{\yyze}_1}^4-
18$ ${{\yyze}_2}^4$ ${{\uuvec^0}_1}^2$ ${{\yyze}_1}-14$ ${{\yyze}_2}^3$ ${{\yyze}_1}^2$ ${{\uuvec^0}_1}^2-
6$ ${{\yyze}_1}$ ${{\yyze}_2}^4$ ${{\uuvec^0}_2}$ ${{\uuvec^0}_1}+4$ ${{\yyze}_1}^2$ ${{\yyze}_2}^3$ ${{\uuvec^0}_2}^2+
18$ ${{\uuvec^0}_1}$ ${{\uuvec^0}_2}$ ${{\yyze}_1}^5-{{\uuvec^0}_2}$ ${{\yyze}_2}^4$ ${{\uuvec^0}_3}$ ${{\yyze}_1}-
20$ ${{\yyze}_1}^3$ ${{\uuvec^0}_1}^2$ ${{\yyze}_2}^2-5$ ${{\yyze}_1}^3$ ${{\yyze}_2}^2$ ${{\uuvec^0}_2}^2-
{{\uuvec^0}_2}$ ${{\yyze}_2}^3$ ${{\uuvec^0}_3}$ ${{\yyze}_1}^2+{{\uuvec^0}_3}^2$ ${{\yyze}_2}^5-
16$ ${{\yyze}_1}^4$ ${{\yyze}_2}$ ${{\uuvec^0}_1}^2+31$ ${{\yyze}_1}^4$ ${{\yyze}_2}$ ${{\uuvec^0}_2}$ ${{\uuvec^0}_1}+
6$ ${{\yyze}_2}^4$ ${{\uuvec^0}_3}^2$ ${{\yyze}_1}-2$ ${{\yyze}_2}^5$ ${{\uuvec^0}_1}$ ${{\uuvec^0}_2}+
9$ ${{\uuvec^0}_3}^2$ ${{\yyze}_2}^3$ ${{\yyze}_1}^2+7$ ${{\yyze}_1}^4$ ${{\yyze}_2}$ ${{\uuvec^0}_2}^2+
2$ ${{\yyze}_1}$ ${{\yyze}_2}^4$ ${{\uuvec^0}_2}^2-8$ ${{\uuvec^0}_2}^2$ ${{\yyze}_1}^5+
7$ ${{\yyze}_1}^3$ ${{\uuvec^0}_1}$ ${{\yyze}_2}^2$ ${{\uuvec^0}_2}-
{{\yyze}_2}^2$ ${{\uuvec^0}_3}$ ${{\yyze}_1}^3$ ${{\uuvec^0}_1}+7$ ${{\uuvec^0}_3}^2$ ${{\yyze}_1}^3$ ${{\yyze}_2}^2+
26$ ${{\yyze}_1}^2$ ${{\uuvec^0}_1}$ ${{\yyze}_2}^3$ ${{\uuvec^0}_2}-
{{\uuvec^0}_1}$ ${{\yyze}_2}^3$ ${{\uuvec^0}_3}$ ${{\yyze}_1}^2)$ ${\Omega}^2+
2$ ${{\yyze}_1}^3$ ${{\yyze}_2}^2$ $({{\yyze}_2}+{{\yyze}_1})$ $({{\yyze}_1}$ ${{\uuvec^0}_1}+
{{\yyze}_2}$ ${{\uuvec^0}_2})$ $(-{{\uuvec^0}_1}$ ${{\yyze}_2}+{{\uuvec^0}_2}$ ${{\yyze}_1})\Big)$ ${{\yyun}_2}+
\Big({{\yyze}_2}$ ${\Omega}^8$ ${{\uuvec^0}_3}-(-{{\yyze}_1}$ ${{\yyze}_2}^2-2$ ${{\yyze}_2}$ ${{\yyze}_1}^2+
{{\yyze}_1}^3-{{\yyze}_2}^3)$ ${{\uuvec^0}_3}$ ${\Omega}^6-{{\yyze}_2}$ $({{\yyze}_2}$ ${{\yyze}_1}^3+
2$ ${{\yyze}_2}^4+5$ ${{\yyze}_1}^2$ ${{\yyze}_2}^2+2$ ${{\yyze}_1}^4)$ ${{\uuvec^0}_3}$ ${\Omega}^4\Big)$ ${{\uuvec^1}_3}+
\Big({{\uuvec^0}_3}$ ${{\uuvec^0}_2}^2$ ${{\yyze}_1}^5+2$ ${{\yyze}_2}^2$ ${{\yyze}_1}^3$ ${{\uuvec^0}_2}^3-
{{\yyze}_2}^5$ ${{\uuvec^0}_1}^2$ ${{\uuvec^0}_3}+{{\yyze}_1}^3$ ${{\yyze}_2}^3$ ${{\uuvec^0}_1}^2-
{{\yyze}_1}^3$ ${{\yyze}_2}^3$ ${{\uuvec^0}_2}^2-5$ ${{\yyze}_1}$ ${{\uuvec^0}_3}$ ${{\yyze}_2}^4$ ${{\uuvec^0}_2}^2-
8$ ${{\yyze}_1}^3$ ${{\uuvec^0}_3}$ ${{\yyze}_2}^2$ ${{\uuvec^0}_1}^2-
15$ ${{\yyze}_1}^4$ ${{\uuvec^0}_3}$ ${{\yyze}_2}$ ${{\uuvec^0}_1}^2+
5$ ${{\yyze}_1}^2$ ${{\uuvec^0}_3}$ ${{\yyze}_2}^3$ ${{\uuvec^0}_2}^2-
{{\yyze}_1}^4$ ${{\uuvec^0}_1}$ ${{\yyze}_2}^2$ ${{\uuvec^0}_2}-{{\yyze}_1}^2$ ${{\uuvec^0}_1}^2$ ${{\yyze}_2}^3$ ${{\uuvec^0}_2}-{{\yyze}_1}^3$ ${{\uuvec^0}_1}$ ${{\yyze}_2}^3$ ${{\uuvec^0}_2}-
{{\uuvec^0}_1}$ ${{\yyze}_1}^3$ ${{\yyze}_2}^2$ ${{\uuvec^0}_2}^2-
3$ ${{\yyze}_1}^3$ ${{\yyze}_2}^2$ ${{\uuvec^0}_2}$ ${{\uuvec^0}_1}^2+
{{\yyze}_1}^2$ ${{\uuvec^0}_1}$ ${{\yyze}_2}^4$ ${{\uuvec^0}_2}-
3$ ${{\yyze}_1}^3$ ${{\yyze}_2}^2$ ${{\uuvec^0}_3}$ ${{\uuvec^0}_2}^2-
11$ ${{\uuvec^0}_1}^2$ ${{\yyze}_1}^2$ ${{\yyze}_2}^3$ ${{\uuvec^0}_3}-
2$ ${{\uuvec^0}_1}^2$ ${{\yyze}_2}^4$ ${{\uuvec^0}_3}$ ${{\yyze}_1}-
4$ ${{\yyze}_1}^2$ ${{\uuvec^0}_2}^2$ ${{\uuvec^0}_1}$ ${{\yyze}_2}^3-
{{\yyze}_1}$ ${{\yyze}_2}^4$ ${{\uuvec^0}_2}$ ${{\uuvec^0}_1}^2+4$ ${{\yyze}_1}^5$ ${{\yyze}_2}$ ${{\uuvec^0}_2}$ ${{\uuvec^0}_3}-
{{\yyze}_1}^5$ ${{\yyze}_2}$ ${{\uuvec^0}_1}$ ${{\uuvec^0}_3}-2$ ${{\yyze}_1}^3$ ${{\yyze}_2}^3$ ${{\uuvec^0}_1}$ ${{\uuvec^0}_3}-
2$ ${{\yyze}_2}^5$ ${{\yyze}_1}$ ${{\uuvec^0}_1}$ ${{\uuvec^0}_3}-
6$ ${{\yyze}_1}^2$ ${{\yyze}_2}^4$ ${{\uuvec^0}_1}$ ${{\uuvec^0}_3}+
{{\yyze}_1}^4$ ${{\yyze}_2}^2$ ${{\uuvec^0}_2}$ ${{\uuvec^0}_3}+
2$ ${{\yyze}_1}^2$ ${{\yyze}_2}^4$ ${{\uuvec^0}_2}$ ${{\uuvec^0}_3}-
5$ ${{\yyze}_1}^4$ ${{\yyze}_2}^2$ ${{\uuvec^0}_1}$ ${{\uuvec^0}_3}+
5$ ${{\yyze}_1}^3$ ${{\uuvec^0}_2}$ ${{\yyze}_2}^3$ ${{\uuvec^0}_3}-
6$ ${{\uuvec^0}_1}$ ${{\uuvec^0}_3}$ ${{\yyze}_1}^4$ ${{\yyze}_2}$ ${{\uuvec^0}_2}+
10$ ${{\uuvec^0}_1}$ ${{\uuvec^0}_2}$ ${{\yyze}_2}^4$ ${{\uuvec^0}_3}$ ${{\yyze}_1}+
9$ ${{\uuvec^0}_1}$ ${{\uuvec^0}_2}$ ${{\uuvec^0}_3}$ ${{\yyze}_2}^2$ ${{\yyze}_1}^3+
{{\yyze}_1}^2$ ${{\uuvec^0}_1}^2$ ${{\yyze}_2}^4-{{\yyze}_1}^2$ ${{\yyze}_2}^4$ ${{\uuvec^0}_2}^2+
2$ ${{\uuvec^0}_3}$ ${{\uuvec^0}_2}^2$ ${{\yyze}_2}^5-2$ ${{\uuvec^0}_1}$ ${{\yyze}_2}^6$ ${{\uuvec^0}_3}-
{{\yyze}_1}^2$ ${{\yyze}_2}^3$ ${{\uuvec^0}_2}^3-{{\yyze}_2}^2$ ${{\yyze}_1}^3$ ${{\uuvec^0}_1}^3+
{{\yyze}_1}$ ${{\uuvec^0}_1}$ ${{\yyze}_2}^5$ ${{\uuvec^0}_2}+{{\yyze}_1}^4$ ${{\yyze}_2}$ ${{\uuvec^0}_3}$ ${{\uuvec^0}_2}^2+
2$ ${{\yyze}_1}$ ${{\yyze}_2}^5$ ${{\uuvec^0}_2}$ ${{\uuvec^0}_3}-
31$ ${{\uuvec^0}_1}$ ${{\uuvec^0}_2}$ ${{\yyze}_2}^3$ ${{\uuvec^0}_3}$ ${{\yyze}_1}^2\Big)$ ${\Omega}^3+
\big(-{{\yyze}_2}$ ${\Omega}^7$ ${{\uuvec^0}_3}+{{\yyze}_2}^3$ ${{\uuvec^0}_3}$ ${\Omega}^5\big)$ ${{\yyun}_2}^2+
\Big(4$ ${\Omega}^8$ ${{\yyze}_1}$ ${{\uuvec^0}_1}-{{\yyze}_2}$ $\big(4$ ${{\yyze}_1}^2$ ${{\uuvec^0}_2}-
4$ ${{\yyze}_2}$ ${{\uuvec^0}_2}$ ${{\yyze}_1}+4$ ${{\uuvec^0}_1}$ ${{\yyze}_2}^2-4$ ${{\uuvec^0}_1}$ ${{\yyze}_2}$ ${{\yyze}_1}+
{{\uuvec^0}_3}$ ${{\yyze}_1}$ ${{\yyze}_2}\big)$ ${\Omega}^6+
{{\yyze}_1}$ ${{\yyze}_2}^2$ $\big({{\yyze}_1}^2$ ${{\uuvec^0}_1}+{{\uuvec^0}_3}$ ${{\yyze}_1}$ ${{\yyze}_2}+
{{\yyze}_2}$ ${{\uuvec^0}_2}$ ${{\yyze}_1}+{{\yyze}_2}^2$ ${{\uuvec^0}_3}\big)$ ${\Omega}^4-
{{\yyze}_1}^4$ ${{\yyze}_2}$ $({{\yyze}_2}+{{\yyze}_1})$ $({{\yyze}_1}$ ${{\uuvec^0}_1}+
{{\yyze}_2}$ ${{\uuvec^0}_2})$ ${\Omega}^2\Big)$ ${{\uuvec^1}_2}+\Big(-4$ ${{\yyze}_2}$ ${\Omega}^8$ ${{\uuvec^0}_1}-
{{\yyze}_1}^2$ $({{\yyze}_1}$ ${{\uuvec^0}_1}+{{\yyze}_2}$ ${{\uuvec^0}_3}+{{\yyze}_2}$ ${{\uuvec^0}_2})$ ${\Omega}^6+
{{\yyze}_2}$ ${{\yyze}_1}^2$ $({{\yyze}_1}^2$ ${{\uuvec^0}_1}+{{\uuvec^0}_3}$ ${{\yyze}_1}$ ${{\yyze}_2}+
{{\yyze}_2}$ ${{\uuvec^0}_2}$ ${{\yyze}_1}+{{\yyze}_2}^2$ ${{\uuvec^0}_3})$ ${\Omega}^4+
{{\yyze}_1}^3$ ${{\yyze}_2}^2$ $({{\yyze}_2}+{{\yyze}_1})$ $({{\yyze}_1}$ ${{\uuvec^0}_1}+
{{\yyze}_2}$ ${{\uuvec^0}_2})$ ${\Omega}^2\Big)$ ${{\uuvec^1}_1}+\Big({{\yyze}_1}^2$ ${{\yyze}_2}^5$ ${{\uuvec^0}_2}^3-
3$ ${{\uuvec^0}_1}^2$ ${{\uuvec^0}_3}$ ${{\yyze}_2}^7-2$ ${{\uuvec^0}_1}^2$ ${{\uuvec^0}_3}$ ${{\yyze}_1}^7+
{{\yyze}_2}^2$ ${{\yyze}_1}^5$ ${{\uuvec^0}_1}^3+10$ ${{\uuvec^0}_2}^2$ ${{\yyze}_1}^5$ ${{\yyze}_2}^2$ ${{\uuvec^0}_3}+
17$ ${{\yyze}_1}^3$ ${{\yyze}_2}^4$ ${{\uuvec^0}_3}$ ${{\uuvec^0}_2}^2+
3$ ${{\yyze}_1}^4$ ${{\uuvec^0}_3}$ ${{\uuvec^0}_2}^2$ ${{\yyze}_2}^3+
17$ ${{\uuvec^0}_1}^2$ ${{\uuvec^0}_3}$ ${{\yyze}_1}^6$ ${{\yyze}_2}+
17$ ${{\uuvec^0}_1}^2$ ${{\uuvec^0}_3}$ ${{\yyze}_2}^6$ ${{\yyze}_1}+
21$ ${{\uuvec^0}_1}^2$ ${{\uuvec^0}_3}$ ${{\yyze}_1}^4$ ${{\yyze}_2}^3+
33$ ${{\uuvec^0}_1}^2$ ${{\uuvec^0}_3}$ ${{\yyze}_2}^4$ ${{\yyze}_1}^3+
3$ ${{\yyze}_2}^4$ ${{\yyze}_1}^3$ ${{\uuvec^0}_2}^2$ ${{\uuvec^0}_1}+
15$ ${{\yyze}_1}^5$ ${{\yyze}_2}^2$ ${{\uuvec^0}_3}$ ${{\uuvec^0}_1}^2+
3$ ${{\yyze}_2}^3$ ${{\yyze}_1}^4$ ${{\uuvec^0}_2}$ ${{\uuvec^0}_1}^2+
6$ ${{\yyze}_2}^6$ ${{\uuvec^0}_3}$ ${{\uuvec^0}_2}^2$ ${{\yyze}_1}+
3$ ${{\yyze}_2}^5$ ${{\uuvec^0}_3}$ ${{\uuvec^0}_2}^2$ ${{\yyze}_1}^2+
12$ ${{\yyze}_1}^2$ ${{\yyze}_2}^5$ ${{\uuvec^0}_3}$ ${{\uuvec^0}_2}$ ${{\uuvec^0}_1}+
11$ ${{\yyze}_1}^5$ ${{\yyze}_2}^2$ ${{\uuvec^0}_3}$ ${{\uuvec^0}_2}$ ${{\uuvec^0}_1}-
3$ ${{\yyze}_1}^6$ ${{\yyze}_2}$ ${{\uuvec^0}_3}$ ${{\uuvec^0}_2}$ ${{\uuvec^0}_1}+
11$ ${{\yyze}_1}^4$ ${{\yyze}_2}^3$ ${{\uuvec^0}_3}$ ${{\uuvec^0}_1}$ ${{\uuvec^0}_2}+
11$ ${{\yyze}_1}^3$ ${{\yyze}_2}^4$ ${{\uuvec^0}_3}$ ${{\uuvec^0}_1}$ ${{\uuvec^0}_2}-
{{\yyze}_2}^7$ ${{\uuvec^0}_3}$ ${{\uuvec^0}_2}^2+{{\yyze}_1}$ ${{\yyze}_2}^6$ ${{\uuvec^0}_2}^3+
{{\yyze}_2}^3$ ${{\yyze}_1}^4$ ${{\uuvec^0}_1}^3+3$ ${{\yyze}_2}^5$ ${{\yyze}_1}^2$ ${{\uuvec^0}_1}$ ${{\uuvec^0}_2}^2+
3$ ${{\yyze}_2}^4$ ${{\yyze}_1}^3$ ${{\uuvec^0}_1}^2$ ${{\uuvec^0}_2}-
2$ ${{\yyze}_1}$ ${{\yyze}_2}^6$ ${{\uuvec^0}_3}$ ${{\uuvec^0}_2}$ ${{\uuvec^0}_1}\Big)$ ${\Omega}+
\Big(2$ ${{\uuvec^0}_3}$ ${{\yyze}_1}^3$ ${{\uuvec^0}_2}^2-
3$ ${{\uuvec^0}_3}$ ${{\uuvec^0}_1}$ ${{\yyze}_2}$ ${{\yyze}_1}^3+
{{\yyze}_2}^4$ ${{\uuvec^0}_3}$ ${{\uuvec^0}_1}-{{\yyze}_1}^2$ ${{\yyze}_2}^2$ ${{\uuvec^0}_1}^2+
11$ ${{\yyze}_2}^2$ ${{\uuvec^0}_3}$ ${{\uuvec^0}_1}$ ${{\yyze}_1}^2-
4$ ${{\uuvec^0}_1}$ ${{\yyze}_1}$ ${{\yyze}_2}^3$ ${{\uuvec^0}_3}+
{{\yyze}_1}$ ${{\yyze}_2}^2$ ${{\uuvec^0}_1}^2$ ${{\uuvec^0}_2}-
3$ ${{\yyze}_2}$ ${{\uuvec^0}_3}$ ${{\uuvec^0}_2}$ ${{\yyze}_1}^2$ ${{\uuvec^0}_1}+
{{\uuvec^0}_1}$ ${{\yyze}_1}^2$ ${{\yyze}_2}$ ${{\uuvec^0}_2}^2-
13$ ${{\uuvec^0}_2}$ ${{\yyze}_1}^3$ ${{\yyze}_2}$ ${{\uuvec^0}_3}-
4$ ${{\uuvec^0}_3}$ ${{\yyze}_1}$ ${{\yyze}_2}^2$ ${{\uuvec^0}_1}^2+
{{\uuvec^0}_3}$ ${{\uuvec^0}_2}$ ${{\yyze}_1}^4-7$ ${{\uuvec^0}_3}$ ${{\yyze}_1}^2$ ${{\yyze}_2}$ ${{\uuvec^0}_2}^2+
2$ ${{\yyze}_1}^2$ ${{\uuvec^0}_3}^3$ ${{\yyze}_2}-
2$ ${{\yyze}_1}$ ${{\uuvec^0}_3}$ ${{\uuvec^0}_1}$ ${{\yyze}_2}^2$ ${{\uuvec^0}_2}-
2$ ${{\yyze}_1}^3$ ${{\uuvec^0}_3}^3-{{\yyze}_1}$ ${{\yyze}_2}^2$ ${{\uuvec^0}_2}^3+
3$ ${{\yyze}_1}^2$ ${{\yyze}_2}^2$ ${{\uuvec^0}_3}$ ${{\uuvec^0}_2}+{{\yyze}_1}^3$ ${{\uuvec^0}_1}$ ${{\yyze}_2}$ ${{\uuvec^0}_2}-
{{\yyze}_1}$ ${{\yyze}_2}^3$ ${{\uuvec^0}_1}$ ${{\uuvec^0}_2}-{{\yyze}_1}$ ${{\yyze}_2}^2$ ${{\uuvec^0}_3}$ ${{\uuvec^0}_2}^2+
{{\yyze}_1}^2$ ${{\yyze}_2}^2$ ${{\uuvec^0}_2}^2-2$ ${{\uuvec^0}_2}$ ${{\yyze}_1}$ ${{\yyze}_2}^3$ ${{\uuvec^0}_3}\Big)$ ${\Omega}^5+
\Big(4$ ${{\yyze}_1}^2$ ${\Omega}^8$ $\fracp{\evec_2}{x_3}$ $-
4$ ${{\yyze}_1}$ ${\Omega}^8$ $\fracp{\evec_1}{x_3}$ ${{\yyze}_2}\Big)$ ${{\yyun}_3}+\Big(4$ ${\Omega}^9$ ${{\yyze}_1}-
{{\yyze}_1}$ $({{\yyze}_2}$ ${{\uuvec^0}_3}-4$ ${{\uuvec^0}_1}$ ${{\yyze}_2}+4$ ${{\uuvec^0}_2}$ ${{\yyze}_1})$ ${\Omega}^7+
(-3$ ${{\uuvec^0}_2}$ ${{\yyze}_2}$ ${{\yyze}_1}^3+5$ ${{\uuvec^0}_2}$ ${{\yyze}_2}^3$ ${{\yyze}_1}-
{{\uuvec^0}_2}$ ${{\yyze}_2}^4+3$ ${{\yyze}_1}$ ${{\yyze}_2}^3$ ${{\uuvec^0}_1}-6$ ${{\yyze}_2}^4$ ${{\uuvec^0}_1}-
4$ ${{\yyze}_1}^2$ ${{\yyze}_2}^2$ ${{\uuvec^0}_2}+{{\yyze}_1}$ ${{\yyze}_2}^3$ ${{\uuvec^0}_3}-{{\yyze}_1}^4$ ${{\uuvec^0}_1}+
{{\yyze}_2}^2$ ${{\yyze}_1}^2$ ${{\uuvec^0}_3}+{{\yyze}_1}^2$ ${{\uuvec^0}_1}$ ${{\yyze}_2}^2)$ ${\Omega}^5+
{{\yyze}_2}$ $(3$ ${{\yyze}_1}^2$ ${{\yyze}_2}^2+{{\yyze}_2}^4+{{\yyze}_1}^4+
{{\yyze}_2}$ ${{\yyze}_1}^3)$ $({{\yyze}_1}$ ${{\uuvec^0}_1}+{{\yyze}_2}$ ${{\uuvec^0}_2})$ ${\Omega}^3\Big)$ ${\evec_3}$ $+
\Big(4$ ${\Omega}^7$ ${{\yyze}_2}$ ${{\uuvec^0}_3}$ ${{\yyze}_1}+{{\yyze}_1}$ $(-8$ ${{\yyze}_1}^2$ ${{\yyze}_2}$ ${{\uuvec^0}_3}+
{{\yyze}_1}$ ${{\yyze}_2}^2$ ${{\uuvec^0}_2}-3$ ${{\yyze}_1}$ ${{\uuvec^0}_3}$ ${{\yyze}_2}^2+3$ ${{\yyze}_1}^3$ ${{\uuvec^0}_3}-
3$ ${{\yyze}_2}^3$ ${{\uuvec^0}_3}+{{\yyze}_1}^2$ ${{\yyze}_2}$ ${{\uuvec^0}_1})$ ${\Omega}^5+
{{\yyze}_1}$ ${{\yyze}_2}$ $(3$ ${{\yyze}_1}^3$ ${{\uuvec^0}_3}$ ${{\yyze}_2}-{{\yyze}_1}^2$ ${{\yyze}_2}^2$ ${{\uuvec^0}_2}+
9$ ${{\yyze}_2}^2$ ${{\yyze}_1}^2$ ${{\uuvec^0}_3}-{{\uuvec^0}_2}$ ${{\yyze}_2}^3$ ${{\yyze}_1}+5$ ${{\yyze}_1}^4$ ${{\uuvec^0}_3}-
{{\yyze}_1}^2$ ${{\uuvec^0}_1}$ ${{\yyze}_2}^2-{{\uuvec^0}_1}$ ${{\yyze}_2}$ ${{\yyze}_1}^3+3$ ${{\uuvec^0}_3}$ ${{\yyze}_2}^4+
2$ ${{\yyze}_1}$ ${{\yyze}_2}^3$ ${{\uuvec^0}_3})$ ${\Omega}^3\Big)$ ${\evec_1}$ $+\Big(-{{\yyze}_2}$ ${\Omega}^7$ ${{\uuvec^0}_3}+
{{\yyze}_2}$ $(8$ ${{\yyze}_1}^2$ ${{\uuvec^0}_3}-{{\yyze}_2}$ ${{\uuvec^0}_2}$ ${{\yyze}_1}-
{{\uuvec^0}_1}$ ${{\yyze}_2}$ ${{\yyze}_1}+5$ ${{\uuvec^0}_3}$ ${{\yyze}_1}$ ${{\yyze}_2}+
2$ ${{\yyze}_2}^2$ ${{\uuvec^0}_3})$ ${\Omega}^5-{{\yyze}_1}$ ${{\yyze}_2}$ $(8$ ${{\yyze}_1}$ ${{\uuvec^0}_3}$ ${{\yyze}_2}^2-
{{\yyze}_1}$ ${{\yyze}_2}^2$ ${{\uuvec^0}_1}+7$ ${{\yyze}_1}^2$ ${{\yyze}_2}$ ${{\uuvec^0}_3}-
{{\yyze}_1}$ ${{\yyze}_2}^2$ ${{\uuvec^0}_2}+6$ ${{\yyze}_2}^3$ ${{\uuvec^0}_3}+7$ ${{\yyze}_1}^3$ ${{\uuvec^0}_3}-
{{\yyze}_1}^2$ ${{\yyze}_2}$ ${{\uuvec^0}_1}-{{\yyze}_2}^3$ ${{\uuvec^0}_2})$ ${\Omega}^3\Big)$ ${{\yyun}_1}^2+
\Big(\big({{\yyze}_2}$ $(-2$ ${{\uuvec^0}_3}+{{\uuvec^0}_2})$ ${\Omega}^7+
{{\yyze}_2}$ $(-{{\yyze}_2}$ ${{\uuvec^0}_2}$ ${{\yyze}_1}+8$ ${{\yyze}_2}^2$ ${{\uuvec^0}_3}+
4$ ${{\yyze}_1}^2$ ${{\uuvec^0}_3}-{{\uuvec^0}_1}$ ${{\yyze}_2}$ ${{\yyze}_1}-2$ ${{\yyze}_2}^2$ ${{\uuvec^0}_2}+
9$ ${{\uuvec^0}_3}$ ${{\yyze}_1}$ ${{\yyze}_2})$ ${\Omega}^5-{{\yyze}_2}^2$ $(7$ ${{\yyze}_1}^2$ ${{\yyze}_2}$ ${{\uuvec^0}_3}-
{{\yyze}_1}$ ${{\yyze}_2}^2$ ${{\uuvec^0}_2}-{{\yyze}_2}^3$ ${{\uuvec^0}_2}-{{\yyze}_1}^2$ ${{\yyze}_2}$ ${{\uuvec^0}_1}+
6$ ${{\yyze}_2}^3$ ${{\uuvec^0}_3}+8$ ${{\yyze}_1}^3$ ${{\uuvec^0}_3}+9$ ${{\yyze}_1}$ ${{\uuvec^0}_3}$ ${{\yyze}_2}^2-
{{\yyze}_1}$ ${{\yyze}_2}^2$ ${{\uuvec^0}_1})$ ${\Omega}^3\big)$ ${{\yyun}_2}+({{\yyze}_1}$ ${{\yyze}_2}$ ${\Omega}^7-
{{\yyze}_1}$ ${{\yyze}_2}^2$ $({{\yyze}_2}+{{\yyze}_1})$ ${\Omega}^5)$ ${{\uuvec^1}_1}+({{\yyze}_2}^2$ ${\Omega}^7-
{{\yyze}_2}^3$ $({{\yyze}_2}+{{\yyze}_1})$ ${\Omega}^5)$ ${{\uuvec^1}_2}+(-{{\yyze}_2}$ $({{\yyze}_2}+
2$ ${{\yyze}_1})$ ${\Omega}^7+{{\yyze}_1}$ ${{\yyze}_2}$ $({{\yyze}_1}$ ${{\yyze}_2}+2$ ${{\yyze}_1}^2+
3$ ${{\yyze}_2}^2)$ ${\Omega}^5)$ ${{\uuvec^1}_3}+\big((6$ ${{\yyze}_2}$ ${{\yyze}_1}^2-2$ ${{\yyze}_1}^3-
4$ ${{\yyze}_2}^3+{{\yyze}_1}$ ${{\yyze}_2}^2)$ ${\Omega}^6-{{\yyze}_2}$ ${{\yyze}_1}^2$ $({{\yyze}_1}$ ${{\yyze}_2}+
2$ ${{\yyze}_1}^2+3$ ${{\yyze}_2}^2)$ ${\Omega}^4\big)$ ${\evec_1}$ $+
\big(-{{\yyze}_2}$ $(-10$ ${{\yyze}_1}$ ${{\yyze}_2}+2$ ${{\yyze}_1}^2-
{{\yyze}_2}^2)$ ${\Omega}^6-{{\yyze}_1}$ ${{\yyze}_2}^2$ $({{\yyze}_1}$ ${{\yyze}_2}+2$ ${{\yyze}_1}^2+
3$ ${{\yyze}_2}^2)$ ${\Omega}^4\big)$ ${\evec_2}$ $+({\Omega}^8$ ${{\yyze}_2}-{{\yyze}_2}^2$ $({{\yyze}_2}+
{{\yyze}_1})$ ${\Omega}^6)$ ${\evec_3}$ $-4$ ${\Omega}^8$ ${{\yyze}_1}$ ${{\yyze}_2}$ $\fracp{\evec_1}{x_1}$ $+
4$ ${\Omega}^8$ ${{\yyze}_1}^2$ $\fracp{\evec_2}{x_1}$ $+({{\uuvec^0}_1}^2+{{\yyze}_2}$ ${{\uuvec^0}_2}-
2$ ${{\uuvec^0}_1}$ ${{\yyze}_2})$ ${\Omega}^8+\big({{\yyze}_2}^3$ ${{\uuvec^0}_1}-2$ ${{\uuvec^0}_2}^2$ ${{\yyze}_1}^2-
2$ ${{\yyze}_1}^2$ ${{\uuvec^0}_1}^2-{{\uuvec^0}_1}$ ${{\yyze}_2}$ ${{\uuvec^0}_2}$ ${{\yyze}_1}+
{{\uuvec^0}_1}^2$ ${{\yyze}_2}^2-{{\yyze}_1}^2$ ${{\yyze}_2}$ ${{\uuvec^0}_2}-2$ ${{\yyze}_2}^3$ ${{\uuvec^0}_2}-
2$ ${{\uuvec^0}_2}$ ${{\yyze}_1}^3+3$ ${{\yyze}_1}$ ${{\yyze}_2}$ ${{\uuvec^0}_2}^2+
15$ ${{\yyze}_1}$ ${{\yyze}_2}$ ${{\uuvec^0}_1}^2+3$ ${{\uuvec^0}_2}$ ${{\uuvec^0}_1}$ ${{\yyze}_2}^2+
4$ ${{\yyze}_1}^2$ ${{\yyze}_2}$ ${{\uuvec^0}_1}+2$ ${{\uuvec^0}_3}^2$ ${{\yyze}_1}$ ${{\yyze}_2}\big)$ ${\Omega}^6+
\big({{\yyze}_1}$ ${{\yyze}_2}^3$ ${{\uuvec^0}_2}^2-10$ ${{\yyze}_2}^3$ ${{\uuvec^0}_3}^2$ ${{\yyze}_1}-
7$ ${{\yyze}_1}$ ${{\yyze}_2}^3$ ${{\uuvec^0}_1}$ ${{\uuvec^0}_2}+6$ ${{\yyze}_1}^2$ ${{\yyze}_2}^2$ ${{\uuvec^0}_2}^2-
{{\uuvec^0}_1}$ ${{\yyze}_2}^4$ ${{\uuvec^0}_2}+11$ ${{\uuvec^0}_2}^2$ ${{\yyze}_1}^3$ ${{\yyze}_2}-
8$ ${{\yyze}_1}^4$ ${{\uuvec^0}_2}$ ${{\uuvec^0}_1}-6$ ${{\yyze}_2}^4$ ${{\uuvec^0}_3}^2-
19$ ${{\yyze}_1}^2$ ${{\uuvec^0}_1}$ ${{\yyze}_2}^2$ ${{\uuvec^0}_2}+2$ ${{\uuvec^0}_3}^2$ ${{\yyze}_1}^4+
{{\yyze}_2}^2$ ${{\uuvec^0}_3}$ ${{\uuvec^0}_1}$ ${{\yyze}_1}^2+
2$ ${{\yyze}_1}^3$ ${{\uuvec^0}_1}$ ${{\yyze}_2}$ ${{\uuvec^0}_2}+{{\yyze}_1}^2$ ${{\yyze}_2}^2$ ${{\uuvec^0}_3}$ ${{\uuvec^0}_2}-
11$ ${{\yyze}_2}^2$ ${{\uuvec^0}_3}^2$ ${{\yyze}_1}^2-16$ ${{\yyze}_2}$ ${{\uuvec^0}_3}^2$ ${{\yyze}_1}^3\big)$ ${\Omega}^4+
\big(-3$ ${{\yyze}_1}^4$ ${{\uuvec^0}_1}^2$ ${{\yyze}_2}^2+3$ ${{\yyze}_1}^4$ ${{\yyze}_2}^2$ ${{\uuvec^0}_2}^2-
10$ ${{\yyze}_1}^3$ ${{\yyze}_2}^3$ ${{\uuvec^0}_1}^2-16$ ${{\yyze}_1}^3$ ${{\yyze}_2}^3$ ${{\uuvec^0}_2}^2+
6$ ${{\yyze}_1}^4$ ${{\uuvec^0}_1}$ ${{\yyze}_2}^2$ ${{\uuvec^0}_2}+
5$ ${{\yyze}_1}^3$ ${{\uuvec^0}_1}$ ${{\yyze}_2}^3$ ${{\uuvec^0}_2}+
11$ ${{\yyze}_1}^2$ ${{\uuvec^0}_1}$ ${{\yyze}_2}^4$ ${{\uuvec^0}_2}+
8$ ${{\yyze}_1}^5$ ${{\yyze}_2}$ ${{\uuvec^0}_2}$ ${{\uuvec^0}_1}-
{{\yyze}_1}^3$ ${{\yyze}_2}^3$ ${{\uuvec^0}_1}$ ${{\uuvec^0}_3}-{{\yyze}_1}^2$ ${{\yyze}_2}^4$ ${{\uuvec^0}_2}$ ${{\uuvec^0}_3}-
{{\yyze}_1}^4$ ${{\yyze}_2}^2$ ${{\uuvec^0}_1}$ ${{\uuvec^0}_3}-{{\yyze}_1}^3$ ${{\uuvec^0}_2}$ ${{\yyze}_2}^3$ ${{\uuvec^0}_3}-
12$ ${{\yyze}_1}^5$ ${{\yyze}_2}$ ${{\uuvec^0}_1}^2+21$ ${{\yyze}_1}^3$ ${{\uuvec^0}_3}^2$ ${{\yyze}_2}^3+
7$ ${{\yyze}_1}$ ${{\uuvec^0}_3}^2$ ${{\yyze}_2}^5+8$ ${{\yyze}_1}^6$ ${{\uuvec^0}_1}$ ${{\uuvec^0}_2}-
2$ ${{\yyze}_1}$ ${{\uuvec^0}_1}^2$ ${{\yyze}_2}^5-5$ ${{\yyze}_1}$ ${{\yyze}_2}^5$ ${{\uuvec^0}_2}^2-
3$ ${{\yyze}_1}^2$ ${{\uuvec^0}_1}^2$ ${{\yyze}_2}^4-6$ ${{\yyze}_1}^5$ ${{\yyze}_2}$ ${{\uuvec^0}_2}^2+
4$ ${{\uuvec^0}_1}$ ${{\uuvec^0}_2}$ ${{\yyze}_2}^6+18$ ${{\yyze}_1}^2$ ${{\yyze}_2}^4$ ${{\uuvec^0}_3}^2+
13$ ${{\yyze}_1}^4$ ${{\uuvec^0}_3}^2$ ${{\yyze}_2}^2+13$ ${{\yyze}_1}^5$ ${{\yyze}_2}$ ${{\uuvec^0}_3}^2+
4$ ${{\yyze}_1}^2$ ${{\yyze}_2}^4$ ${{\uuvec^0}_2}^2-4$ ${{\yyze}_1}$ ${{\uuvec^0}_1}$ ${{\yyze}_2}^5$ ${{\uuvec^0}_2}+
3$ ${{\uuvec^0}_1}^2$ ${{\yyze}_2}^6+6$ ${{\yyze}_2}^6$ ${{\uuvec^0}_3}^2+3$ ${{\uuvec^0}_1}^2$ ${{\yyze}_1}^6\big)$ ${\Omega}^2+
2$ ${{\yyze}_1}^4$ ${{\yyze}_2}$ $({{\yyze}_2}+{{\yyze}_1})$ $(-{{\uuvec^0}_1}$ ${{\yyze}_2}+
{{\uuvec^0}_2}$ ${{\yyze}_1})$ $({{\yyze}_1}$ ${{\uuvec^0}_1}+{{\yyze}_2}$ ${{\uuvec^0}_2})\Big)$ ${{\yyun}_1}+
\Big(-4$ ${\Omega}^7$ ${{\yyze}_1}^2$ ${{\uuvec^0}_3}+{{\yyze}_2}$ $(-8$ ${{\yyze}_1}^2$ ${{\yyze}_2}$ ${{\uuvec^0}_3}+
{{\yyze}_1}$ ${{\yyze}_2}^2$ ${{\uuvec^0}_2}-3$ ${{\yyze}_1}$ ${{\uuvec^0}_3}$ ${{\yyze}_2}^2+3$ ${{\yyze}_1}^3$ ${{\uuvec^0}_3}-
3$ ${{\yyze}_2}^3$ ${{\uuvec^0}_3}+{{\yyze}_1}^2$ ${{\yyze}_2}$ ${{\uuvec^0}_1})$ ${\Omega}^5+
{{\yyze}_2}^2$ $(3$ ${{\yyze}_1}^3$ ${{\uuvec^0}_3}$ ${{\yyze}_2}-{{\yyze}_1}^2$ ${{\yyze}_2}^2$ ${{\uuvec^0}_2}+
9$ ${{\yyze}_2}^2$ ${{\yyze}_1}^2$ ${{\uuvec^0}_3}-{{\uuvec^0}_2}$ ${{\yyze}_2}^3$ ${{\yyze}_1}+5$ ${{\yyze}_1}^4$ ${{\uuvec^0}_3}-
{{\yyze}_1}^2$ ${{\uuvec^0}_1}$ ${{\yyze}_2}^2-{{\uuvec^0}_1}$ ${{\yyze}_2}$ ${{\yyze}_1}^3+3$ ${{\uuvec^0}_3}$ ${{\yyze}_2}^4+
2$ ${{\yyze}_1}$ ${{\yyze}_2}^3$ ${{\uuvec^0}_3})$ ${\Omega}^3\Big)$ ${\evec_2}$ $+
6$ ${\Omega}^7$ ${{\yyze}_1}$ ${{\yyze}_2}^2$ $\fracp{\evec_1}{x_2}$ ${{\uuvec^0}_3}+
6$ ${\Omega}^7$ ${{\yyze}_1}^2$ ${{\yyze}_2}$ $\fracp{\evec_1}{x_1}$ ${{\uuvec^0}_3}-
6$ ${{\yyze}_1}$ ${{\yyze}_2}$ $({{\yyze}_1}$ ${{\uuvec^0}_1}+{{\yyze}_2}$ ${{\uuvec^0}_2})$ ${\Omega}^7$ $\fracp{\evec_1}{x_3}$ $-
2$ ${{\yyze}_1}$ $(-{{\yyze}_2}^2+2$ ${{\yyze}_1}^2)$ ${{\uuvec^0}_3}$ ${\Omega}^7$ $\fracp{\evec_2}{x_1}$ $-
2$ ${{\yyze}_2}$ $(-{{\yyze}_2}^2+2$ ${{\yyze}_1}^2)$ ${{\uuvec^0}_3}$ ${\Omega}^7$ $\fracp{\evec_2}{x_2}$ $\bigg)$ $
\bigg/$ ${\Omega}^{10},
$

~

$
\frac{d{\uuvec^1}_3}{dt}$ $=$ $
$
$
1/4$ $\bigg(4$ ${\Omega}^7$ ${{\yyze}_2}$ $\fracp{\evec_2}{t}$ $+4$ ${\Omega}^7$ $\fracp{\evec_1}{t}$ ${{\yyze}_1}+
2$ ${{\yyze}_1}$ $\Big(-3$ ${{\uuvec^0}_1}$ ${{\yyze}_2}$ ${{\yyze}_1}-{{\yyze}_2}^2$ ${{\uuvec^0}_2}+
2$ ${{\yyze}_1}^2$ ${{\uuvec^0}_2}\Big)$ ${\Omega}^5$ $\fracp{\evec_1}{x_2}$ $+
2$ ${\Omega}^7$ ${{\uuvec^0}_3}$ ${{\yyze}_2}$ $\fracp{\evec_3}{x_2}$ $-
2$ ${\Omega}^7$ ${{\uuvec^0}_3}$ ${{\yyze}_1}$ $\fracp{\evec_1}{x_3}$ $+
2$ ${\Omega}^7$ ${{\uuvec^0}_3}$ $\fracp{\evec_3}{x_1}$ ${{\yyze}_1}-
2$ ${\Omega}^7$ ${{\uuvec^0}_3}$ ${{\yyze}_2}$ $\fracp{\evec_2}{x_3}$ $-
2$ ${{\yyze}_2}$ $\Big(-2$ ${{\uuvec^0}_1}$ ${{\yyze}_2}^2+{{\yyze}_1}^2$ ${{\uuvec^0}_1}+
3$ ${{\yyze}_2}$ ${{\uuvec^0}_2}$ ${{\yyze}_1}\Big)$ ${\Omega}^5$ $\fracp{\evec_2}{x_1}$ $+
2$ ${{\yyze}_2}$ $\Big(-3$ ${{\uuvec^0}_1}$ ${{\yyze}_2}$ ${{\yyze}_1}-{{\yyze}_2}^2$ ${{\uuvec^0}_2}+
2$ ${{\yyze}_1}^2$ ${{\uuvec^0}_2}\Big)$ ${\Omega}^5$ $\fracp{\evec_2}{x_2}$ $-
2$ ${{\yyze}_1}$ $\Big(-2$ ${{\uuvec^0}_1}$ ${{\yyze}_2}^2+{{\yyze}_1}^2$ ${{\uuvec^0}_1}+
3$ ${{\yyze}_2}$ ${{\uuvec^0}_2}$ ${{\yyze}_1}\Big)$ ${\Omega}^5$ $\fracp{\evec_1}{x_1}$ $+
\Big(-2$ ${{\uuvec^0}_2}$ ${{\yyze}_1}^5$ ${{\uuvec^0}_3}^2-2$ ${{\yyze}_2}^5$ ${{\uuvec^0}_1}^3-
2$ ${{\yyze}_2}^2$ ${{\yyze}_1}^3$ ${{\uuvec^0}_2}^3+2$ ${{\yyze}_2}$ ${{\uuvec^0}_2}^3$ ${{\yyze}_1}^4+
8$ ${{\yyze}_2}^4$ ${{\uuvec^0}_1}^3$ ${{\yyze}_1}+{{\yyze}_1}$ ${{\uuvec^0}_3}$ ${{\yyze}_2}^4$ ${{\uuvec^0}_2}^2+
{{\yyze}_1}^3$ ${{\uuvec^0}_3}$ ${{\yyze}_2}^2$ ${{\uuvec^0}_1}^2+
{{\yyze}_1}^4$ ${{\uuvec^0}_3}$ ${{\yyze}_2}$ ${{\uuvec^0}_1}^2+
{{\yyze}_1}^2$ ${{\uuvec^0}_3}$ ${{\yyze}_2}^3$ ${{\uuvec^0}_2}^2-
8$ ${{\yyze}_1}^2$ ${{\uuvec^0}_1}^2$ ${{\yyze}_2}^3$ ${{\uuvec^0}_2}+
2$ ${{\uuvec^0}_1}$ ${{\yyze}_1}^3$ ${{\yyze}_2}^2$ ${{\uuvec^0}_2}^2+
4$ ${{\yyze}_1}^3$ ${{\yyze}_2}^2$ ${{\uuvec^0}_2}$ ${{\uuvec^0}_1}^2-
{{\yyze}_1}^3$ ${{\yyze}_2}^2$ ${{\uuvec^0}_3}$ ${{\uuvec^0}_2}^2-
{{\uuvec^0}_1}^2$ ${{\yyze}_1}^2$ ${{\yyze}_2}^3$ ${{\uuvec^0}_3}-
{{\uuvec^0}_1}^2$ ${{\yyze}_2}^4$ ${{\uuvec^0}_3}$ ${{\yyze}_1}+2$ ${{\yyze}_1}^2$ ${{\uuvec^0}_2}^2$ ${{\uuvec^0}_1}$ ${{\yyze}_2}^3+
2$ ${{\yyze}_1}$ ${{\yyze}_2}^4$ ${{\uuvec^0}_2}$ ${{\uuvec^0}_1}^2-2$ ${{\yyze}_2}$ ${{\uuvec^0}_1}$ ${{\yyze}_1}^4$ ${{\uuvec^0}_2}^2+
4$ ${{\uuvec^0}_1}$ ${{\uuvec^0}_2}$ ${{\uuvec^0}_3}$ ${{\yyze}_2}^2$ ${{\yyze}_1}^3+
4$ ${{\yyze}_1}^2$ ${{\yyze}_2}^3$ ${{\uuvec^0}_2}^3+
4$ ${{\yyze}_2}^2$ ${{\yyze}_1}^3$ ${{\uuvec^0}_1}^3-{{\yyze}_1}^4$ ${{\yyze}_2}$ ${{\uuvec^0}_3}$ ${{\uuvec^0}_2}^2-
3$ ${{\uuvec^0}_2}$ ${{\uuvec^0}_3}^2$ ${{\yyze}_1}^3$ ${{\yyze}_2}^2-5$ ${{\uuvec^0}_2}$ ${{\yyze}_2}$ ${{\uuvec^0}_3}^2$ ${{\yyze}_1}^4+
2$ ${{\uuvec^0}_1}$ ${{\yyze}_2}^2$ ${{\yyze}_1}^3$ ${{\uuvec^0}_3}^2+2$ ${{\uuvec^0}_1}$ ${{\yyze}_2}$ ${{\yyze}_1}^4$ ${{\uuvec^0}_3}^2-
3$ ${{\uuvec^0}_2}$ ${{\yyze}_2}^4$ ${{\uuvec^0}_3}^2$ ${{\yyze}_1}-3$ ${{\uuvec^0}_2}$ ${{\uuvec^0}_3}^2$ ${{\yyze}_2}^3$ ${{\yyze}_1}^2-
6$ ${{\yyze}_2}^4$ ${{\uuvec^0}_1}$ ${{\yyze}_1}$ ${{\uuvec^0}_2}^2+
4$ ${{\uuvec^0}_1}$ ${{\uuvec^0}_2}$ ${{\yyze}_2}^3$ ${{\uuvec^0}_3}$ ${{\yyze}_1}^2-2$ ${{\yyze}_2}^3$ ${{\uuvec^0}_1}^3$ ${{\yyze}_1}^2+
2$ ${{\yyze}_2}^5$ ${{\uuvec^0}_2}$ ${{\uuvec^0}_1}^2\Big)$ ${\Omega}+\Big(4$ ${{\yyze}_1}^2$ ${{\uuvec^0}_2}$ ${{\uuvec^0}_1}-
4$ ${{\yyze}_1}$ ${{\yyze}_2}$ ${{\uuvec^0}_1}^2+2$ ${{\yyze}_1}^2$ ${{\uuvec^0}_3}^2+2$ ${{\yyze}_1}^2$ ${{\uuvec^0}_3}$ ${{\uuvec^0}_2}-
2$ ${{\yyze}_1}$ ${{\uuvec^0}_3}$ ${{\uuvec^0}_1}$ ${{\yyze}_2}+
{{\yyze}_1}$ ${{\uuvec^0}_1}^2$ ${{\uuvec^0}_3}+{{\uuvec^0}_3}$ ${{\yyze}_1}$ ${{\uuvec^0}_2}^2+
4$ ${{\yyze}_1}$ ${{\yyze}_2}$ ${{\uuvec^0}_2}^2-4$ ${{\uuvec^0}_2}$ ${{\uuvec^0}_1}$ ${{\yyze}_2}^2\Big)$ ${\Omega}^5+
\Big({{\uuvec^0}_3}$ $(-3$ ${{\uuvec^0}_2}+{{\yyze}_1})$ ${\Omega}^6-{{\yyze}_1}$ ${{\uuvec^0}_3}$ $(2$ ${{\yyze}_1}$ ${{\uuvec^0}_1}+
3$ ${{\yyze}_2}$ ${{\uuvec^0}_2})$ ${\Omega}^4+{{\uuvec^0}_3}$ $(-{{\uuvec^0}_2}$ ${{\yyze}_2}$ ${{\yyze}_1}^3+
7$ ${{\yyze}_1}^2$ ${{\uuvec^0}_1}$ ${{\yyze}_2}^2+4$ ${{\yyze}_1}^4$ ${{\uuvec^0}_1}-3$ ${{\yyze}_1}$ ${{\yyze}_2}^3$ ${{\uuvec^0}_1}+
3$ ${{\uuvec^0}_2}$ ${{\yyze}_1}^4+10$ ${{\yyze}_1}^2$ ${{\yyze}_2}^2$ ${{\uuvec^0}_2}+2$ ${{\uuvec^0}_2}$ ${{\yyze}_2}^3$ ${{\yyze}_1}+
4$ ${{\uuvec^0}_2}$ ${{\yyze}_2}^4)$ ${\Omega}^2+{{\yyze}_2}$ ${{\uuvec^0}_3}$ $\big(-{{\yyze}_2}^5$ ${{\uuvec^0}_2}-
{{\yyze}_1}^3$ ${{\yyze}_2}^2$ ${{\uuvec^0}_2}+{{\yyze}_1}^2$ ${{\uuvec^0}_1}$ ${{\yyze}_2}^3+2$ ${{\yyze}_1}^5$ ${{\uuvec^0}_2}+
2$ ${{\yyze}_1}^4$ ${{\yyze}_2}$ ${{\uuvec^0}_2}+{{\yyze}_1}$ ${{\yyze}_2}^4$ ${{\uuvec^0}_2}-
3$ ${{\yyze}_1}^2$ ${{\yyze}_2}^3$ ${{\uuvec^0}_2}-{{\yyze}_1}$ ${{\yyze}_2}^4$ ${{\uuvec^0}_1}-
3$ ${{\yyze}_1}^4$ ${{\yyze}_2}$ ${{\uuvec^0}_1}-5$ ${{\yyze}_1}^3$ ${{\yyze}_2}^2$ ${{\uuvec^0}_1}\big)\Big)$ ${{\yyun}_2}+
\Big(4$ ${\Omega}^7$ ${{\yyze}_1}+2$ ${{\yyze}_1}$ $(6$ ${{\uuvec^0}_2}$ ${{\yyze}_1}+{{\yyze}_2}$ ${{\uuvec^0}_3}-
6$ ${{\uuvec^0}_1}$ ${{\yyze}_2})$ ${\Omega}^5-2$ ${{\yyze}_2}^2$ $({{\yyze}_2}+{{\yyze}_1})$ $({{\yyze}_1}$ ${{\uuvec^0}_3}+
{{\uuvec^0}_2}$ ${{\yyze}_1}-{{\uuvec^0}_1}$ ${{\yyze}_2})$ ${\Omega}^3+2$ ${{\yyze}_2}^2$ $({{\yyze}_1}^3+2$ ${{\yyze}_2}^3+
{{\yyze}_2}$ ${{\yyze}_1}^2)$ $(-{{\uuvec^0}_1}$ ${{\yyze}_2}+{{\uuvec^0}_2}$ ${{\yyze}_1})$ ${\Omega}\Big)$ ${\evec_2}$ $+
\Big(\big((-2$ ${{\yyze}_1}$ ${{\uuvec^0}_1}+3$ ${{\yyze}_2}$ ${{\uuvec^0}_2}+3$ ${{\uuvec^0}_2}$ ${{\yyze}_1})$ ${\Omega}^5+
(-2$ ${{\yyze}_2}^3$ ${{\uuvec^0}_2}-6$ ${{\yyze}_1}$ ${{\yyze}_2}^2$ ${{\uuvec^0}_2}-5$ ${{\yyze}_1}^2$ ${{\yyze}_2}$ ${{\uuvec^0}_2}-
2$ ${{\yyze}_1}^2$ ${{\yyze}_2}$ ${{\uuvec^0}_1}+2$ ${{\yyze}_1}$ ${{\yyze}_2}^2$ ${{\uuvec^0}_1}-
3$ ${{\uuvec^0}_2}$ ${{\yyze}_1}^3)$ ${\Omega}^3+{{\yyze}_2}$ $({{\yyze}_2}$ ${{\yyze}_1}^2+{{\yyze}_1}^3-
{{\yyze}_2}^3+3$ ${{\yyze}_1}$ ${{\yyze}_2}^2)$ $({{\yyze}_1}$ ${{\uuvec^0}_1}+
{{\yyze}_2}$ ${{\uuvec^0}_2})$ ${\Omega}\big)$ ${{\yyun}_2}+
({{\yyze}_2}$ ${\Omega}^5$ ${{\yyze}_1}-{{\yyze}_1}^2$ ${{\yyze}_2}$ $({{\yyze}_2}+{{\yyze}_1})$ ${\Omega}^3)$ ${{\uuvec^1}_1}+
({\Omega}^5$ ${{\yyze}_2}^2-{{\yyze}_1}$ ${{\yyze}_2}^2$ $({{\yyze}_2}+{{\yyze}_1})$ ${\Omega}^3)$ ${{\uuvec^1}_2}+({\Omega}^7-
{{\yyze}_2}$ $({{\yyze}_2}+{{\yyze}_1})$ ${\Omega}^5)$ ${{\uuvec^1}_3}+(-{{\yyze}_1}$ ${\Omega}^6+
{{\yyze}_1}$ ${{\yyze}_2}$ $({{\yyze}_2}+{{\yyze}_1})$ ${\Omega}^4)$ ${\evec_1}$ $+(-{{\yyze}_2}$ ${\Omega}^6+
{{\yyze}_2}^2$ $({{\yyze}_2}+{{\yyze}_1})$ ${\Omega}^4)$ ${\evec_2}$ $+({{\yyze}_2}$ ${\Omega}^6-
{{\yyze}_1}$ ${{\yyze}_2}$ $({{\yyze}_2}+{{\yyze}_1})$ ${\Omega}^4)$ ${\evec_3}$ $+\big({{\uuvec^0}_1}$ ${{\yyze}_2}-
4$ ${{\uuvec^0}_1}$ ${{\uuvec^0}_3}-{{\yyze}_2}$ ${{\uuvec^0}_3}-{{\uuvec^0}_2}$ ${{\yyze}_1}-
{{\yyze}_1}$ ${{\uuvec^0}_3}\big)$ ${\Omega}^6+\big(3$ ${{\uuvec^0}_3}$ ${{\uuvec^0}_2}$ ${{\yyze}_2}^2-
{{\yyze}_1}$ ${{\yyze}_2}^2$ ${{\uuvec^0}_1}+{{\yyze}_1}^2$ ${{\yyze}_2}$ ${{\uuvec^0}_2}+
2$ ${{\uuvec^0}_1}$ ${{\yyze}_2}$ ${{\uuvec^0}_2}$ ${{\yyze}_1}+{{\yyze}_1}^3$ ${{\uuvec^0}_3}-
{{\uuvec^0}_1}^2$ ${{\yyze}_2}^2-{{\uuvec^0}_2}^2$ ${{\yyze}_1}^2+{{\yyze}_1}^2$ ${{\yyze}_2}$ ${{\uuvec^0}_3}-
{{\yyze}_2}$ ${{\uuvec^0}_3}$ ${{\uuvec^0}_2}$ ${{\yyze}_1}-{{\yyze}_2}^3$ ${{\uuvec^0}_1}+
{{\yyze}_1}$ ${{\yyze}_2}^2$ ${{\uuvec^0}_2}+{{\yyze}_1}^2$ ${{\uuvec^0}_3}$ ${{\uuvec^0}_2}-
2$ ${{\yyze}_1}$ ${{\uuvec^0}_3}$ ${{\uuvec^0}_1}$ ${{\yyze}_2}+2$ ${{\yyze}_2}^2$ ${{\uuvec^0}_3}$ ${{\uuvec^0}_1}\big)$ ${\Omega}^4+
\big({{\uuvec^0}_1}^2$ ${{\yyze}_2}^4+5$ ${{\uuvec^0}_3}$ ${{\uuvec^0}_1}$ ${{\yyze}_2}$ ${{\yyze}_1}^3-
5$ ${{\uuvec^0}_3}$ ${{\uuvec^0}_2}$ ${{\yyze}_1}^4+{{\uuvec^0}_2}^2$ ${{\yyze}_1}^3$ ${{\yyze}_2}+
{{\uuvec^0}_1}^2$ ${{\yyze}_1}$ ${{\yyze}_2}^3-2$ ${{\yyze}_1}^2$ ${{\uuvec^0}_1}$ ${{\yyze}_2}^2$ ${{\uuvec^0}_2}+
8$ ${{\yyze}_2}^2$ ${{\uuvec^0}_3}$ ${{\uuvec^0}_1}$ ${{\yyze}_1}^2+{{\yyze}_1}^2$ ${{\yyze}_2}^2$ ${{\uuvec^0}_2}^2+
2$ ${{\uuvec^0}_2}$ ${{\yyze}_1}$ ${{\yyze}_2}^3$ ${{\uuvec^0}_3}+2$ ${{\uuvec^0}_2}$ ${{\yyze}_1}^3$ ${{\yyze}_2}$ ${{\uuvec^0}_3}-
2$ ${{\yyze}_1}$ ${{\yyze}_2}^3$ ${{\uuvec^0}_1}$ ${{\uuvec^0}_2}-2$ ${{\yyze}_1}^2$ ${{\yyze}_2}^2$ ${{\uuvec^0}_3}$ ${{\uuvec^0}_2}-
{{\uuvec^0}_1}$ ${{\yyze}_1}$ ${{\yyze}_2}^3$ ${{\uuvec^0}_3}-3$ ${{\uuvec^0}_2}$ ${{\uuvec^0}_3}$ ${{\yyze}_2}^4+
4$ ${{\yyze}_2}^4$ ${{\uuvec^0}_3}$ ${{\uuvec^0}_1}+4$ ${{\uuvec^0}_3}$ ${{\uuvec^0}_1}$ ${{\yyze}_1}^4\big)$ ${\Omega}^2+
{{\yyze}_1}$ ${{\uuvec^0}_3}$ $\big(-{{\yyze}_2}^5$ ${{\uuvec^0}_2}-{{\yyze}_1}^3$ ${{\yyze}_2}^2$ ${{\uuvec^0}_2}+
{{\yyze}_1}^2$ ${{\uuvec^0}_1}$ ${{\yyze}_2}^3+2$ ${{\yyze}_1}^5$ ${{\uuvec^0}_2}+2$ ${{\yyze}_1}^4$ ${{\yyze}_2}$ ${{\uuvec^0}_2}+
{{\yyze}_1}$ ${{\yyze}_2}^4$ ${{\uuvec^0}_2}-3$ ${{\yyze}_1}^2$ ${{\yyze}_2}^3$ ${{\uuvec^0}_2}-
{{\yyze}_1}$ ${{\yyze}_2}^4$ ${{\uuvec^0}_1}-3$ ${{\yyze}_1}^4$ ${{\yyze}_2}$ ${{\uuvec^0}_1}-
5$ ${{\yyze}_1}^3$ ${{\yyze}_2}^2$ ${{\uuvec^0}_1}\big)\Big)$ ${{\yyun}_1}+\Big(-2$ ${{\uuvec^0}_3}^2$ ${{\yyze}_1}^4+
2$ ${{\uuvec^0}_1}^2$ ${{\yyze}_2}^4-2$ ${{\yyze}_2}$ ${{\uuvec^0}_3}^2$ ${{\yyze}_1}^3-
2$ ${{\uuvec^0}_1}^2$ ${{\yyze}_1}$ ${{\yyze}_2}^3-2$ ${{\yyze}_1}^2$ ${{\yyze}_2}^2$ ${{\uuvec^0}_3}$ ${{\uuvec^0}_2}-
2$ ${{\yyze}_1}$ ${{\yyze}_2}^3$ ${{\uuvec^0}_1}$ ${{\uuvec^0}_2}-2$ ${{\uuvec^0}_2}$ ${{\yyze}_1}^3$ ${{\yyze}_2}$ ${{\uuvec^0}_3}-
4$ ${{\yyze}_1}^2$ ${{\uuvec^0}_1}^2$ ${{\yyze}_2}$ ${{\uuvec^0}_2}+2$ ${{\yyze}_2}^2$ ${{\uuvec^0}_3}$ ${{\uuvec^0}_1}$ ${{\yyze}_1}^2+
2$ ${{\yyze}_1}^2$ ${{\uuvec^0}_1}$ ${{\yyze}_2}^2$ ${{\uuvec^0}_2}+2$ ${{\uuvec^0}_1}$ ${{\yyze}_1}$ ${{\yyze}_2}^3$ ${{\uuvec^0}_3}-
{{\uuvec^0}_3}$ ${{\yyze}_1}^2$ ${{\yyze}_2}$ ${{\uuvec^0}_2}^2+2$ ${{\yyze}_1}^3$ ${{\uuvec^0}_1}$ ${{\uuvec^0}_2}^2+
{{\uuvec^0}_3}$ ${{\yyze}_1}^3$ ${{\uuvec^0}_2}^2-{{\yyze}_1}^3$ ${{\uuvec^0}_1}^2$ ${{\uuvec^0}_3}-
4$ ${{\yyze}_2}$ ${{\uuvec^0}_3}$ ${{\uuvec^0}_2}$ ${{\yyze}_1}^2$ ${{\uuvec^0}_1}-
2$ ${{\yyze}_1}$ ${{\yyze}_2}^2$ ${{\uuvec^0}_3}$ ${{\uuvec^0}_2}^2-{{\yyze}_1}^2$ ${{\uuvec^0}_1}^2$ ${{\uuvec^0}_3}$ ${{\yyze}_2}-
2$ ${{\yyze}_2}^2$ ${{\uuvec^0}_2}^2$ ${{\uuvec^0}_1}$ ${{\yyze}_1}+{{\yyze}_2}$ ${{\uuvec^0}_2}$ ${{\yyze}_1}^2$ ${{\uuvec^0}_3}^2+
2$ ${{\yyze}_2}^2$ ${{\uuvec^0}_1}$ ${{\uuvec^0}_3}^2$ ${{\yyze}_1}+3$ ${{\yyze}_2}^2$ ${{\uuvec^0}_2}$ ${{\uuvec^0}_3}^2$ ${{\yyze}_1}-
4$ ${{\yyze}_2}$ ${{\uuvec^0}_1}$ ${{\uuvec^0}_3}^2$ ${{\yyze}_1}^2+2$ ${{\yyze}_2}^3$ ${{\uuvec^0}_2}$ ${{\uuvec^0}_1}^2-
2$ ${{\yyze}_1}$ ${{\uuvec^0}_1}^3$ ${{\yyze}_2}^2+4$ ${{\uuvec^0}_2}$ ${{\uuvec^0}_3}^2$ ${{\yyze}_1}^3\Big)$ ${\Omega}^3+
\Big(-{{\yyze}_2}$ $({{\yyze}_2}+{{\yyze}_1})$ ${{\uuvec^0}_3}$ ${\Omega}^5+
{{\yyze}_2}$ $(-2$ ${{\yyze}_1}^2$ ${{\yyze}_2}$ ${{\uuvec^0}_2}+{{\yyze}_1}^3$ ${{\uuvec^0}_3}-
2$ ${{\yyze}_1}^2$ ${{\yyze}_2}$ ${{\uuvec^0}_1}+2$ ${{\yyze}_1}$ ${{\yyze}_2}^2$ ${{\uuvec^0}_1}+
3$ ${{\yyze}_1}^2$ ${{\yyze}_2}$ ${{\uuvec^0}_3}-{{\yyze}_1}$ ${{\uuvec^0}_3}$ ${{\yyze}_2}^2+
2$ ${{\uuvec^0}_2}$ ${{\yyze}_1}^3+{{\yyze}_2}^3$ ${{\uuvec^0}_3})$ ${\Omega}^3\Big)$ ${\evec_3}$ $+
\Big(({{\yyze}_1}$ ${{\uuvec^0}_1}+
{{\yyze}_2}$ ${{\uuvec^0}_2})$ ${\Omega}^5-{{\yyze}_2}^2$ $({{\yyze}_1}$ ${{\uuvec^0}_1}+
{{\yyze}_2}$ ${{\uuvec^0}_2})$ ${\Omega}^3\Big)$ ${{\yyun}_2}^2+\Big(-{\Omega}^6$ ${{\yyze}_1}$ ${{\uuvec^0}_3}+
{{\yyze}_1}$ ${{\yyze}_2}$ $({{\yyze}_2}+{{\yyze}_1})$ ${{\uuvec^0}_3}$ ${\Omega}^4\Big)$ ${{\uuvec^1}_3}+
(-2$ ${{\yyze}_1}$ ${{\uuvec^0}_1}-2$ ${{\yyze}_2}$ ${{\uuvec^0}_2})$ ${\Omega}^7$ $\fracp{\evec_3}{x_3}$ $+
\Big({{\yyze}_1}^2$ $(-{{\yyze}_2}+{{\yyze}_1})$ ${{\uuvec^0}_3}$ ${\Omega}^4-
{{\yyze}_1}^2$ $({{\yyze}_2}+{{\yyze}_1})$ $(-{{\yyze}_2}+{{\yyze}_1})^2$ ${{\uuvec^0}_3}$ ${\Omega}^2\Big)$ ${{\uuvec^1}_1}+
\Big(-4$ ${\Omega}^7$ ${{\yyze}_2}+(2$ ${{\yyze}_1}^2$ ${{\uuvec^0}_3}+12$ ${{\uuvec^0}_1}$ ${{\yyze}_2}^2-
12$ ${{\yyze}_2}$ ${{\uuvec^0}_2}$ ${{\yyze}_1})$ ${\Omega}^5-2$ ${{\yyze}_1}$ ${{\yyze}_2}$ $({{\yyze}_2}+
{{\yyze}_1})$ $({{\yyze}_1}$ ${{\uuvec^0}_3}+{{\uuvec^0}_2}$ ${{\yyze}_1}-{{\uuvec^0}_1}$ ${{\yyze}_2})$ ${\Omega}^3+
2$ ${{\yyze}_1}$ ${{\yyze}_2}$ $({{\yyze}_1}^3+2$ ${{\yyze}_2}^3+{{\yyze}_2}$ ${{\yyze}_1}^2)$ $(-{{\uuvec^0}_1}$ ${{\yyze}_2}+
{{\uuvec^0}_2}$ ${{\yyze}_1})$ ${\Omega}\Big)$ ${\evec_1}$ $+
(2$ ${{\yyze}_2}$ ${{\uuvec^0}_2}+2$ ${{\yyze}_1}$ ${{\uuvec^0}_1})$ ${\Omega}^7+
\Big({{\yyze}_1}^2$ $({{\yyze}_2}+{{\yyze}_1})$ ${{\uuvec^0}_3}$ ${\Omega}^4-
{{\yyze}_1}$ $({{\yyze}_2}+{{\yyze}_1})$ $({{\yyze}_2}^3-
{{\yyze}_1}$ ${{\yyze}_2}^2+{{\yyze}_2}$ ${{\yyze}_1}^2+{{\yyze}_1}^3)$ ${{\uuvec^0}_3}$ ${\Omega}^2\Big)$ ${{\uuvec^1}_2}+
\Big(({{\yyze}_2}$ ${{\uuvec^0}_2}+4$ ${{\yyze}_1}$ ${{\uuvec^0}_1}+2$ ${{\uuvec^0}_1}$ ${{\yyze}_2})$ ${\Omega}^5+
(-4$ ${{\yyze}_1}^3$ ${{\uuvec^0}_1}-{{\yyze}_2}^3$ ${{\uuvec^0}_1}-2$ ${{\yyze}_2}^3$ ${{\uuvec^0}_2}-
6$ ${{\yyze}_1}$ ${{\yyze}_2}^2$ ${{\uuvec^0}_1}+{{\yyze}_1}$ ${{\yyze}_2}^2$ ${{\uuvec^0}_2}-
{{\yyze}_1}^2$ ${{\yyze}_2}$ ${{\uuvec^0}_1}-{{\yyze}_1}^2$ ${{\yyze}_2}$ ${{\uuvec^0}_2})$ ${\Omega}^3+
{{\yyze}_1}$ ${{\yyze}_2}$ $({{\yyze}_1}^2+2$ ${{\yyze}_1}$ ${{\yyze}_2}-{{\yyze}_2}^2)$ $({{\yyze}_1}$ ${{\uuvec^0}_1}+
{{\yyze}_2}$ ${{\uuvec^0}_2})$ ${\Omega}\Big)$ ${{\yyun}_1}^2\bigg)$ $\bigg/$ ${\Omega}^8,
$

~

$\uuvec^1(s;\xvec,\vvec,s)=$ $0.$
}

\end{appendix}
{\small
\bibliographystyle{plain}
\bibliography{biblio}
}

\end{document}